\documentclass[11pt]{amsart}
\usepackage{enumerate}
\usepackage{amsopn,amsfonts,amsmath,amssymb}
\usepackage[dvips]{graphicx}
\usepackage{float}
\usepackage{subfigure}
\usepackage{dcolumn}
\usepackage{color, colortbl}
\usepackage{booktabs}
\usepackage[usenames,dvipsnames,svgnames,table]{xcolor}
\usepackage[colorlinks=true,
            linkcolor=red,
            urlcolor=blue,
            citecolor=gray]{hyperref}
\allowdisplaybreaks
\usepackage{geometry}  
\def\E{{\mathbb E} }	
\def\indic{{\rm {\large 1}\hspace{-2.3pt}{\large
l}}}

\def\R{{\mathbb R}}
\def\Z{{\mathbb Z}}
\def\N{{\mathbb N}}

\def\C{{\mathbb C}}
\def\supp{{\rm supp}}

\def\argmin{{\rm argmin}}

\newcommand{\norm}[1]{\left\|#1 \right\|}
\newcommand{\braket}[1]{\left\langle #1 \right\rangle }
\newcommand{\abs}[1]{ \left|  #1\right| }
\newcommand{\mt}[1]{ \boldsymbol{ #1 } }

\theoremstyle{definition}
\newtheorem{theorem}{Theorem}
\newtheorem{proposition}{Proposition}
\newtheorem{lemma}{Lemma}

\newtheorem{assumption}{Assumption}

\newtheorem{remark}{Remark}

\usepackage[foot]{amsaddr}

\begin{document}
\title[]{Adaptive estimation in the linear random coefficients model when regressors have limited variation}

\author[Gaillac]{Christophe Gaillac$^{{(1),(2)}}$}
\address{$ ^{(1)}$ Toulouse School of Economics, 1 esplanade de l'universit\'e, 31000 Toulouse, France}
\address{$ ^{(2)}$ CREST, 5 avenue Henry Le Chatelier, 91764 Palaiseau, France}
\email{\href{mailto:christophe.gaillac@tse-fr.eu}
{christophe.gaillac@tse-fr.eu}}

\author[Gautier]{Eric Gautier$^{(1)}$}
\email{\href{mailto:eric.gautier@tse-fr.eu}{eric.gautier@tse-fr.eu}}
\date{This version: \today.} 

\thanks{\emph{Keywords}: Adaptation, 
Ill-posed Inverse Problem, Minimax, 
Random Coefficients.}
\thanks{\emph{AMS 2010 Subject Classification}: Primary 62P20 ; secondary
42A99, 62C20, 62G07, 62G08, 62G20.}
\thanks{The authors acknowledge financial support from the grants ERC POEMH 337665 and ANR-17-EURE-0010. They are grateful to the seminar participants at Berkeley, Brown, CREST, Duke, Harvard-MIT, Rice, TSE, ULB, University of Tokyo, those of 2016 SFDS, ISNPS, Recent Advances in Econometrics, and 2017 IAAE conferences for comments.}

\begin{abstract}
We consider a linear model where the coefficients - intercept and slopes - are random with a law in a nonparametric class and independent from the regressors. Identification often requires the regressors to have a support which is the whole space. This is hardly ever the case in practice. Alternatively, the coefficients can have a compact support but this is not compatible with unbounded error terms as usual in regression models. In this paper, the regressors can have a support which is a proper subset but the slopes (not the intercept) do not have heavy-tails. Lower bounds on the supremum risk for the estimation of the joint density of the random coefficients density are obtained for a wide range of smoothness, where some allow for polynomial and nearly parametric rates of convergence. We present a minimax optimal estimator, a data-driven rule for adaptive estimation, and made available a \textsf{R} package.
\end{abstract}
\maketitle

\section{Introduction}\label{s1}
Inferring causal effects from a data set is of great importance for applied researchers. 
This paper assumes that the explanatory variables are determined outside the model (\emph{e.g.}, a treatment is randomly assigned)
and adresses the question of the heterogeneity of the effects. The linear regression with random coefficients (\emph{i.e.}, a continuous mixture of linear regressions) allows for heterogeneous effects across observations.  For example, a researcher interested in the effect of the income of the parents on pupils' achievements might want to allow different effects for different pupils. Maintaining parametric assumptions on the mixture density is open to criticism because these assumptions can drive the results (see \cite{HS}). 
For this reason, this paper considers a nonparametric setup. Unfortunately, most of the estimation theory for this model 
has relied on assumptions on either the data or the model which are almost never satisfied. This is probably the  reason why, up to now,  
applied researchers have preferred models such as the quantile regression. However, the assumption on the linearity of the conditional quantiles at the basis of quantile regression hold if the underlying model is a linear random coefficients model where the coefficients are functions of a scalar uniform distribution but it is hard to argue for such degeneracy. 

For a random variable $\alpha$ and random vectors $\mt{X}$ and $\mt{\beta}$ of dimension $p$,
the linear random coefficients model is
\begin{align}
&Y=\alpha+\mt{\beta}^{\top}\mt{X},\label{eRC}\\
&(\alpha,\mt{\beta}^{\top})\ \text{and}\ \mt{X}\ \text{are independent}.\label{eindep}
\end{align}
The researcher has at her disposal $n$ observations $(Y_i,\mt{X}_i^{\top})_{i=1}^n$ of $(Y,\mt{X}^{\top})$ but does not observe the realizations $(\alpha_i,\mt{\beta}_i^{\top})_{i=1}^n$ of $(\alpha,\mt{\beta}^{\top})$. $\alpha$ subsumes the intercept and error term and the vector of slope coefficients $\mt{\beta}$ is heterogeneous (\emph{i.e.}, varies across $i$). 
$(\alpha,\mt{\beta}^{\top})$ correspond to multidimensional unobserved heterogeneity and $\mt{X}$ to observed heterogeneity. Restricting unobserved heterogeneity to a scalar, as when only $\alpha$ is random or to justify quantile regression, can have undesirable implications such as monotonicity in the literature on policy evaluation (see \cite{Gautier2}). Model \eqref{eRC} is a linear model with homogeneous slopes and heteroscedastic errors, hence the averages of the coefficients are easy to obtain. 
However, the law of coefficients, 
prediction intervals for $Y$ given $\mt{X}=\mt{x}$ (see \cite{Beran3}), welfare measures, treatment and counterfactual effects, which depend on the law of the coefficients can also be of great interest. Other   random coefficients models have been analyzed recently in econometrics (see, \emph{e.g.}, \cite{breunig2018specification,GK,hoderlein2014triangular,masten2017random} and references therein). 

Estimation of the density of random coefficients $f_{\alpha,\mt{\beta}}$ has similarities with tomography problems  
involving the Radon transform (see \cite{Beran2,BM,HKM}). Indeed, the density of $Y/\sqrt{1+|\mt{X}|_2^2}$ given $\mt{S}=(1, \mt{X}^{\top})^{\top}/\sqrt{1+|\mt{X}|_2^2}$, where $|\cdot|_2$ is the Euclidian norm, at point $u$ given $\mt{s}$ is the 
integral of $f_{\alpha,\mt{\beta}}$ on the affine hyperplane defined via the pair $(u,\mt{s})$. But the random coefficients model \eqref{eRC}-\eqref{eindep} is not an inverse problem over functions or sequences with an additive Gaussian white noise. Treating it requires to allow the dimension to be larger than in tomography due to more than one or two regressors and the directions to have an unknown but estimable density. $(\alpha,\mt{\beta}^{\top})$ should also have a noncompact support  to 
allow for usual unbounded  errors. 

 To obtain rates of convergence, \cite{HKM} assumes the density of $\mt{S}$ is bounded from below. When $p=1$, this holds when $\mt{X}$ has tails at least as fat as the Cauchy distribution. Recently, \cite{dunker2017tests} motivates testing large features of the density by the possible slow rates of convergence of density estimation and 
\cite{holzmann2019rate} obtains rates of convergence for density estimation with less heavy tails on $\mt{X}$. 
But assuming the support of $\mt{X}$ is $\R^p$ 
is unrealistic for nearly all applications. 
In our motivating example, the income of the parents has limited variation. It is positive and probably bounded. 

The tomography problem corresponding to $p=1$, $\mt{S}$ has a support 
 a known cap (\emph{i.e.}, the support of the angle is an interval), and the object has support in a ball is limited angle tomography. 
\cite{frikel2013sparse} proposes a soft-thresholded curvelet 
regularization for the problem with an additive bounded noise but does not obtain results for the statistical problem (\emph{e.g.}, consistency). Importantly, \cite{Holtzman} shows 
the rate of the minimax risk in Sobolev type ellipsoids relative to the right-singular functions of the Radon transform is logarithmic and obtain that projection estimators are adaptive. It gives the analogy with a random coefficients model where $p=1$, $(\alpha,\mt{\beta}^{\top})$ have support in the unit ball, and some known densities of the regressors. It concludes that a lot remains to be done to handle $p>1$ and estimable densities of the regressors.
The random coefficients model when the support of $\mt{X}$ can be a proper (\emph{i.e.} strict) subset is considered in \cite{BM}. It assumes $p=1$ and $(\alpha,\mt{\beta}^{\top})$ have compact support, and shows that a minimum distance estimator is consistent. Section 2 in the online appendix of \cite{hoderlein2014triangular} proposes a consistent estimator in a similar situation with a single regressor for random coefficients with support in a known ball. 
  
 This paper is directly applicable to \eqref{eRC}-\eqref{eindep}.  
It allows for arbitrary $p$, estimable density of the regressors, densities of the random coefficients for which the researcher does not have prior knowledge on the support and which support can be noncompact. 
We assume the marginals of $\mt{\beta}$ (but not of $\alpha$) do not have heavy tails but can have noncompact support. This allows for many parametric families which are used in mixture modelling, while leaving unspecified the parametric family. We do not rely on the Radon transform 
but on the truncated Fourier transform (see, \emph{e.g.}, \cite{Marechal1}). 
Due to \eqref{eindep},  the conditional characteristic function of $Y$ given $\mt{X}=\mt{x}$ at $t$ is the Fourier transform of $f_{\alpha,\mt{\beta}}$ at $(t,t\mt{x}^{\top})^{\top}$. Hence, the family of conditional characteristic functions indexed by $\mt{x}$ in the support of $\mt{X}$ gives access to the Fourier transform of $f_{\alpha,\mt{\beta}}$ on a double cone of axis $(1,0,\dots,0)\in\R^{p+1}$  and apex 0.  When $\alpha=0$, the support of $\mt{\beta}$ and $\mt{X}$ are compact with nonempty interior, this is the problem of out-of-band extrapolation or super-resolution (see, \emph{e.g.}, \cite{Iv_imagaing}). 
Because we do not restrict $\alpha$ and the support of $\mt{\beta}$ can be noncompact, we generalize this approach.  

A related problem is extrapolation. It is used in \cite{Meister07} to perform deconvolution of compactly supported densities allowing the characteristic function of the error to vanish. 
This paper does not use extrapolation or assume densities are analytic. 
Rather, the operator of the inverse problem is a composition of two operators based on partial Fourier transforms. One involves a truncated Fourier transform and we make use of properties of the singular value decomposition. 

Unlike \cite{BM,hoderlein2014triangular}, we go beyond consistency 
and provide a full analysis of the general case. Similar to \cite{GLP,Holtzman,holzmann2019rate}, we study minimax optimality. But, we obtain lower bounds under a wide variety of assumptions. We show that polynomial and nearly parametric rates can be attained. Hence, we can lose little in terms of rates of convergence from going from a parametric model to a nonparametric one. This contrasts with the pessimistic logarithmic rates in \cite{Holtzman} (also mentioned in \cite{hoderlein2014triangular}) and message to avoid estimating densities  
 in \cite{dunker2017tests}.  
We present an estimator involving: series based estimation of the partial Fourier transform of the density with respect to the first variable, interpolation around zero, and inversion of the partial Fourier transform. The orthonormal systems are tensor products of the Prolate Spheroidal Wave Functions  (henceforth PSWF, see  \cite{Osipov}) when the law of $\boldsymbol{\beta}$ has a support included in a known bounded set and else are a new system introduced for this paper and analyzed in \cite{Note}. These systems can also be used in a wide range of applications such as for stable analytic continuation by Hilbert space techniques (see \cite{Note}). We give rates of convergence and use a Goldenshluger-Lepski type method to obtain data-driven estimators. 
We consider estimation of the marginal $f_{\mt{\beta}}$ in Appendix \ref{Amarginals}. We present a numerical method to compute the estimator which is implemented in the \textsf{R} package \href{https://CRAN.R-project.org/package=RandomCoefficients}{RandomCoefficients} with practical details in \cite{randomcoefficients}. The numerical procedure is a fast alternative to the EM algorithm for parametric mixtures of regression models and is robust to misspecification of the parametric family.

\section{Notations}\label{s12}
The notations $\cdot$, $\cdot_1$, $\cdot_2$, $\star$ are used to denote variables in a function. 
$a\wedge b$ (resp. $ a\vee b$) is used for the minimum (resp. maximum) between $a$ and $b$, $(\cdot)_+$ for $0\vee\cdot$,  and $\indic\{A\}$ for the indicator function of set $A$. 
$\N$ and $\N_0$ stand for the positive and nonnegative integers. Bold letters are used for vectors. For all $r\in\R$, $\underline{\mt{r}}$ is the vector, which dimension will be clear from the text, where each entry is $r$. For $x\ge1$ we denote by $\ln_2(x) = \ln(\ln(x))$. $\mathcal{W}$ is the inverse of $x\in[0,\infty)\mapsto xe^x$. 
$|\cdot|_q$ for $q\in[1,\infty]$ stands for the $\ell_q$ norm of a vector or sequence. 
For all $\mt{\beta}\in\mathbb{C}^d$, $(f_m)_{m\in\N_0}$ functions with values in $\C$, and $\boldsymbol{m}\in\mathbb{N}_0^d$, denote by $\mt{\beta}^{\boldsymbol{m}}=\prod_{k=1}^d\mt{\beta}_k^{\mt{m}_k}$, $|\mt{\beta}|^{\boldsymbol{m}}=\prod_{k=1}^d|\mt{\beta}_k|^{\mt{m}_k}$, and $f_{\boldsymbol{m}}=\prod_{k=1}^df_{\mt{m}_k}$.
For a function $f$ of real variables, $\mathrm{supp}(f)$  denotes its support.  
The inverse of a mapping $f$, when it exists, is denoted by $f^I$.  We denote the interior of $\mathcal{S}\subseteq\R^d$ by $\overset{\circ}{\mathcal{S}}$. 
When $\mathcal{S}$ is measurable and $\mu$ a nonnegative function from $\mathcal{S}$ to $[0,\infty]$, $L^2(\mu)$ is the space of complex-valued square integrable functions equipped with 
$ \braket{f,g}_{L^2(\mu)}=\int_{\mathcal{S}} f(\mt{x})\overline{g}(\mt{x}) \mu(\mt{x})d\mt{x}$. This is denoted by  $L^2(\mathcal{S})$ when $\mu=1$. When $W_{\mathcal{S}}=\indic\{\mathcal{S}\}+\infty\ \indic\{\mathcal{S}^c\}$, we have $L^2(W_{\mathcal{S}})= \{ f \in L^2(\R^d): \  \supp(f) \subseteq \mathcal{S} \}$ and $\braket{f,g}_{L^2(W_{\mathcal{S}})} = \int_{\mathcal{S}} f(\mt{x})\overline{g}(\mt{x}) d\mt{x}$. Denote by $\mathcal{D}$ the set of densities  and by $\otimes$ the product of functions (\emph{e.g.}, $ W^{\otimes d}(\mt{b}) = \prod_{j=1}^{d} W(\mt{b}_j)$) or measures.
The Fourier transform of $f\in L^1\left(\mathbb{R}^{d}\right)$ is $\mathcal{F}\left[f\right](\mt{x})=\int_{\mathbb{R}^{d}}e^{i\mt{b}^{\top}\mt{x}}f(\mt{b})d\mt{b}$ and $\mathcal{F}\left[f\right]$ is also the Fourier transform in  $L^2\left(\mathbb{R}^{d}\right)$. For all $ c>0 $, denote the Paley-Wiener space 
by $PW(c) :=  \left\{ f \in L^2(\R): \mbox{supp}\left(\mathcal{F}\left[f\right]\right) \subseteq [-c,c]\right\} $, by $\mathcal{P}_c$ the projector from $ L^2(\R) $ to $ PW(c) $ ($\mathcal{P}_c[f]  = \mathcal{F}^I\left[ \indic\{[-c,c]\} \mathcal{F}\left[ f\right] \right]$), and 
\begin{equation}\label{eqdef} \begin{array}{cccc}\forall c\ne0,\quad
\mathcal{F}_{c} : & L^{2}\left(W^{\otimes d}\right) & \rightarrow &  L^{2}\left([-1,1]^d\right)			 \\
& f 		& \rightarrow &   \mathcal{F}\left[  f \right](c \ \cdot)    .
\end{array} 
\end{equation}
$\mathcal{F}_{1\mathrm{st}}\left[f\right](t,\cdot_2)$ denotes the partial Fourier transform of $f$ with respect to the first variable. For a random vector $\mt{X}$, $\mathbb{P}_{\mt{X}}$ is its law, $f_{\mt{X}}$ its density, $f_{\mt{X}|\mathcal{X}}$ the truncated density of $\mt{X}$ given $\mt{X}\in\mathcal{X}$, $\mathbb{S}_{\mt{X}}$ its support, and $f_{Y|\mt{X}=\mt{x}}$ 
the conditional density. For a sequence of random variables $ \left(X_{n_0,n}\right)_{(n_0,n)\in\N_0^2} $,  $X_{n_0,n} = \underset{\mathcal{U}}{O_p}(1) $ means that, for all $ \epsilon >0 $, there exists $M$ such that $ \mathbb{P}(|X_{n_0,n}| \ge M) \le \epsilon $ for all $ (n_0,n)\in\N_0^2$ such that $\mathcal{U}$ holds. In the absence of constraint, we drop the notation   $\mathcal{U}$. With a single index ${O_p}(1)$ is the usual  notation. 

\section{Preliminaries}\label{s2}
\begin{assumption}\label{ass:compact}
	\begin{enumerate}[\textup{(}{H1.}1\textup{)}] 
		\item\label{E3} $f_{\mt{X}}$ and $f_{\alpha,\mt{\beta}}$ exist;
		\item\label{E2}$f_{\alpha,\mt{\beta}}\in L^2(w\otimes \overline{W}^{\otimes p})$, where $w\ge1$ and $\overline{W}=e^{|\cdot|/R}$, where $R>0$;
		\item\label{E5} There exists $x_0>0$ and $\mathcal{X}=[-x_0,x_0]^p \subseteq \mathbb{S}_{\mt{X}}$ and we have at our disposal i.i.d $ (Y_i, \mt{X}_i)_{i=1}^n $ and an estimator $ \widehat{f}_{\mt{X}|\mathcal{X}} $ based on $ \mathcal{G}_{n_0} = (\mt{X}_i)_{i=-n_0+1}^{0} $ independent of $ (Y_i, \mt{X}_i)_{i=1}^n$;
		\item\label{E4}  $\mathcal{E}$ is a set of densities on $\mathcal{X}$ such that, for $c_{\mt{X}},C_{\mt{X}}\in(0,\infty)$, for all $f\in\mathcal{E}$,  $\left\|f\right\|_{L^{\infty}(\mathcal{X})}\le C_{\mt{X}}$ and
		$\left\|1/f\right\|_{L^{\infty}(\mathcal{X})}\le c_{\mt{X}}$, and, for $(v({n_0},\mathcal{E}))_{n_0\in\N}\in(0,1)^{\N}$ which tends to 0, we have
		\begin{equation*}\label{eq:unknown}
		\frac{1}{ v({n_0},\mathcal{E})}\sup_{f_{\mt{X}|\mathcal{X}}\in \mathcal{E}}  \left\|\widehat{f}_{\mt{X}|\mathcal{X}}- f_{\mt{X}|\mathcal{X}}  \right\|_{L^{\infty}(\mathcal{X})}^2 = O_p\left(1\right). 
		\end{equation*}
	\end{enumerate}
\end{assumption}
We maintain this assumption for all 
upper bounds. 
If $w^{-1}\in L^1(\R)$, (H\ref{ass:compact}.\ref{E2}) implies that the slopes of $\mt{\beta}$ do not have heavy tails. This means that their tails are not heavier than that of the exponential distribution (\emph{i.e.}, their Laplace transform is finite near 0). Indeed, we have, for all $\epsilon\in(0,1)$ and $k=1,\dots,p$, for $\lambda=(1-\epsilon)/(2R)$, by the Cauchy-Schwarz inequality,  
\begin{align*}
\E\left[ e^{\lambda \mt{\beta}_k}\right] &\leq  \E\left[e^{\lambda \left|\mt{\beta}_k\right|}\right]
\le \left\|f_{\alpha,\mt{\beta}}\right\|_{L^2\left(w\otimes W^{\otimes p}\right)}
\left\|w^{-1}\right\|_{L^1(\R)}^{1/2} \left(
2R/\epsilon\right)^{p/2}<\infty. 
\end{align*}
Now on, $W$ is either $W_{[-R,R]}$ or $\cosh(\cdot/R)$. It is 
such that $L^2\left(w\otimes W^{\otimes p}\right)\subseteq L^2(w\otimes \overline{W}^{\otimes p})$. When $W=W_{[-R,R]}$, $f_{\alpha,\mt{\beta}}\in L^2(w\otimes W^{\otimes p})$ implies that $\mathbb{S}_{\mt{\beta}}\subseteq[-R,R]^p$.
The condition $\mathcal{X}=[-x_0,x_0]^p \subseteq \mathbb{S}_{\mt{X}}$ in (H1.\ref{E4}) is not restrictive because 
$Y=\alpha+\mt{\beta}^{\top}\underline{\mt{x}}+\mt{\beta}^{\top}(\mt{X}-\underline{\mt{x}})$, we can take $ \underline{\mt{x}}$ and $ x_0$ such that  $\mathcal{X}\subseteq \mathbb{S}_{\mt{X}- \underline{\mt{x}}}$, and there is a one-to-one mapping between $f_{\alpha+\mt{\beta}^{\top}\underline{\mt{x}},\mt{\beta}}$ and $f_{\alpha,\mt{\beta}}$. We assume (H1.\ref{E4}) because the estimator involves estimators of $f_{\mt{X}|\mathcal{X}}$ in denominators. Alternative solutions exist when $p=1$ (see,  \emph{e.g.}, \cite{kerkyacharian2004regression}) only. Assuming the availability of an estimator of  $f_{\mt{X}|\mathcal{X}}$ using the preliminary sample $\mathcal{G}_{n_0}$ is common in the deconvolution literature (see, \emph{e.g.}, \cite{Comte_Lacour}). By using estimators of $f_{\mt{X}|\mathcal{X}}$ for a well chosen $\mathcal{X}$ rather than of $f_{\mt{X}}$, the assumption $\|f_{\mt{X}|\mathcal{X}}\|_{L^{\infty}(\mathcal{X})}\le C_{\mt{X}}$ and $\|1/f_{\mt{X}|\mathcal{X}}\|_{L^{\infty}(\mathcal{X})}\le c_{\mt{X}}$ in (H1.\ref{E4}) becomes mild. 

\subsection{Inverse problem in Hilbert spaces}\label{s21}
Estimation of $f_{\alpha,\mt{\beta}}$ is a statistical ill-posed inverse problem. The operator depends on $w$ and $W$. 
We have, for all $t\in\R$ and $\mt{u}\in[-1,1]^p$,
$\mathcal{K} f_{\alpha,\mt{\beta}}(t,\mt{u}) = \mathcal{F}\left[f_{Y|\mt{X}=x_0\mt{u}}\right](t)x_0|t|^{p/2} $, where \begin{equation}\label{opglobal} 
\begin{array}{ccccc}
\mathcal{K} : &  L^2\left(w\otimes W^{\otimes p}\right)& \rightarrow &  L^{2}(\R\times[-1,1]^p) &   \\
& f & \rightarrow &  (t,\mt{u}) \mapsto \mathcal{F}\left[f\right](t,x_0t\mt{u})x_0|t|^{p/2} .
\end{array}
\end{equation}
\begin{proposition}\label{prop:notcompact} 
	$L^2\left(w\otimes W^{\otimes p}\right)$ is continuously embedded into $L^2(\R^{p+1})$. Moreover, $\mathcal{K}$ is injective and continuous,  and not compact if $w=1$. 
\end{proposition}
The case $w=1$ corresponds to mild integrability in the first variable, there the SVD of $\mathcal{K}$ does not exist. This makes it difficult to prove rates of convergence even for estimators which do not rely explicitly on the SVD such as the Tikhonov and Landweber method (Gerchberg algorithm in out-of-band extrapolation, see \cite{Iv_imagaing}). Rather than work with $\mathcal{K}$ directly, we use that $\mathcal{K}$ 
is the composition of operators which are easier to analyze
\begin{equation}\label{e:decomp}\mathrm{for}\ t\in\R,\ \mathcal{K}[f](t,\star)=\mathcal{F}_{t x_0}\left[\mathcal{F}_{1\mathrm{st}}\left[f\right](t,\cdot_2)\right](\star)x_0|t|^{p/2} \ 
\mathrm{in}\ L^2([-1,1]^p).
\end{equation}
For all $f\in L^2\left(w\otimes W^{\otimes p}\right)$ 
and $t \in \R$, $\mathcal{F}_{1\mathrm{st}}\left[f\right](t,\cdot_2)$ belongs to $L^2(W^{\otimes p})$ and, for $c\ne0$, $\mathcal{F}_{c}:\ L^2(W^{\otimes p})\to L^2([-1,1]^p)$ admits a SVD, where both orthonormal systems are complete. This is a tensor product of the SVD when $p=1$ that we denote by $(\sigma_m^{W,c},\varphi_m^{W,c},g_m^{W,c})_{m\in\N_0}$, where $(\sigma_m^{W,c})_{m\in \N_0}\in(0,\infty)^{\N_0}$ is in decreasing order repeated according to multiplicity. 


\begin{proposition}\label{sec:upper:extension} 
	For all $c\ne0$, $(g_m^{W,c})_{m\in \N_0}$ and $(\varphi_m^{W,c})_{m\in \N_0}$ are bases of, respectively, $L^2([-1,1])$ and  $L^2(W	)$. 
\end{proposition}
The singular functions $(g_m^{W_{[-1,1]},c})_{m\in \N_0}$ are the PSWF. They can be extended as entire functions in $L^2(\R)$ and form a complete orthogonal system of $PW(c)$ for which we use the same notation. They are useful to carry interpolation and extrapolation (see, \emph{e.g.}, \cite{lindberg2012mathematical}) with Hilbertian techniques. In this paper, for all $t\ne0$, $\mathcal{F}_{\rm{1st}}\left[f_{\alpha,\boldsymbol{\beta}}\right](t,\cdot_2)$ plays the role of the Fourier transform in the definition of $PW(c)$.  The weight $\cosh(\cdot/R )$ allows for larger classes than $PW(c)$ and noncompact $\mathbb{S}_{\mt{\beta}}$. This is useful even if $\mathbb{S}_{\mt{\beta}}$ is compact when the researcher does not know a superset containing $\mathbb{S}_{\mt{\beta}}$. The results on the corresponding SVD and a numerical algorithm to compute it are 
in \cite{Note}. 

\subsection{Sets of smooth and integrable functions}\label{s22}
Define, for $q\in\{1,\infty\}$, 
\begin{align*}
b_{\boldsymbol{m}}(t):=\braket{ 
	\mathcal{F}_{1\mathrm{st}}\left[f\right](t, \cdot_2),\varphi_{\boldsymbol{m}}^{W,x_0t}
}_{L^2\left(W^{\otimes p}\right)},
\theta_{q,k}(t) := \left(\sum_{\mt{m}\in\N_0^p: \ \abs{\mt{m}}_{q} = k} \left | b_{\boldsymbol{m}}(t) \right |^2\right)^{1/2},
\end{align*} 
and, for all $(\phi(t))_{t\ge0}$ and $\left(\omega_m\right)_{m\in\N_0}$ increasing, $ \phi(0)=\omega_0 = 1$, $l,M>0$, $ t\in\R,\ \mt{m} \in\N_0^p$, $k\in\N_0$, $\mathcal{I}_{w,W}(M) := \{f : \ \|f\|_{L^2\left(w\otimes W^{\otimes p}\right)}\le M\}$,  and 
$$\mathcal{H}^{q,\phi,\omega}_{w,W}(l,M):= \left\{f:
\sum_{k\in \N_0}  \int_{\R} \phi^2(|t|)\theta_{q,k}^2(t) dt \bigvee
\sum_{k\in\N_0}\omega_{k}^2\|\theta_{q,k}\|_{L^2(\R)}^2
\leq 2 \pi l^2\right\}\bigcap\mathcal{I}_{w,W}(M).$$

We use the notation $\mathcal{H}^{q,\phi,\omega}_{w,W}(l)$ when we require $ \|f\|_{L^2\left(w\otimes W^{\otimes p}\right)}< \infty $ rather than $ \|f\|_{L^2\left(w\otimes W^{\otimes p}\right)} \leq M$.  
The set $\mathcal{I}_{w,W}(M)$ imposes the integrability discussed in the beginning fo the section. 
The first set in the definition of $\mathcal{H}^{q,\phi,\omega}_{w,W}(l, M)$ defines the notion of smoothness 
analyzed in this paper. It involves a maximum, thus two inequalities: the first for smoothness in the first variable
and the second for smoothness in the other variables. The asymmetry in the treatment of the first and remaining variables is due to the fact that 
only the random slopes are multiplied by regressors which have limited variation and we make integrability assumptions in the first variable which are as mild as possible. The smoothness classes in the analysis of the Radon transform usually involve nonstandard weight functions well suited to the operator. In contrast, the ones in this paper are not too hard to interpret.  The first inequality can be written as
$$ \int_{\R} \phi^2(|t|)\left\|\mathcal{F}_{1\mathrm{st}}\left[f\right](t, \cdot_2)\right\|_{L^2\left(W^{\otimes p}\right)}^2 dt \le 2 \pi l^2,$$
so it is the usual Sobolev smoothness of functions in $L^2(\R;L^2(W^{\otimes p}))$. 
We show in Appendix \ref{ASob} that when, for almost every $a$, $\boldsymbol{b}\mapsto f(a,\boldsymbol{b})$ has compact support, it is possible to use, instead of the PSWF basis, the basis giving rise to Fourier series (the characters of the torus) and $(\omega_m)_{m\in\N_0}$ are the usual weights used for Sobolev smoothness. So when $\phi$ and $\omega_m$ are monomials, smoothness corresponds to the function having bounded sum of squared $L^2$ norm of partial derivatives of degree that of the monomial. When these are exponentials, it implies that all partial derivatives are square integrable. It corresponds to supersmooth classes (see, \emph{e.g.}, \cite{Cavalier2}). In the main text we replace the characters of the torus by the bases $(\varphi_{\boldsymbol{m}}^{W,x_0t})_{\mt{m}\in\N_0^p}$ for $t\ne0$. So we can treat, using the same framework involving a sum, the case where $\mathbb{S}_{\boldsymbol{\beta}}$ is included in a known bounded set and where such set is unknown or the support can be noncompact. The second inequality can be written, there exists a density $\overline{\phi}$  on $\R$ such that, for almost every $t\in\R$,   
$$\sum_{\mt{m}\in\N_0^p}\omega_{\abs{\mt{m}}_{q}}^2\left | b_{\boldsymbol{m}}(t) \right |^2
\leq \overline{\phi}(t) 2 \pi l^2.$$
For fixed $t$ this is a source condition based on the truncated Fourier operator. Proceeding as in (21) and (22) in \cite{Note} allows to rewrite this condition as a smoothness condition on the Fourier transform. Because, for all $m\in\N$, we have $\omega_m\ge 1$, nonsmooth functions have analytic Fourier transforms. Smooth functions involve weights $\omega_m$ which are monomials or exponentials, they have a Fourier transform in a smaller class of analytic functions. 
We analyze all types of smoothness. Because smoothness is unknown anyway, we provide an adaptive estimator. 
We analyze two values of $q$ and show its value 
matters for the rates of convergence for supersmooth functions. 
\begin{remark}\label{rmrk1} The next model is related to \eqref{eRC} 
when $f_{\mt{X}}$ is known: 
	\begin{equation}\label{eq:modalt}
	dZ(t) = \mathcal{K}\left[f\right](t,\cdot_2) dt + \frac{\sigma}{\sqrt{n}}dG(t),\quad t\in\R,
	\end{equation}
	where $f$ plays the role of $f_{\alpha,\mt{\beta}}$, $\sigma>0$ is known and the partial derivative in the sense of distributions with respect to time  of $G$ is the space time white noise in $L^{2}(\R\times[-1,1]^p)$. A usual mathematical formalisation of the space time white noise (see \cite{da2014stochastic}) is that and $(G(t))_{t\in\R}$ is a complex two-sided cylindrical Gaussian process on $ L^2([-1,1]^p)$. This means,  for $\Phi$ Hilbert-Schmidt from $L^2([-1,1]^p)$ to a separable Hilbert space $H$, $(\Phi G(t))_{t\in\R}$ is a Gaussian process in $H$ of covariance $\Phi\Phi^*$. Taking $\Phi G (t)=\sum_{\mt{m}\in\N_0^p} \Phi[g_{\mt{m}}^{W,x_0t}]B_{\mt{m}}(t)$, where $ B_{\mt{m}}(t) =B_{\mt{m}}^{\mathfrak{R}}(t) + i B_{\mt{m}}^{\mathfrak{I}}(t)$, $ (B_{\mt{m}}^{\mathfrak{R}}(t))_{t\in\R}  $ and $ (B_{\mt{m}}^{\mathfrak{I}}(t))_{t\in\R}  $ are  independent two-sided Brownian motions, the system of independent equations
	\begin{equation}\label{eq:modalt2}
	Z_{\mt{m}}(t) = \int_{0}^{t} \sigma_{\mt{m}}^{W,x_0s} b_{\mt{m}}(s) ds + \frac{\sigma}{\sqrt{n}} B_{\mt{m}}(t),\quad t\in\R,
	\end{equation}  
	where,  $Z_{\mt{m}}: = \langle Z(\star),  g_{\mt{m}}^{W,x_0\star} \rangle_{L^2([-1,1]^p)}$ and $\mt{m} \in \N_0^p$, 
	is equivalent to \eqref{eq:modalt}. Because $\sigma_{\mt{m}}^{W,x_0s} $ is small when $|\mt{m}|_{q}$ is large or $s$ is small (see Lemma \ref{upper_bound}), the estimator of Section \ref{s23} truncates large values of $|\mt{m}|_{q}$ and does not rely on small values of $|s|$ but uses interpolation. 
\end{remark}

\subsection{Interpolation}\label{sec:interpol}
Denote, for all $m\in \N_0$ and $c\neq0$, $\rho_{m}^{W,c}:  =2\pi(\sigma_{m}^{W,c})^2/\abs{c} $. Define, for all $\underline{a},\epsilon>0$ the operator on $L^2(\R)$ with domain $PW(\underline{a})$
\begin{align}\label{eq:inter}
\mathcal{I}_{\underline{a},\epsilon}\left[f\right]&:=\sum_{m \in \N_0} \frac{\rho_{m} ^{W_{[-1,1]},\underline{a}\epsilon}\braket{f,\mathcal{C}_{1/\epsilon}\left[g^{W_{[-1,1]},\underline{a}\epsilon}_m\right]}_{L^{2}(\R\setminus(-\epsilon,\epsilon))}}{\left(1 - \rho_{m} ^{W_{[-1,1]},\underline{a}\epsilon}\right)\epsilon}  \mathcal{C}_{1/\epsilon}\left[g^{W_{[-1,1]},\underline{a}\epsilon}_m\right]. 
\end{align} 
\begin{proposition}\label{prop:interp}
	For all $\underline{a},\epsilon>0$, 
	$ \mathcal{I}_{\underline{a},\epsilon}( L^2(\R)) \subseteq L^2([-\epsilon,\epsilon])$  
	and, for all $g\in PW(\underline{a})$, $\mathcal{I}_{\underline{a},\epsilon}[g]=g$ in $L^2(\R)$ and, for $C_0:=4\cdot/(\pi(1-\rho_0^{W_{[-1,1]},\cdot})^2)$ and all $f,h\in L^2(\R)$, 
	\begin{align}\label{eq:interpol}
	\left\|f- \mathcal{I}_{\underline{a},\epsilon}\left[ h\right] \right\|_{L^2([-\epsilon,\epsilon])}^2
	\le &2(1 + C_0(\underline{a}\epsilon))\left\|f-\mathcal{P}_{\underline{a}}[f]\right\|_{L^2(\R)}^2 
	\\&
	+2 C_0(\underline{a}\epsilon)\left\|f-h\right\|_{L^2(\R\setminus(-\epsilon,\epsilon))}^2.\notag 
	\end{align}
\end{proposition} If $f\in PW(\underline{a})$, 
$\mathcal{I}_{\underline{a},\epsilon}[f]$ only relies on $f\indic\{\R\setminus(-\epsilon,\epsilon)\}$ 
and $\mathcal{I}_{\underline{a},\epsilon}[f]=f$ on $\R\setminus(-\epsilon,\epsilon)$, 
so \eqref{eq:inter} provides an analytic formula to carry interpolation on $[-\epsilon,\epsilon]$ of functions in $PW(\underline{a})$. 
Else, \eqref{eq:interpol} provides an upper bound on the error made by approximating $f$ by $\mathcal{I}_{\underline{a},\epsilon}\left[ h\right]$ on $[-\epsilon,\epsilon]$ 
when $h$ approximates $f$ outside $[-\epsilon,\epsilon]$.  We use interpolation when the variance of an initial estimator $\widehat{f}^{0}$ of $f$ is large due to its values near 0 but $\|f-\widehat{f}^{0}\|_{L^2(\R\setminus(-\epsilon,\epsilon))}^2$ is small and work with 
$
\widehat{f}(t)=\widehat{f}^{0}(t)\indic\{|t|\ge\epsilon\}+\mathcal{I}_{\underline{a},\epsilon}[\widehat{f}^{0}](t)\indic\{|t|<\epsilon\}.$ 
Then, \eqref{eq:interpol} yields
\begin{align}\notag
\left\|f-\widehat{f}\right\|_{L^2(\R)}^2\le&
(1+2C_0(\underline{a}\epsilon))\left\|f-\widehat{f}^{0}\right\|_{L^2(\R\setminus(-\epsilon,\epsilon))}^2 \\
&+2(1 + C_0(\underline{a}\epsilon))\left\|f-\mathcal{P}_{\underline{a}}[f]\right\|_{L^2(\R)}^2.\label{eq:interpol6}
\end{align}
When $\mathrm{supp}\left(\mathcal{F}[f]\right)$ is compact, $\underline{a}$ is taken such that $\mathrm{supp}\left(\mathcal{F}[f]\right)\subseteq[-\underline{a},\underline{a}]$. Else, 
$\underline{a}$ goes to infinity so the second term in \eqref{eq:interpol6} goes to 0. $\epsilon$ is taken such that $ \underline{a}\epsilon$ is constant because, due to (3.87) in \cite{Osipov}, $\lim_{t\to\infty}C_0(t)=\infty$ and \eqref{eq:interpol} and \eqref{eq:interpol6} become useless. 
We set 
$C = 2(1 + C_0(\underline{a}\epsilon))$. 

\subsection{Risk}
The risk 
is the mean integrated squared error (MISE)
\begin{equation*}
\mathcal{R}_{n_0}^{W}\left(\widehat{f}_{\alpha,\mt{\beta}},f_{\alpha,\mt{\beta}}\right)  : = \E \left[\left\|  \widehat{f}_{\alpha,\mt{\beta}} -f_{\alpha,\mt{\beta}} \right\|_{L^2\left(1\otimes W^{\otimes p}\right)}^2  \middle| \mathcal{G}_{n_0}  \right].
\end{equation*}
It is $\E [\|\widehat{f}_{\alpha,\mt{\beta}}-f_{\alpha,\mt{\beta}}\|_{L^2(\R^{p+1})}^2 | \mathcal{G}_{n_0}]$ when $W= W_{[-R,R]}$ and $\supp(\widehat{f}_{\alpha,\mt{\beta}})\subseteq \R\times [-R,R]^p$, 
else it is
\begin{equation}\label{eq:R_control}
\E \left[\left\|  \widehat{f}_{\alpha,\mt{\beta}} -f_{\alpha,\mt{\beta}} \right\|_{L^2(\R^{p+1})}^2  \middle| \mathcal{G}_{n_0}  \right]\le   \left\|1/W\right\|_{L^{\infty}(\R)}^p\mathcal{R}_{n_0}^{W}\left(\widehat{f}_{\alpha,\mt{\beta}},f_{\alpha,\mt{\beta}}\right).
\end{equation}

\section{Lower bounds}
\label{s41}
The lower bounds involve a function $r$ (for rate) and take the form 
\begin{equation}\label{eq:lbtype}\exists \nu>0:\  
 \underline{\lim}_{n\to\infty} 
	\underset{\widehat{f}_{\alpha,\mt{\beta}}}{\inf} \ \  \underset{f_{\alpha,\mt{\beta}}  \in  \mathcal{H}^{q,\phi,\omega}_{w,W}(l)\cap\mathcal{D}}{\sup}\frac{ \E \left[\left\|  \widehat{f}_{\alpha,\mt{\beta}} -f_{\alpha,\mt{\beta}} \right\|_{L^2(\R^{p+1})}^2 \right]}{r(n)^2} \geq \nu.
\end{equation}
When we replace $f_{\alpha,\mt{\beta}}$ by $f$, $\widehat{f}_{\alpha,\mt{\beta}}$ by $\widehat{f}$, remove $\mathcal{D}$ from the nonparametric class, and consider model \eqref{eq:modalt2}, we refer to (\ref{eq:lbtype}'). We use 
$k_q := \indic\{q=1\}+p\indic\{q=\infty\}$ and $k_q'= p+1-k_q$. We focus on the lower bounds for polynomial and exponential weights $(\omega_k)_{k\in \N_0}$ which yield the usual smooth and supersmooth cases. 
To be comparable to rates in other inverse problems, the exponential weight $(\omega_k)_{k\in \N_0}$ has the same form as the rate of decay as the singular values (see Lemma \ref{upper_bound} and Theorem 7 in \cite{Note}), hence the different forms in (T\ref{theo:lower}.1\ref{it:expindic}) and  (T\ref{theo:lower}.2\ref{it:expch}) due to the different values of $W$. 

\begin{theorem}\label{theo:lower}
	For all $q\in \{1,\infty\}$, $ \phi$ increasing on $[0,\infty)$,  
	$0< l,s,R, \kappa< \infty$,  
	and $w$ such that $\int_1^{\infty}w(a)/a^4<\infty$, if $W=W_{[-R,R]}$, 
	\begin{enumerate} [\textup{(}{T1.1}a\textup{)}] 
		\item\label{it:lowerindic} $(\omega_k)_{k\in \N_0} =(k^{\sigma})_{k\in \N_0}  $, $\phi$ is such that $\lim_{\tau\to\infty}\int_{0}^{\infty}\phi(t)^2e^{-2\tau t}dt=0$,  
		$f_{\mt{X}}$ is known, $\mathbb{S}_{\mt{X}}=\mathcal{X}$, and $\left\|f_{\mt{X}}\right\|_{L^{\infty}(\mathcal{X})}<\infty$, then (\ref{eq:lbtype}) holds with $r(n)=\left(\ln(n)/\ln_2(n)\right)^{-(2+k_q/2)\vee\sigma}$,
		\item\label{it:expindic} we consider model \eqref{eq:modalt2}, $(\omega_k)_{k\in \N_0} =  \left( e^{\kappa k\ln\left(k+1 \right)}\right)_{k\in \N_0} $, then (\ref{eq:lbtype}') holds with $r(n)=n^{-\kappa/(2\kappa+2 k_q)}/\ln(n)$,
	\end{enumerate}
	else if $W =\cosh(\cdot/R)$, we consider model \eqref{eq:modalt2}, 
	\begin{enumerate} [\textup{(}{T1.2}a\textup{)}] 
		\item\label{it:lowerch}  $(\omega_k)_{k\in \N_0} =(k^{\sigma})_{k\in \N_0}  $, for all $\overline{\sigma}>1/2$, then  (\ref{eq:lbtype}') holds with $r(n)=\ln\left(n/\ln(n)\right)^{-\overline{\sigma}\vee\sigma}$
		\item\label{it:expch} $(\omega_k)_{k\in \N_0} =  \left( e^{\kappa k}\right)_{k\in \N_0} $, then (\ref{eq:lbtype}') holds with $r(n)=n^{-\kappa/(2\kappa+2k_q)}$.
	\end{enumerate}
\end{theorem}
%
%
%

By \eqref{eq:R_control}, (T1.2\ref{it:lowerch}), and (T1.2\ref{it:expch}), we obtain lower bounds involving $\mathcal{R}_{n_0}^{W}$. 
We relate those rates to those in other inverse problems after theorems \ref{theo:compact} and \ref{theo:non_compact}. Importantly, for sufficiently smooth classes of functions, polynomial rates can be attained, for this severely ill-posed inverse problem. 

\section{Estimation}
\subsection{Estimator}\label{s23}
For all $q\in \{1,\infty\}$, $0<\epsilon<1<T$, $\underline{N}\in\R^{\R}$, $N(t) =\left\lfloor \underline{N}(t)\right\rfloor $ for $\epsilon \leq |t|\leq T$, $N(t)= N(\epsilon)$ for $|t|\le\epsilon$ and  $N(t)= N(T)$ for $|t|>T$, a regularized inverse is obtained by: 
\begin{enumerate} [\textup{(}{S.}1\textup{)}] 
	\item\label{S1} for all $t \ne0$, obtain a first approximation of $F_1(t, \cdot):= \mathcal{F}_{1\mathrm{st}}\left[f_{\alpha,\mt{\beta}}\right](t, \cdot)$ 
	\begin{align}
	F_1^{q,N,T,0}(t,\cdot_2)&:=\indic\{|t|\le T\}\sum_{\abs{\mt{m}}_{q}\leq N(t)}\frac{c_{\mt{m}}(t)}{\sigma_{\mt{m}}^{W,x_0t}}\varphi_{\mt{m}}^{W,x_0t},\notag\\
	c_{\mt{m}}(t)&:= \big\langle \mathcal{F}\left[f_{Y|\mt{X}=x_0\cdot}\right](t), g_{\mt{m}}^{W,x_0t}\big\rangle_{L^2([-1,1]^p)},\notag
	\end{align}
	\item\label{S2} for all $t \in [-\epsilon, \epsilon]$, we use $$F_1^{q,N,T,\epsilon}(t, \cdot):=F_1^{q,N,T,0}(t, \cdot)\indic\{|t|\ge \epsilon\}+ \mathcal{I}_{\underline{a},\epsilon}\left[F_1^{q,N,T,0}(\star, \cdot)\right](t)\indic\{ |t|< \epsilon\},$$
	\item\label{S3} $f_{\alpha,\mt{\beta}}^{q,N,T,\epsilon}(\cdot_1,\cdot_2):=\mathcal{F}_{1\mathrm{st}}^I\left[F_1^{q,N,T,\epsilon}(\star, \cdot_2)\right](\cdot_1)$.
\end{enumerate}
For $|t|\le T$, (S.\ref{S1}) is obtained from \eqref{e:decomp} and a regularized inverse of the truncated Fourier operator $\mathcal{F}_{tx_0}$. It involves spectral cut-off. The indicator $\indic\{|t|\le T\}$ is used in (S.\ref{S3}) as a standard regularization device when inverting the Fourier transform which consists in removing high frequencies.  
To deal with the statistical problem, we 
replace $c_{\mt{m}}$ by 
\begin{equation}\label{eq:change2}\widehat{c}_{\mt{m}}:=\frac{1}{n} \sum_{j=1}^n \frac{e^{i\star Y_j}}{x_0^p \widehat{f}^{\delta}_{\mt{X}|\mathcal{X}}(\mt{X}_j)} \overline{g_{\mt{m}}^{W,x_0\star}}\left(\frac{\mt{X}_j}{x_0}\right)\indic\left\{\mt{X}_j\in\mathcal{X}\right\},\end{equation}
where $\widehat{f}^{\delta}_{\mt{X}| \mathcal{X}}(\mt{X}_j):=\widehat{f}_{\mt{X}| \mathcal{X}}(\mt{X}_j)\vee \sqrt{\delta(n_0)}$ and $ \delta(n_0)$ is a trimming factor converging to zero. 
This yields the estimators $\widehat{F}_1^{q,N,T,0}$, $\widehat{F}_1^{q,N,T,\epsilon}$, and $ \widehat{f}^{q,N,T, \epsilon}_{\alpha,\mt{\beta}} $. 
Because inverting the truncated Fourier operator $\mathcal{F}_{tx_0}$ is more ill-posed near $0$ (see 
	Lemma \ref{upper_bound} and Theorem 7 in \cite{Note}), 
	$\widehat{F}_1^{q,N,T,0}$ has a large variance for $t\in[-\epsilon,\epsilon]$. Hence we use interpolation (see  
	Section \ref{sec:interpol}). 
	
We use $(\widehat{f}^{q,N,T, \epsilon}_{\alpha,\mt{\beta}})_+$ as a final estimator of $f_{\alpha,\mt{\beta}}$  which has a smaller risk than $\widehat{f}^{q,N,T, \epsilon}_{\alpha,\mt{\beta}}$ (see \cite{GK,Tsybakov}). 
We use $n_e = n\wedge \left\lfloor \delta(n_0)/v(n_0,\mathcal{E})\right\rfloor$ for the sample size required for an ideal estimator where $f_{\mt{X}|\mathcal{X}}$ is known to achieve the rate of the plug-in estimator. The upper bounds below take the form 
\begin{equation}\label{eq:weight2_res2}
\sup_{f_{\alpha,\mt{\beta}} \in \mathcal{H}^{q,\phi,\omega}_{w,W}(l, M)\cap\mathcal{D},\ f_{\mt{X}|\mathcal{X}} \in \mathcal{E}} \frac{\mathcal{R}_{n_0}^{W}\left(\widehat{f}_{\alpha,\mt{\beta}}^{q,N,T,\epsilon},f_{\alpha,\mt{\beta}}\right)}{r(n_e)^2}  = O_p(1).
\end{equation} 
With 
the restriction $f_{\alpha,\mt{\beta}} \in \mathcal{H}^{q,\phi,\omega}_{w,W}(l)\cap\mathcal{D} $, we refer to (\ref{eq:weight2_res2}'). 

\subsection{Upper bounds}
\label{s43}
We use $ T = \phi^{I}(\omega_{\underline{N}})$, $\underline{a} =w^{I}(\omega_{\underline{N}}^2) $ when $w\neq W_{[-\underline{a},\underline{a}]}$, 
for $u>0$, $K_{\underline{a}}(u):=\underline{a} \indic\{ w\neq W_{[-\underline{a},\underline{a}]}\} +u \indic\{ w= W_{[-\underline{a},\underline{a}]}\}$.
Below $\underline{N}$, hence $N$, is constant.
\begin{theorem}\label{theo:compact} 
	Let  $W=W_{[-R,R]}$. For all $q\in \{1,\infty\}$,  $ l,M ,s,R,\sigma, \kappa,\mu, \gamma, \nu  >0$,  $\mathbb{S}_{\mt{\beta}}\subseteq [-R,R]^p$, $\underline{N}$ solution of 
	$2k_q\underline{N}\ln\left(\underline{N}K_{\underline{a}}(1)\right) + \ln(\omega_{\underline{N}}^2) +(p-1)\ln\left(\underline{N}\right)= \ln(n_e)$, and $\epsilon=7e\pi/(Rx_0K_{\underline{a}}(1))$,
	\begin{enumerate} [\textup{(}{T2.}1\textup{)}] 
		\item\label{t_comp_1} if $\phi= 1\vee\abs{\cdot}^s$, $(\omega_k)_{k\in \N_0} = \left(k^{\sigma}\right)_{k\in \N_0}$, and $w= 1\vee\abs{\cdot}^{\mu}$, then \eqref{eq:weight2_res2} holds with $r(n_e)= \left(\ln\left(n_e\right)/\ln_2\left(n_e\right)\right)^{-\sigma}$,
		\item\label{t_comp_2} if $\phi= 1\vee\abs{\cdot}^s$, $ (\omega_k)_{k\in \N_0} = (e^{\kappa k \ln\left(k+1 \right)})_{k\in \N_0}$,  $s \geq \kappa(p+1)/(2k_q(\nu\indic\{W\neq W_{[-\underline{a},\underline{a}]}\}+1))$, $\Lambda: = (p-1)(1 - \left(\kappa(p+1)/\left(2s\left(k_q(\cdot+1)+\kappa\right)\right)\right)/2$,
		\begin{enumerate} [\textup{(}{T2.2}a\textup{)}] 
			\item\label{poly00}  
		and	$w^{I}(e^{2\kappa\abs{\cdot}\ln(\abs{\cdot}+1)}) = \cdot^{\nu}$, then (\ref{eq:weight2_res2}) holds with $r(n_e)= n_e^{-\kappa/(2\kappa +2(\nu+1)k_q)}\ln(n_{e})^{\Lambda(\nu)}$,
			\item\label{poly01} and  $\underline{a}$ such that $\mathbb{S}_{{\alpha}}\subseteq [-\underline{a},\underline{a}]$,  $w = W_{[-\underline{a},\underline{a}]}$, then (\ref{eq:weight2_res2}') holds with $r(n_e)= n_e^{-\kappa/(2\kappa+2k_q)}\ln(n_{e})^{\Lambda(0)}$,
		\end{enumerate}
		\item\label{t_comp_3} if $\phi= e^{\gamma \abs{\cdot}}$, $ r >1$,  $ (\omega_k)_{k\in \N_0} = (e^{\kappa \left(k \ln\left(k+1 \right)\right)^r})_{k\in \N_0}$, $w$ such that $w^{I}(e^{2\kappa\left(\abs{\cdot}\ln(\abs{\cdot}+1)\right)^{r}}) = \cdot^{\nu}$,  	
		$d_0 = 2\kappa\left(1 + (p-1)/(p+1)^r\right) +  4\kappa k_q(1+\nu)/((p+1)\ln(p+2))^{r-1}$, 	$k_0 := \lfloor r/(r-1)\rfloor$, and for $k\in \{1,\dots,k_0\}$, 
		\begin{align*} 
		d_k: =& \left(\frac{k_q(1+\nu)(2\kappa)^{1-1/r}\indic\{k\equiv  0  (\text{mod} \ 2)\}}{\kappa (1+1/((p+1)\ln(p+2)))^{r}}\right)^{k} \\
		&+\left(\left(\frac{(k_q+1)p+k_q-1}{p+1}+k_q\nu\right)\frac{\indic\{k\equiv  1  (\text{mod} \ 2)\} }{\kappa d_0^{1/r-1}}\right)^{k},
		\end{align*} 
	 and $\varphi :=  \exp(-\sum_{k=1}^{k_0} \left(-1\right)^k d_k \ln(\cdot)^{(1/r-1)k+1})\bigvee \ln(\cdot)^{p+1+(p-1)/r}$   then  (\ref{eq:weight2_res2}) holds with $r(n_e)=\sqrt{\varphi\left(n_e\right)/ n_e}$.
	\end{enumerate}
\end{theorem}

Theorem \ref{theo:lower} shows the rate in (T2.\ref{t_comp_2})  is optimal when $f_{\mt{X}}$ is known and $\mathbb{S}_{\mt{X}}=\mathcal{X}$. It is the same as in \cite{Meister07} for deconvolution with a known characteristic function of the noise on an interval when the signal has compact support. The rates in Theorem \ref{theo:compact} are independent of $p$ as common for severely ill-posed problems (see \cite{Chen07, Note}). The rates in (T2.\ref{t_comp_2}) and (T2.\ref{t_comp_3}) are polynomial and nearly parametric even if the problem is severely ill-posed.

\begin{theorem}\label{theo:non_compact} 
	Let  $ W=\cosh(\cdot/R )$. For all $q\in \{1,\infty\}$,  $ l,M ,s,R,\sigma, \kappa, \mu >0$, $\phi=1\vee\abs{\cdot}^s$, $\underline{N}$ solution of 
	$2k_q\underline{N}\ln\left(
	K_{\underline{a}}(e)\right) + \ln(\omega_{\underline{N}}^2)+(p-1)\ln(\underline{N})/q = \ln(n_e)$, and  $\epsilon=7e^2\pi/(2Rx_0
	K_{\underline{a}}(14e^2))$, 
	\begin{enumerate} [\textup{(}{T3.}1\textup{)}] 
		\item\label{t_noncomp_1} if  $ (\omega_k)_{k\in \N_0} = \left(k^{\sigma}\right)_{k\in \N_0}$, and $w= 1\vee\abs{\cdot}^{\mu}$, then \eqref{eq:weight2_res2} holds with $r(n_e)= (\ln\left(n_e\right)/\ln_2\left(n_e\right))^{-\sigma}$,  
		\item\label{t_noncomp_2} if $ (\omega_k)_{k\in \N_0} = (e^{\kappa k })_{k\in \N_0}$, $\underline{a}$ such that $\mathbb{S}_{{\alpha}}\subseteq [-\underline{a},\underline{a}]$,  and $w = W_{[-\underline{a},\underline{a}]}$,  then (\ref{eq:weight2_res2}') holds with $r(n_e)=n_e^{-\kappa/(2\kappa+ 2k_q)}\ln(n_{e})^{(p-1)\kappa/(2q(\kappa+k_q))}$.
	\end{enumerate}
\end{theorem}

In (T\ref{theo:non_compact}.\ref{t_noncomp_2}), we relax the assumption that $\mathbb{S}_{\mt{\beta}}$ is compact 
in (T\ref{theo:compact}.\ref{poly00}). 
The results of theorems \ref{theo:compact} and \ref{theo:non_compact} are related to those for ``2exp-severely ill-posed problems" (see \cite{cavalier2004block} and \cite{tsybakov2000best} which obtains the same 
rates up to logarithmic factor as in (T\ref{theo:compact}.\ref{poly01}) when $1/v(n_0,\mathcal{E}) \geq n$ and $p=1$). When $1/v(n_0,\mathcal{E}) \geq n  $, the rate in (T\ref{theo:non_compact}.2) matches  the lower bound for model \eqref{eq:modalt2}. 

\subsection{Data-driven estimator}\label{s4}
We use a Goldenshluger-Lepski method (see \cite{Gold_Lepski2,lacour2016minimal}). Let $R,\epsilon>0$, $q\in\{1,\infty\}$,
$ \zeta_0 = 1/(1+4p(1+\indic\{W=\cosh(\cdot/R)\}))$, 
$K_{{\rm max}} := \lfloor \zeta_0\ln(n)/\ln(2) \rfloor$, 
$T_{\max}:=2^{K_{{\rm max}}}$, $\mathcal{T}_n:= \{2^{k}:k = 1, \dots, K_{{\rm max}},2^{k}\ge\epsilon\}$, $p_n:=3\bigvee 6(1+\zeta_0)\ln(n)$, and, for all $N\in\N_0^{\R}$, $N_0,T\in\N_0$, $ t\ne 0 $, $N_{\max,q}^{W}  =\lfloor \underline{N_{\max,q}^W} \rfloor$,  $Q^W_{q}(N_0) :=  \indic\{q=\infty\} (2^p \indic\{W=W_{[-R,R]}\} + \indic\{W=\cosh(\cdot/R)\}) + (N_0+p-1)^{p-1}\indic\{q=1\}/(p-1)!$, 
\begin{align*}
B_1\left(t,N_0\right) &:= \underset{N_0\leq N'\leq N_{\max,q}^{W}}{ \max}\left( \sum_{N_0 \leq \abs{\mt{m}}_{q}\leq N'}\left(\frac{\abs{\widehat{c}_{\mt{m}}(t)}}{\sigma_{\mt{m}}^{W,x_0t}}\right)^2  - \Sigma\left(t, N'\right)\right)_{+}, \notag \\
B_2\left(T,N\right)& := \underset{ T'\in \mathcal{T}_n,T'\ge T}{\max}\left(\int_{T \leq \abs{t} \leq T'} \sum_{  \abs{\mt{m}}_{q}\leq N(t)}\left(\frac{\abs{\widehat{c}_{\mt{m}}(t)}}{\sigma_{\mt{m}}^{W,x_0t}}\right)^2  - \Sigma\left(t, N(t)\right)dt\right)_{+},\\
\Sigma\left(t,N_0\right) &:= 8(2+\sqrt{5})(1+2p_n) \frac{c_{\boldsymbol{X}}}{n} \left(\frac{\abs{t}}{2\pi}\right)^p  \nu_q^{W}(x_0t,N_0), \\
 \Sigma_2(T,N)&:=\int_{\epsilon \leq \abs{t} \leq T} \Sigma(t,N(t)) dt;
\end{align*} 
(N.1) if $W=W_{[-R,R]}$,  $\underline{N_{\max,q}^W}$ solution of $ 2k_q N_0\ln( 7e\pi N_0/(Rx_0 \epsilon)) = \ln(n) $,
\begin{align*}
	\nu_{q}^{W}(t,N_0)&= (N_0+1)^{k_q}Q^W_{q}(N_0)\left( 1\bigvee \frac{7e\pi (N_0+1)}{R\abs{t}} \right)^{2  k_qN_0}.
\end{align*} 
(N.2) 
if $W=\cosh(\cdot/R)$,
	\begin{align*}
	\underline{N_{\max,q}^{W}}= &\frac{\ln(n)}{2 k_q}  \indic\left\{\epsilon = \frac{\pi}{4Rx_0}\right\}  +\frac{\ln(n)}{2k_q\ln\left(7e^2\pi/(2Rx_0\epsilon)\right)}  \indic\left\{ \epsilon<  \frac{\pi}{4Rx_0}\right\},\\
	\nu_q^{W}(t,N_0) =  & 2^{k_q}\left( \frac{e\pi}{2} \right)^{2p} Q^W_{q}(N_0)\left(\frac{7e^2\pi}{2R|t|}\right)^{2k_qN_0} \indic\left\{\abs{t} \leq \frac{\pi}{4R} \right\} \\
	&  +2^p\left(\frac{ 2e R\abs{t}}{\pi}\right)^{k_q} Q^W_{q}(N_0)\exp\left( \frac{\pi k_q (N_0+k_q')}{2R\abs{t}} \right) \indic\left\{\abs{t} > \frac{\pi}{4R}\right\} ;
	\end{align*} 
$\widehat{N}$ and $\widehat{T}$ are defined, using $c_1 = 1+1/(2+\sqrt{5})^2$%
, as 
\begin{align}
\forall t\in \R\setminus(-\epsilon,\epsilon), \quad& \widehat{N}(t) \in   \underset{0\le N \leq N_{\max,q}^{W}}{\argmin}  \left( B_1(t,N) + c_1 \Sigma(t,N)\right),\label{eq:choiceN}\\
& \widehat{T} \in  \underset{T \in \mathcal{T}_n}{\argmin}  \left( B_2\left(T,\widehat{N}\right) +\Sigma_2\left(T,\widehat{N}\right) \right).\label{eq:choiceT}
\end{align}
The upper bounds below take the form
\begin{align}\label{eq:adapt}
\sup_{ \substack{f_{\alpha,\mt{\beta}} \in \mathcal{H}^{q,\phi,\omega}_{w,W}(l, M)\cap\mathcal{D} \\ f_{\mt{X}|\mathcal{X}} \in \mathcal{E}}} \frac{\mathcal{R}_{n_0}^{W}\left(\widehat{f}_{\alpha,\mt{\beta}}^{q,\widehat{N}, \widehat{T}, \epsilon} , f_{\alpha,\mt{\beta}} \right)}{r(n)^2} = \underset{\substack{v(n_0, \mathcal{E})/\delta(n_0) \leq n^{-2}\ln(n)^{-p}\\ n_e\ge 3}}{O_p}(1),
\end{align}
and we refer to (\ref{eq:adapt}') when we use the restriction $f_{\alpha,\mt{\beta}} \in \mathcal{H}^{q,\phi,\omega}_{w,W}(l)\cap\mathcal{D}$. 
\begin{theorem}\label{cor:adp_rate} 
	For all $l,M,s,R, \mu, \sigma>0$,  $H\in \N$,  $q\in \{1,\infty\}$, $\phi = 1\vee\abs{\cdot}^{s} $, if
	\begin{enumerate} [\textup{(}{T4.}1\textup{)}] 
		\item\label{a00} $(\omega_k)_{k\in\N_0}=(k^{\sigma})_{k\in\N_0}$, $  \underline{a} =1/\epsilon$, $w =1\vee\abs{\cdot}$, (\ref{eq:adapt}) holds with $r(n)= (\ln\left(n\right)/\ln_2\left(n\right))^{-\sigma}$ when
		\begin{enumerate} [\textup{(}{T4.1}a\textup{)}] 
			\item\label{a0} $W=W_{[-R,R]}$, $\mathbb{S}_{\mt{\beta}}\subseteq [-R,R]^p$,  and $\epsilon =  7e\pi/(Rx_0\ln(n))$,
			\item\label{b0} $W=\cosh(\cdot/R )$, and $\epsilon =  7e^2\pi/(2Rx_0\ln(n))$,
		\end{enumerate}
		\item\label{a01} $(\omega_k)_{k\in\N_0}=(e^{\kappa k \ln(1+ k)})_{k\in\N_0}$, $\underline{a}$ such that $\mathbb{S}_{{\alpha}}\subseteq [-\underline{a},\underline{a}]$, $w = W_{[-\underline{a},\underline{a}]}$, $W=W_{[-R,R]}$, $\mathbb{S}_{\mt{\beta}}\subseteq [-R,R]^p$, 
		$\epsilon =7e\pi/(Rx_0)$, and $s>(2p+1/2)\vee(\kappa(p+1)/(2k_q))$, (\ref{eq:adapt}') holds with $r(n)=n^{-\kappa/(2\kappa+2k_q)}\ln(n)^{1/2+ \Lambda(0)}$ and $\Lambda$ defined in (T2.\ref{t_comp_2}), 
		\item\label{a02} $(\omega_k)_{k\in\N_0}=(e^{\kappa k})_{k\in\N_0}$, $\underline{a}$ such that $\mathbb{S}_{{\alpha}}\subseteq [-\underline{a},\underline{a}]$, $w = W_{[-\underline{a},\underline{a}]}$, $W=\cosh(\cdot/R)$, $\epsilon =\pi/(4Rx_0)$, and $s> 4p+1/2$,  (\ref{eq:adapt}') holds with $r(n)=n^{-\kappa/(2\kappa+2k_q)}\ln(n)^{1/2+(p-1)\kappa/(2q(\kappa+k_q))}$.
	\end{enumerate}
\end{theorem}
The results in Theorem \ref{cor:adp_rate} are for $v(n_0, \mathcal{E})/\delta(n_0) \leq n^{-2}\ln(n)^{-p}$, in which case $n_e=n$. 
Theorem \ref{theo:lower} and (T\ref{cor:adp_rate}.\ref{a0})  show that $\widehat{f}_{\alpha,\mt{\beta}}^{q,\widehat{N}, \widehat{T}, \epsilon}$ is adaptive.  The rate in (T\ref{cor:adp_rate}.\ref{a01}) matches, up to a logarithmic factor, the lower bound in Theorem (T.\ref{theo:lower}.1\ref{it:expindic}) for model \eqref{eq:modalt2}. For the other cases, the risk is different for the lower bounds and the upper bounds in Theorem \ref{cor:adp_rate}, but using \eqref{eq:R_control} we obtain the same rate up to logarithmic factors for the risk involving the weight $W$. (T4.\ref{a01}) and (T4.\ref{a01}) rely on $\mathbb{S}_{{\alpha}}\subseteq [-\underline{a},\underline{a}]$ because, else, the choice $\underline{a} =w^{I}(\omega_{\underline{N}}^2) $ in Section \ref{s43} depends on the parameters of the smoothness class. However, it is possible to check that we can obtain the rate in (T2.\ref{poly00}) up to a $\sqrt{\ln(n)}$ factor when $\nu=1$ for a choice of $\underline{a}$ independent of the parameters of the smoothness class.

To gain insight, let us sketch the proof when $ \widehat{f}_{\mt{X}|\mathcal{X}}^{\delta}=f_{\mt{X}|\mathcal{X}}$ (hence we simply write $\mathcal{R}^{W}$).  
Let $T\in \mathcal{T}_{n}$, for all  $t\in \R$, $N\in\N_0^{\R}$, and $T'\in [0,\infty)$, 
\begin{align*}
 \mathcal{L}_{q}^{W}(t,N, T') &:= \left\|\left(\widehat{F}_1^{q,N,T',0} -  \mathcal{F}_{\rm{1st}}\left[f_{\alpha,\mt{\beta}}\right]\right)(t,\cdot_2) \right\|^2_{L^2\left(W^{\otimes p}\right)},   
\end{align*} and $\widetilde{w}:= \indic\{w\neq W_{[-\underline{a},\underline{a}]}\}/w$. The Plancherel identity and \eqref{eq:K1K2} yield
\begin{align}
\mathcal{R}^{W}\left(\widehat{f}_{\alpha,\mt{\beta}}^{q,\widehat{N},\widehat{T},\epsilon},f_{\alpha,\mt{\beta}}\right)
&\le \frac{C}{2\pi} \int_{\epsilon \leq \abs{t}}
\E \left[ \mathcal{L}_{q}^{W}\left(t,\widehat{N}(t), \widehat{T}\right)   \right]dt +  C M^2 \widetilde{w}(\underline{a}).\label{eRisk}
\end{align}
The upper bound in \eqref{eRisk} with nonrandom $\widehat{N}$ and $\widehat{T}$ is the one we use to obtain theorems \ref{theo:compact} and \ref{theo:non_compact}. \eqref{eq:choiceT}  allows to obtain an upper bound with a similar quantity but with arbitrary nonrandom $\widehat{T}$.  By arguments in the proof of Lemma \ref{adp2} for the first inequality and \eqref{eq:choiceT} for the second,  
\begin{align*}
&\int_{\epsilon \leq \abs{t}} \E \left[ \mathcal{L}_{q}^{W}\left(t,\widehat{N}(t),\widehat{T}\right)   \right]dt\\
&\leq  \frac{2 + \sqrt{5}}{\sqrt{5}}\int_{\epsilon \leq \abs{t}}
\E \left[ \mathcal{L}_{q}^{W}\left(t,\widehat{N}(t),T\right)   \right]dt \\
& \quad + (2+\sqrt{5}) \left( \E\left[B_2\left(\widehat{T}, \widehat{N}\right)  + \Sigma_2\left(T,\widehat{N}\right)  \right]+ \E\left[B_2\left(T, \widehat{N}\right)  + \Sigma_2\left(\widehat{T},\widehat{N}\right)  \right]\right)  \\
&\leq  \frac{2 + \sqrt{5}}{\sqrt{5}}\int_{ \epsilon \leq \abs{t}}
\E \left[ \mathcal{L}_{q}^{W}\left(t,\widehat{N}(t),T\right)   \right]dt + 2(2+\sqrt{5}) \E\left[B_2\left(T, \widehat{N}\right)  + \Sigma_2\left(T,\widehat{N}\right)  \right]
\end{align*}
and, by a 
Talagrand's 
inequality, 
$$\E\left[B_2\left(T, \widehat{N}\right)\right] \leq \left(1+\frac{2}{\sqrt{5}}\right)\int_{ \epsilon \leq \abs{t}}
\E \left[ \mathcal{L}_{q}^{W}(t,\widehat{N}(t),T)   \right]dt + O\left(\frac{1}{n}\right),$$
where the $O\left(1/n\right)$ term is independent of $T$ and $\widehat{N}$. \eqref{eq:choiceN} allows to obtain yet another upper bound which replaces $\widehat{N}$ by an arbitrary nonrandom $N$. We conclude because the final upper bound \eqref{eq:start_adpt0} has a similar form as \eqref{eq:START} 
when we deal with nonrandom $\widehat{N}$ and $\widehat{T}$. 

\section{Simulations}\label{s3}
Let $p=1$, $q=\infty$, and $  (\alpha,\beta)^{\top} =\xi_1D+\xi_2(1-D)$ with $\mathbb{P}(D=1) = \mathbb{P}(D=0) =0.5$. The law of $X$ is a truncated normal based on a normal of mean $ 0 $ and variance $ 2.5$ and truncated to $\mathcal{X}$ with $x_0=1.5$. The laws of $\xi_1$ and $\xi_2$ are either: (Case 1) truncated normals based on
normals with means $ \mu_1 =
\left(
\begin{array}{c}
-2  \\
3 
\end{array}
\right)$ and $ \mu_2 = \left(
\begin{array}{c}
3  \\
0 
\end{array}
\right)$, same covariance
$\left(\begin{array}{cc}
2 & 1 \\ 
1 & 2
\end{array} \right)$, and truncated to $[-6,6]^{p+1}$ or (Case 2) nontruncated. 

\begin{figure}[H]
	\centering
	\subfigure[True density]{\includegraphics[width=0.4\linewidth, height=0.14\textheight]{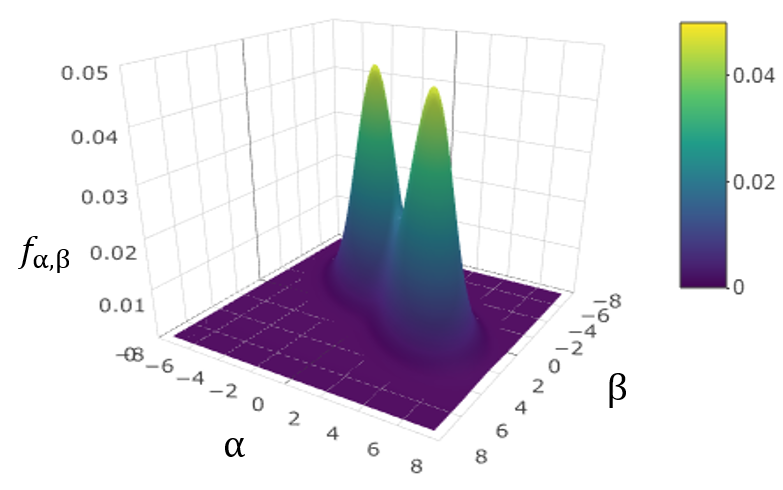}
		\label{fig:subfigure1}}
	\quad
	\subfigure[Mean of estimates]{\includegraphics[width=0.4\linewidth, height=0.14\textheight]{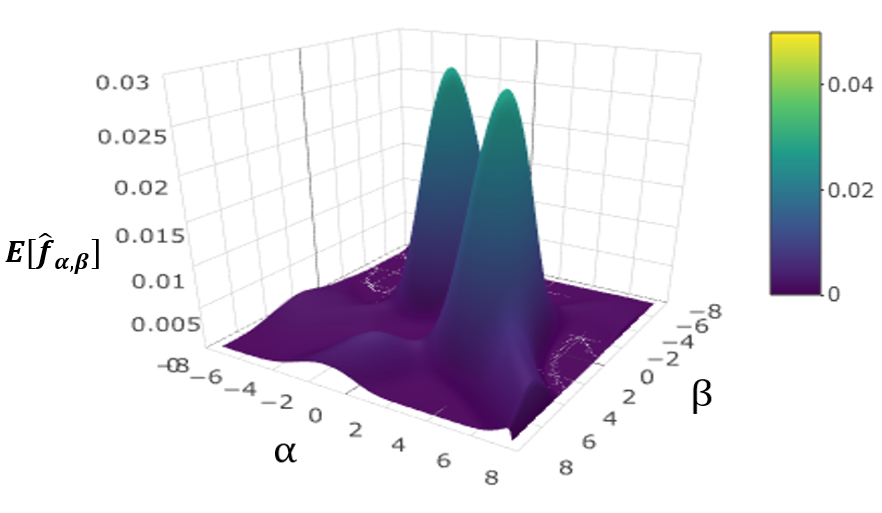}
		\label{fig:subfigure2}}
	\subfigure[97.5\% quantile of estimates]{\includegraphics[width=0.4\linewidth, height=0.14\textheight]{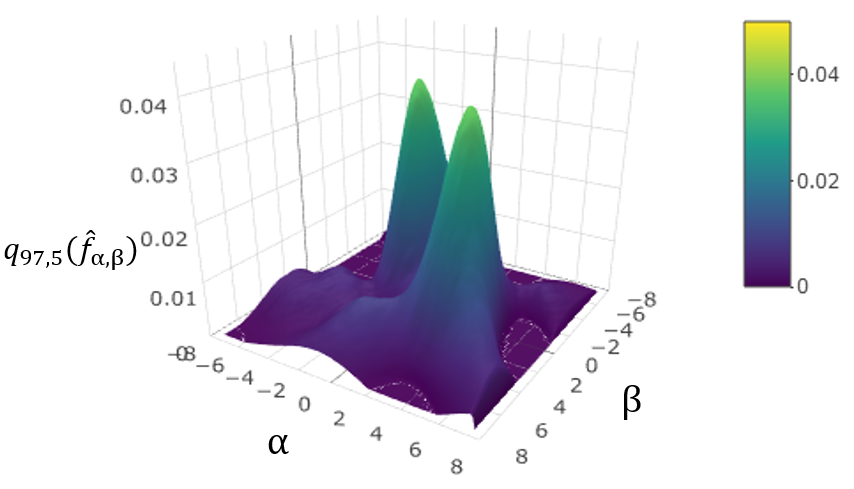}
		\label{fig:subfigure3}}
	\quad
	\subfigure[2.5\% quantile of estimates]{\includegraphics[width=0.4\linewidth, height=0.14\textheight]{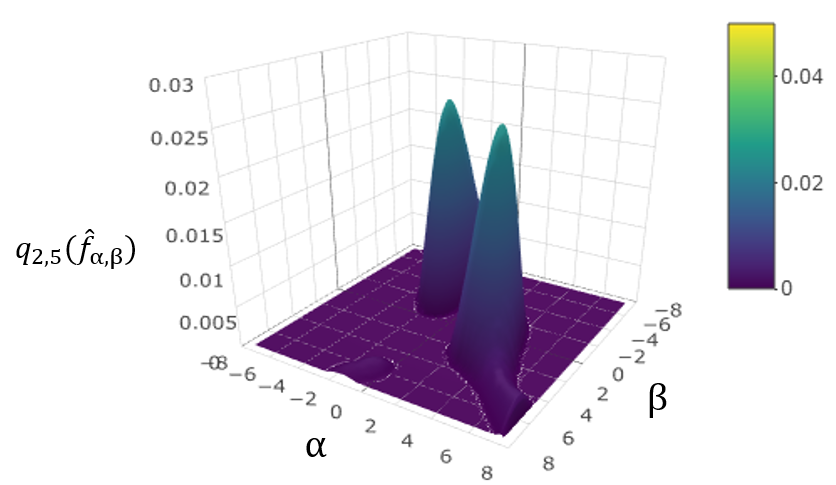}
		\label{fig:subfigur&e4}}
	\caption{Case 1, $W = W_{[-7.5,7.5]}$}
	\label{fig:Sp2}
	\end{figure}

\begin{figure}[H]
	\centering
	\subfigure[True density]{\includegraphics[width=0.4\linewidth, height=0.14\textheight]{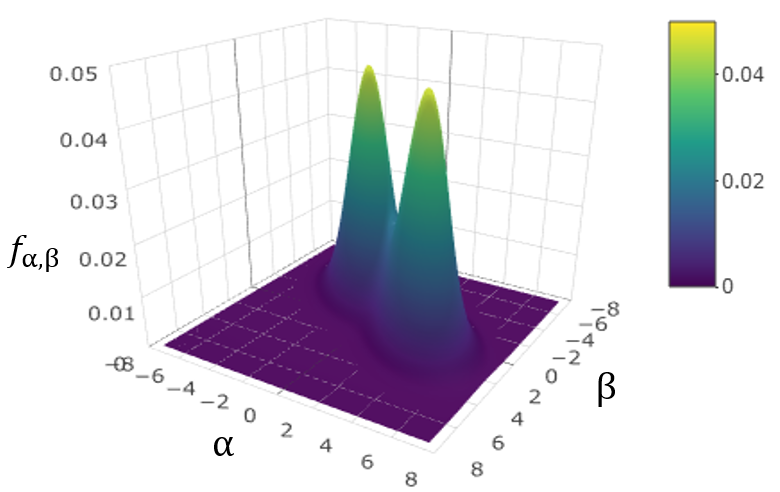}
		\label{fig:subfigure11}}
	\quad
	\subfigure[Mean of estimates]{\includegraphics[width=0.4\linewidth, height=0.14\textheight]{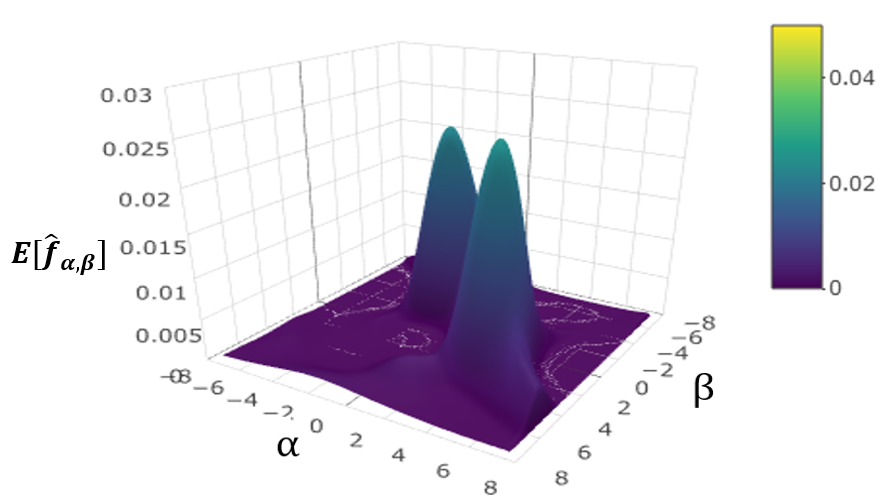}
		\label{fig:subfigure22}}
	\subfigure[97.5\% quantile of estimates]{\includegraphics[width=0.4\linewidth, height=0.14\textheight]{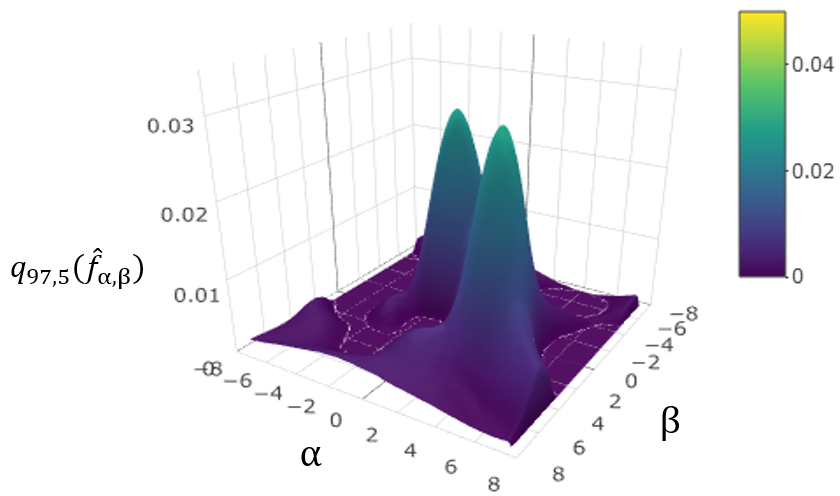}
		\label{fig:subfigure33}}
	\quad
	\subfigure[2.5\% quantile of estimates]{\includegraphics[width=0.4\linewidth, height=0.14\textheight]{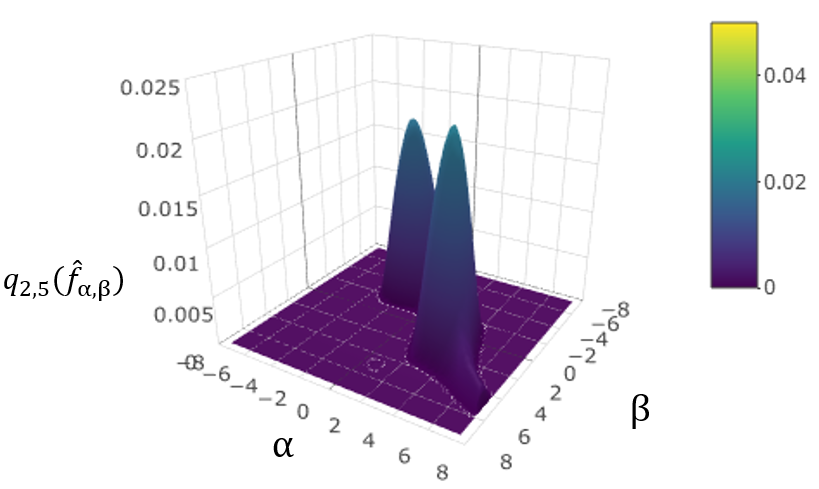}
		\label{fig:subfigure44}}
	\caption{Case 2, $W=\cosh\left(\cdot/7.5\right) $}
	\label{fig:Sp3}
\end{figure}

Table \ref{fig:table3} compares $\E [\| \widehat{f}_{\alpha,\beta}^{\infty,\widehat{N},\widehat{T},\epsilon} - f_{\alpha,\beta} \|^2_{L^2([-7.5,7.5]^2) }]$ and the risk of the oracle $\underset{T \in \mathcal{T}_n,N\in\mathcal{N}_{n,H} }{\min}    \E [\| \widehat{f}_{\alpha,\beta}^{\infty,N,T,\epsilon} - f_{\alpha,\beta} \|^2_{L^2([-7.5,7.5]^2)}]$ for cases 1 and 2. The Monte-Carlo experiment uses 1000 simulations. 
\begin{table}[H]
	\centering \footnotesize
	\begin{tabular}{@{}rrrrcrrr@{}}\toprule	
  	& \multicolumn{3}{c}{$W = W_{[-7.5,7.5]}$, Case 1} &  \phantom{abc} & \multicolumn{3}{c}{ $W=\cosh\left(\cdot/7.5\right) $, Case 2} \\
		MISE 	&$n=300$  	& $n=500$  & $n=1000$  & &$n=300$  & $n=500$ &$n=1000$  \\ \midrule
	data-driven &   $0.092$ &  $ 0.086$ & 0.083 & &  $   	0.089   $ &  0.087   &0.085 \\
	oracle & $ 0.091 $   &  $0.086 $  & 0.082  & & $ 0.088 $  &  0.087 & 0.085  \\
		\bottomrule 
	\end{tabular}
	\caption{Risk}\label{fig:table3}
\end{table}
Figure \ref{fig:Sp2} (resp. Figure \ref{fig:Sp3}) displays summaries of the law of the estimator  for $W = W_{[-7.5,7.5]}$ (resp. $W =\cosh(\cdot/7.5)$)  in Case 1 (resp. Case 2) and $n=1000$. As standard in the literature (see, \emph{e.g.}, \cite{comte2013nonparametric, dion2014new}), 
the multiplicative constant appearing in $\Sigma$ 
 is in practice calibrated from a simulation study. 
$\widehat{f}_{X|X\in\mathcal{X}}$ is obtained with the same data to illustrate that sample splitting is unnecessary in practice and only used for the mathematical argument. For $\widehat{f}_{X|X\in\mathcal{X}}$ we use a Gaussian kernel density estimator using the \textsf{R} package ks and the multivariate plug-in bandwidth selector of \cite{wand1994multivariate}.  $\epsilon$ is chosen as in (T4.\ref{a0}) and (T4.\ref{b0}) respectively for Case 1 and Case 2. The estimator requires the SVD of $\mathcal{F}_c$.  By Proposition \ref{pdebut}, we have $g_m^{W(\cdot/R),c}=g_m^{W,Rc}$ for all $m\in\N_0$. When $W=W_{[-1,1]}$, the first coefficients of the decomposition on the Legendre polynomials are obtained by solving for the eigenvectors of two tridiagonal symmetric Toeplitz matrices (see Section 2.6 in \cite{Osipov}). When $W=\mathrm{cosh}$, we refer to Section 7 in \cite{Note}.   We use the image of $g_m^{W,Rc}$ by the adjoint of $\mathcal{F}_c$ (see Appendix \ref{subsec:estim}) and 
that $\varphi_m^{W,Rc}$ has norm 1 to get the rest of the SVD. We obtain the Fourier inverse by fast Fourier transform. We use a resolution of $2^{13}$, which appears on simulations to realise a good trade-off between computational time and precision. For more details about the implementation, see the vignette \cite{randomcoefficients} of the package \href{https://CRAN.R-project.org/package=RandomCoefficients}{RandomCoefficients}.%

\appendix

\renewcommand{\theequation}{A.\arabic{equation}}
\renewcommand{\thelemma}{A.\arabic{lemma}}
\renewcommand{\thecorollary}{A.\arabic{corollary}}
\renewcommand{\thedefinition}{A.\arabic{definition}}
\renewcommand{\theproposition}{A.\arabic{proposition}}
\renewcommand{\theremark}{A.\arabic{remark}}
\renewcommand{\thetheorem}{A.\arabic{theorem}}
\renewcommand{\theassumption}{A.\arabic{assumptio}}
\renewcommand{\thesubsection}{A.\arabic{subsection}}
 \setcounter{equation}{0}  
 \setcounter{lemma}{0}
 \setcounter{corollary}{0}
 \setcounter{proposition}{0}
 \setcounter{remark}{0}
 \setcounter{definition}{0}
 \setcounter{lemma}{0}
 \setcounter{theorem}{0}
 \setcounter{assumption}{0}
 \setcounter{section}{0}
  \setcounter{subsection}{0}
 \setcounter{footnote}{0}
 \setcounter{figure}{0}

\section{Proofs of the main results}\label{sec:estim}

\subsection{Notations and preliminaries}\label{subsec:estim}
$\mathfrak{R}$ and $\mathfrak{I}$ denote the real and imaginary parts. For a differentiable function $f$ of real variables, $f^{(\mt{m})}$ denotes $\prod_{j=1}^d\frac{\partial^{\mt{m}_j}}{\partial x_j^{\mt{m}_j}}f$.  $C^{\infty}\left(\R^d\right) $ is the space of infinitely differentiable functions. 
Abusing notations, we sometimes use $\mathcal{F}_{c}[f]$ for the function in $L^2(\R)$ and $\mathcal{E}xt[f]$  assigns the value 0 outside $[-1,1]^d$. Denote by  $\Pi:L^2(\R^d)\to L^2(\R^d)$ such that $\Pi f(\mt{x})=f(-\mt{x})$ and by \begin{equation}\label{eqdef} \begin{array}{cccc}
\mathcal{C}_{c}: & L^{2}\left(\R^{d}\right) & \rightarrow &  L^{2}\left(\R^{d}\right)  \\
	      &      f		      & \rightarrow &   |c|^df(c\ \cdot).
\end{array} 
\end{equation}
 Because $\mathcal{F}_c=\mathcal{F}\mathcal{C}_{1/c}=(1/|c|)\mathcal{C}_{c}\mathcal{F}$, $\Pi\mathcal{F}_c=\mathcal{F}_c\Pi$, $\mathcal{F}_c^*=(1/W)\Pi\mathcal{F}_c\mathcal{E}xt$, and $W$ is even, we obtain $\mathcal{F}_c^*=\Pi\left((1/W)\mathcal{F}_c\mathcal{E}xt\right)$ and  
 \begin{align*}
 	\mathcal{F}_c\mathcal{F}_c^*&=\Pi\mathcal{F}_c\left(\frac{1}{W}\mathcal{F}_c\mathcal{E}xt\right)\\
 	&=\frac{2\pi}{|c|}\mathcal{F}^I\left(\mathcal{C}_{1/c}\left(\frac{1}{W}\mathcal{C}_c\mathcal{F}\mathcal{E}xt\right)\right)\\
&=2\pi\mathcal{F}^I\left(\mathcal{C}_{1/c}\left(\frac{1}{W}\right)\mathcal{F}\mathcal{E}xt\right).
\end{align*} 
The operator $\mathcal{Q}_c^W=(\abs{c}/(2\pi))\mathcal{F}_{c}\mathcal{F}_{c}^*$ is a compact positive definite self-adjoint operator (see \cite{Osipov} and \cite{Widom} for the two choices of $W$). Its eigenvalues in decreasing order repeated according to multiplicity are denoted by $\left(\rho_{m}^{W,c}\right)_{m\in\N_0}$ and a basis of eigenfunctions by $\left(g_{m}^{W,c}\right)_{m\in\N_0}$. The other elements of the SVD are $\sigma_{{m}}^{W,c}  =\sqrt{2\pi\rho_{{m}}^{W,c}/\abs{c}} $ and $\varphi_{{m}}^{W,c}=\mathcal{F}_c^*g_{m}^{W,c}/\sigma_{{m}}^{W,c}$. We denote, for all $m\in\N_0$, 
by $\psi_m^{c}$ the function $g_{m}^{W_{[-1,1]},c}$ and $ \mu_{m}^{c} =  i^{m} \sigma_{m}^{W_{[-1,1]},c}$. Because $\psi_m^{c}=\mathcal{F}_c( \mathcal{E}xt[\psi_m^c])/\mu_m^{c}$ in $L^2([-1,1])$, $\psi_m^{c}$ can be extended as an entire function which we denote with the same notation.  
Using the injectivity of $\mathcal{F}_c$ (see the proof of Proposition \ref{prop:notcompact}), we have  $\varphi_{m}^{W_{[-1,1]},c}= i^{-m}  \mathcal{E}xt[\psi_m^{c}]$.  \\
We make use of
\begin{align}\label{elog}\forall a,b>0,\ \sup_{t\ge1} \frac{\ln(t)^a}{t^b} = \left(\frac{a}{eb}\right)^a,\\
\label{Young}\forall c>0,\ \forall a,b\in \R,\ ab\leq \frac{a^2}{2c} + \frac{b^2c}{2}. 
\end{align} 
Expectations are conditional on $\mathcal{G}_{n_0}$ when $f_{\mt{X}|\mathcal{X}}$ is unknown and we rely on $\mathcal{G}_{n_0}$ to estimate it. We remove the conditioning in the notations for simplicity.  

\subsection{Proofs of Proposition \ref{prop:notcompact}, \ref{sec:upper:extension}, and \ref{prop:interp}}$\ $\\
\noindent {\bf Proof of Proposition \ref{prop:notcompact}.} 
The first assertion comes from the fact that $W$ is bounded from below. 
The second uses Theorem IX.13 in \cite{ReedSimon} which implies that, for all $c\ne0$, $\mathcal{F}_c$ defined in \eqref{eqdef} is injective. 
We now show that $\mathcal{K}$ is continuous at 0. Let $f\in L^2\left(w\otimes W^{\otimes p}\right)$. The change of variables, the Plancherel identity, and the lower bounds on the weights yield
\begin{align*}
\left\|\mathcal{K}[f]\right\|_{L^{2}(\R\times[-1,1]^p)}^2&\le\int_{\R^{p+1}}\left|\mathcal{F}[f](t,\mt{v})\right|^2(t,\mt{v})dtd\mt{v}
\le\left(\frac{2\pi}{W(0)}\right)^{p}\left\|f\right\|_{
	L^2\left(w\otimes W^{\otimes p}\right)}^2.
\end{align*}
Let $w=1$. We exhibit a bounded sequence $ \left( f_{k}\right)_{k\in\N_0} $ in $ L^2(1\otimes W^{\otimes p}) $ for which 
there does not exist a convergent subsequence of $ \left( \mathcal{K}\left[ f_{k}\right]\right)_{k\in\N_0} $. 
Take $v_0$ such that $\mathrm{supp}(v_0)\subset  [1, 2]$, $ \norm{v_0}_{L^2(\R)}=1 $ and, for all $  k \in\N_0 $ and $(a,\mt{b}^{\top})^{\top}\in\R^{p+1}$, $v_k(\cdot)= 2^{-k/2}v_0(2^{-k}\cdot) $ and $f_{k}(a,\mt{b})= \mathcal{F}^I\left[ v_k(\cdot) \varphi_{\underline{\mt{0}}}^{W,x_0\cdot }(\mt{b}) \right](a)$. $\left( f_{k}\right)_{k\in\N_0}$ is bounded by the Plancherel identity and 
\begin{align*}
\norm{f_{k}}^2_{L^2\left(1\otimes W^{\otimes p}\right)}  &= \frac{1}{2\pi} \int_{\R} v_k(t)^2 \int_{\R^p} \left|\varphi_{\underline{\mt{0}}}^{W,x_0t}(\mt{b}) \right|^2 W^{\otimes p}(\mt{b}) dtd\mt{b}\leq \frac{1}{2\pi}.
\end{align*}
Using $ \mathcal{K}\left[f_{k}\right](\cdot,\cdot_2) = \sigma_{\underline{\mt{0}}}^{W,x_0\cdot}  v_k(\cdot) g_{\underline{\mt{0}}}^{W,x_0\cdot }(\cdot_2) \abs{x_0\cdot}^{p/2} $,
$ c\in (0,\infty) \mapsto \rho_{0}^{W,c} $ is nondecreasing (by Lemma 1 in \cite{Note} which holds for all $W$ satisfying (H1.\ref{E2})), and using, for all $j\in\N_0$, $\norm{v_j}_{L^2(\R)}=1$, we obtain, for all $(j,k) \in \N_0^2$, $j< k$, 
\begin{align*}
\left\| \mathcal{K}\left[ f_{j}\right] - \mathcal{K}\left[f_{k}\right]\right\|^2_{L^2(\R \times[-1,1]^p)} 
&\geq \rho_{\underline{\mt{0}}}^{W,2^jx_0}(2\pi)^p
\int_{\R} \left(v_j(t)^2+v_k(t)^2\right)  dt \\
&\geq 2(2\pi)^p\rho_{\underline{\mt{0}}}^{W,x_0}>0,
\end{align*}
so $\mathcal{K}$ is not compact. 
\hfill $\square$\vspace{0.3cm}


\noindent {\bf Proof of Proposition  \ref{sec:upper:extension}.} This holds by Theorem 15.16 in \cite{Kress} and the injectivity of $\mathcal{F}_c$.\hfill $\square$\vspace{0.3cm}

\noindent {\bf Proof of Proposition  \ref{prop:interp}.} 
Take $f\in L^2(\R)$ and start by showing that $ \mathcal{I}_{\underline{a},\epsilon}[f]\in L^2([-\epsilon,\epsilon])$. The terms $1 - \rho_{m} ^{W_{[-1,1]},\underline{a}\epsilon}$ in the denominator of \eqref{eq:inter}  are nonzero because $ \left(\rho^{W_{[-1,1]},\underline{a}\epsilon}_m\right)_{m\in\N_0} $ is nonincreasing and $ \rho^{W_{[-1,1]},\underline{a}\epsilon}_0 <1$ (see (3.49) in \cite{Osipov}). 
Using that $\left(g_{m}^{W_{[-1,1]},\underline{a}\epsilon}\left(\cdot/\epsilon\right)/\sqrt{\epsilon}\right)_{m\in\N_0}$  is a basis of $L^2([-\epsilon,\epsilon])$, that $ \left(\rho^{W_{[-1,1]},\underline{a}\epsilon}_m\right)_{m\in\N_0} $ is nonincreasing, and the Cauchy-Schwarz inequality for the first inequality, using $ \sum_{m\in\N_0} \rho_m^{W_{[-1,1]},\underline{a}\epsilon} = 2\underline{a}\epsilon/\pi  $ (see  (3.55) in \cite{Osipov}) and $\left\|g_{m}^{W_{[-1,1]},\underline{a}\epsilon}\right\|^2_{L^2(\R)}=1/ \rho_m^{W_{[-1,1]},\underline{a}\epsilon}$ (see (3) in \cite{PSWF}) for the second yield
\begin{align}
&\sum_{m \in \N_0} \left(\frac{\rho_{m} ^{W_{[-1,1]},\underline{a}\epsilon}}{\left(1 - \rho_{m} ^{W_{[-1,1]},\underline{a}\epsilon}\right)\epsilon}\right)^2 \left|\braket{f,g^{W_{[-1,1]},\underline{a}\epsilon}_m\left(\frac{\star}{\epsilon}\right) }_{L^{2}(\R\setminus[-\epsilon,\epsilon])}\right|^2 \left\|g_{m}^{W_{[-1,1]},\underline{a}}\left(\frac{\cdot}{\epsilon}\right)\right\|^2_{L^2([-\epsilon,\epsilon])}\notag\\
&\leq  \frac{\left\|f\right\|_{L^{2}(\R\setminus[-\epsilon,\epsilon])}^2  }{\left(1 - \rho_0^{W_{[-1,1]},\underline{a}\epsilon}\right)^2} \sum_{m \in \N_0} \left(\rho_{m}^{W_{[-1,1]},\underline{a}\epsilon}\right)^2  \left\|g_{m}^{W_{[-1,1]},\underline{a}\epsilon} \right\|^2_{L^{2}(\R)} \leq  \frac{2\underline{a}\epsilon \left\|f\right\|_{L^{2}(\R\setminus[-\epsilon,\epsilon])}^2}{\pi\left(1 - \rho_0^{W_{[-1,1]},\underline{a}\epsilon}\right)^2}\label{e:maj}.
\end{align}
Let us now show the second statement. Take $ \epsilon >0 $ and $ g\in PW(\underline{a}) $. Let $(\alpha_m)_{m\in\N}$ be the sequence of coefficients of $g(\epsilon\cdot)\in PW(\underline{a}\epsilon)$ on the complete orthogonal system $\left(g_{m}^{W_{[-1,1]},\underline{a}\epsilon}\right)_{m\in\N_0}$. Because $\left(g_{m}^{W_{[-1,1]},\underline{a}\epsilon}\right)_{m\in\N_0}$ is a basis of $L^2([-1,1])$, we have 
$ \sum_{m \in \N_0} \alpha_m g^{W_{[-1,1]},\underline{a}\epsilon}_m   =  g(\epsilon\cdot)  \indic\{\abs{\cdot}\geq 1\} + \sum_{m \in \N_0}\alpha_m  g^{W_{[-1,1]},\underline{a}\epsilon}_m \indic\{\abs{\cdot}\leq 1\}.$ 
We identify the coefficients by taking the Hermitian product in $ L^{2}(\R) $ with $  g^{W_{[-1,1]},\underline{a}\epsilon}_m$ and obtain $\mathcal{I}_{\underline{a},\epsilon}[g]=g$ in $L^2(\R)$ and, for all $f,h\in L^2(\R)$, 
\begin{align}
\left\|f-\mathcal{I}_{\underline{a},\epsilon}\left[ h\right]\right\|_{L^2([-\epsilon,\epsilon])}^2
&
\le 2\left(\left\|f-\mathcal{P}_{\underline{a}}[f]\right\|_{L^2([-\epsilon,\epsilon])}^2+\left\|\mathcal{I}_{\underline{a},\epsilon}\left[ \mathcal{P}_{\underline{a}}\left[f\right]-  h\right]\right\|_{L^2([-\epsilon,\epsilon])}^2 \right). \label{eq:interp1}
\end{align}
Replacing $f$ by $\mathcal{P}_{\underline{a}}\left[f\right]-  h$ in \eqref{e:maj} yields
\begin{equation}\label{eq:interp2}
\left\|\mathcal{I}_{\underline{a},\epsilon}\left[ \mathcal{P}_{\underline{a}}\left[f\right]-  h\right]\right\|_{L^2([-\epsilon,\epsilon])}^2 \le  \frac{C_0(\underline{a}\epsilon)}{2} \left\|\mathcal{P}_{\underline{a}}\left[f\right]-  h\right\|_{L^2(\R\setminus[-\epsilon,\epsilon])}^2 .
\end{equation}
Using \eqref{eq:interp1} and \eqref{eq:interp2} for the first display, $\mathcal{P}_{\underline{a}}\left[f\right]-  h = \left(\mathcal{P}_{\underline{a}}[f] - f\right) + \left(f-h\right) $ and the Jensen inequality for the second display, we obtain
\begin{align*}
&\left\|f-\mathcal{I}_{\underline{a},\epsilon}\left[h\right]\right\|_{L^2([-\epsilon,\epsilon])}^2\\
 &\le 2\left\|f-\mathcal{P}_{\underline{a}}[f]\right\|_{L^2([-\epsilon,\epsilon])}^2+C_0(\underline{a}\epsilon) \left\|\mathcal{P}_{\underline{a}}[f] -h \right\|_{L^2(\R\setminus[-\epsilon,\epsilon])}^2 \\
&\le 2(1 + C_0(\underline{a}\epsilon) )\left\|f-\mathcal{P}_{\underline{a}}[f]\right\|_{L^2(\R)}^2+ 2C_0(\underline{a}\epsilon) \left\|f-h\right\|_{L^2(\R\setminus[-\epsilon,\epsilon])}^2.\quad\square
\end{align*}

\subsection{Lower bounds}\label{app:lower} We denote by $\mathbb{P}_{j}$ the law of density $f_{j,n}$ and by $\mathbb{P}_{j,n}$ the law of an iid sample of size $n$,  and use
\begin{align*} \underset{\widehat{f}}{\inf} \ \  \underset{f\in   \mathcal{H}}{\sup} \E \left[\left\|\widehat{f} -f\right\|_{L^2(\R^{p+1})}^2  \right] \geq  \underset{\widehat{f}}{\inf} \ \   \max_{f_{j,n}  \in  \mathcal{H}, \ j\in\{1,2\}} \E \left[ \left\| \widehat{f} -  f_{j,n}\right\|_{L^2(\R^{p+1})}^2 \right]
\end{align*}
and the next lemma (see Theorem 2.2, (2.5), and (2.9) in \cite{Tsybakov}). 
\begin{lemma}\label{lem:lower}
	If there exists $\xi<\sqrt{2}$ such that
	\begin{enumerate}[\textup{(}i\textup{)}] 
		\item \label{i} $ \forall  j \in\{1, 2\}, \ f_{j,n} \in   \mathcal{H}$,
		\item \label{ii}$\norm{f_{1,n} - f_{2,n}}_{L^2(\R^{p+1})}^2 \geq 4r(n)^2 >0 $,
		\item \label{iii}$\chi_2(\mathbb{P}_{2, n}, \mathbb{P}_{1,n}) \leq \xi^2$ or $K(\mathbb{P}_{2, n},\mathbb{P}_{1,n}) \leq \xi^2$,
	\end{enumerate}
	then we have 
	$$\frac{1}{r(n)^{2}}\underset{\widehat{f}}{\inf} \ \  \underset{  f_{j,n} \in   \mathcal{H},\ j\in\{1,2\}}{\max}  \E \left[\left\|\widehat{f} -  f_{j,n}\right\|_{L^2(\R^{p+1})}^2  \right]   \geq \frac{1}{2}\left( \frac{e^{-\xi^2}}{2}\bigvee \left(1-\frac{\xi}{\sqrt{2}}\right)\right).$$
\end{lemma}
The proofs of the lower bounds consist in defining $f_{j,n} \in   \mathcal{H}$, for  $j \in\{1, 2\}$ as a function of parameters, establishing the conditions on the parameters so that $f_{j,n}$ satisfy the three conditions \eqref{i}-\eqref{iii} of Lemma \ref{lem:lower}, and finally choosing the value of these parameters as a function of $n$ to deduce the lower bounds.\\
\noindent {\bf Proof of (T\ref{theo:lower}.1\ref{it:lowerindic}).} For $j=1,2$, $ f_{j,n}$ is a possible $f_{\alpha,\mt{\beta}}$, $\left(b_{\boldsymbol{m}}^j\right)_{\mt{m}\in\N_0^p}$ the sequence of its coefficients (see Section \ref{s22}), 
Steps 1-3 give conditions under which \eqref{i}-\eqref{iii} in Lemma \ref{lem:lower} are satisfied when $f_{1,n} := f_0 $ and 
\begin{align}
&
f_{2,n} := f_0+F,\ f_0(a,\mt{b}) :=  \frac{1}{\pi \tau\left(1+\left(a/\tau\right)^2\right)} \frac{ \indic\{\abs{\mt{b}}_{\infty} \leq R\} }{(2R)^{p/2} },
\label{eq:hyp}
\end{align}
for all $(a,\mt{b})\in\R^{p+1}$,
\begin{align}
& F(a,\mt{b}): = \gamma \mathcal{F}_{1\mathrm{st}}^{I}\left[\left( \frac{c(\abs{\star})}{2\pi} \right)^{p/2} \lambda(\star)\psi_{\widetilde{\mt{N}}(q)}^{Rc(\star)}\left( \frac{\mt{b}}{R}\right) \right](a) \indic\{\abs{\mt{b}}_{\infty} \leq R\},\label{eq:H_N}
\end{align}
for all $U/2\le|t|\le U$
\begin{align}
&\lambda(t) := 
\exp\left( 1 - \frac{1}{1 - \left( 4\abs{t} - 3U\right)^2/U^2}\right),\ \textrm{else }\lambda(t):=0\label{eq:phi},
\\
&\widetilde{\mt{N}}(1):=\left(N,\underline{\mt{N(Rx_0U)}}^{\top}\right)^{\top},\ \widetilde{\mt{N}}(\infty):=\underline{\mt{N}}\in\N^p,\label{eq:Nq}
\end{align}
$N(Rx_0U) := \lceil H(Rx_0U)\rceil$, for $H$ from Section \ref{app:lower1}, $ n $ large enough, $ N$ (odd), $ \gamma  $, $\tau\ge 1$, and $U$ from Step 4 and such that $N\ge N(Rx_0U)$, hence $N \geq Rx_0U\vee 2$ by the discussion before Lemma \ref{lem:ratio}. Note $\left\|\lambda\right\|_{L^{\infty}(\R)} \leq 1$.\\
\textbf{Step 1.1.} We prove that $ f_{1,n} $ and $ f_{2,n} $ are nonnegative when 
\begin{align}
\gamma UN^{k_q/2}\frac{\left((Rx_0U)/\pi\right)^{\frac{p}{2}} }{1+p/2}   \left(\frac{5}{4}\right)^{\frac{k_q}{2}} \left(\frac{5}{4}N(Rx_0U)\right)^{\frac{p-1}{2q}}\leq    \frac{1}{\tau+1/\tau},\label{eq:pos}  \\
\gamma UN^{k_q/2}\frac{2^{\frac{p}{2}}}{2}C_{8}(Rx_0U, p, U)N^{2}\leq    \frac{1}{\tau+1/\tau},\label{eq:pos1}
\end{align}
where 
$C_{8}(Rx_0U, p, U)$ 
is defined in  Lemma \ref{lem:secondH}. Let 
$ (a,\mt{b}) \in \R \times[-R,R]^p $. 
We show that 
\eqref{eq:pos} and \eqref{eq:pos1} yield $ f_0(a,\mt{b}) \geq  \abs{F(a,\mt{b})} $ which ensures that $  f_{2,n}$ is nonnegative. 
\eqref{eq:pos} yields the result when $\abs{a}< 1$ because, by the third assertion in Lemma \ref{rem:psi}, 
\begin{align}
\left|F(a,\mt{b}) \right| \leq & \frac{\gamma  }{2\pi }  \left( \frac{x_0}{2\pi}\right)^{p/2}  \left(N+\frac{1}{2}\right)^{k_q/2} \left(N(Rx_0U)+1/2\right)^{(p-1)/(2q)} \int_{\R}  \abs{t}^{p/2} \lambda(t) dt  \notag\\
\leq & \frac{  \gamma  U }{\pi(1+p/2) }  \left( \frac{Ux_0}{2\pi}\right)^{p/2} \left(\frac{5}{4}N\right)^{k_q/2} \left(\frac{5}{4}N(Rx_0U)\right)^{(p-1)/(2q)}.\label{ee} 
\end{align} 
Because $ t \mapsto \psi_{\widetilde{\mt{N}}(q)}^{Rx_0t}\left(\mt{b}/R\right)$ is analytic (see \cite{Fuchs} page 320), the function $$ t \mapsto \left(\frac{x_0\abs{t}}{2\pi} \right)^{p/2} \lambda(t) \psi_{\widetilde{\mt{N}}(q)}^{Rx_0t}\left( \frac{\mt{b}}{R}\right) \in
C^{\infty}(\R)$$ and its derivatives have compact support. 
By integration by parts, we obtain, for all $a\neq0$,
\begin{align} \abs{F(a,\mt{b})}\leq  \frac{\gamma}{\pi a^2 R^{p/2} }\int_{U/2}^U \abs{\frac{\partial^2}{\partial t^2} \left(\left( \frac{Rx_0t}{2\pi} \right)^{p/2} \lambda(t) \psi_{\widetilde{\mt{N}}(q)}^{Rx_0t}\left( \frac{\mt{b}}{R}\right)  \indic\{\abs{\mt{b}}_{\infty}  \leq R\}\right) }dt. \notag
\end{align}  
The result when $|a|\ge 1$ is obtained by  $1+(a/\tau)^2\le(1+1/\tau^2)a^2$, so by \eqref{eq:pos1},\\ 
$\gamma UC_{8}(Rx_0U, p, U)   N^{2+k_q/2} /(2a^2 )  \leq    1/\left(2^{p/2} \tau(1+(a/\tau)^2)\right)$, and by
Lemma \ref{lem:secondH}, for all $(a,\mt{b})\in\R^{p+1}$ such that $|a|\ge1$,
\begin{equation} 
\abs{F(a,\mt{b})}\leq \frac{\gamma UC_{8}(Rx_0U, p, U)}{2\pi a^2 R^{p/2} }  N^{2+k_q/2}\indic\{\abs{\mt{b}}_{\infty} \leq R\}.   \label{eq:H_N_IPP}
\end{equation}
$f_{1,n}=f_0$ has integral 1 and so has $f_{2,n}$ by Fubini's theorem and that $\psi_{N}^{c}$ is odd when $N$ is odd.\\
\textbf{Step 1.2.} We give conditions for $ f_{1,n} ,f_{2,n}\in\mathcal{H}^{q,\phi,\omega}_{w,W}(l) $.
By \eqref{eq:hyp}-\eqref{eq:H_N}, and because, by 
Step 1.1, for all $(a,\mt{b})\in\R^{p+1}$,  $f_{2,n}(a,\mt{b})^2 \leq 4 f_{1,n}(a,\mt{b})^2 $, $f_{1,n}$ and $f_{2,n}$ belong to $L^2\left(w \otimes W_{[-R,R]}^{\otimes p}\right)$. Let us show that $ f_{2,n}$, hence $f_{1,n}$ ($f_{2,n}$ with $\gamma =0$), satisfy the first condition in $\mathcal{H}^{q,\phi,\omega}_{w,W}(l) $  if
\begin{align}
&  2\left(  \int_{0}^{\infty}  \phi(t)^2e^{-2\tau t}dt + \gamma ^2\left(\frac{Rx_0U}{2\pi}\right)^p\frac{ \phi(U)^2 U }{p+1}  \right) \leq  \pi \l^2 \label{eq:cond31}\\
& \frac{C_{12}(\sigma,p)}{\tau p^{2\sigma/q}} +  \gamma ^2\frac{2U p^{2\sigma/q} N^{2\sigma}}{p+1}\left(\frac{Rx_0U}{2\pi}\right)^p \leq  \pi \l^2. 
\label{eq:cond32}
\end{align}
Let $\boldsymbol{m}\in\N_0^p$ and $c^P_{\mt{m}}(t):=  \braket{ 1/2^{p/2} ,\psi_{\mt{m}}^{Rx_0t}}_{L^2([-1,1]^p)}$. By 
Proposition \ref{pdebut} \eqref{pdebutiii}, change of variables, and for all $ t\in\R $, $ \mathcal{F}_{1\mathrm{st}}\left[f_0(\cdot ,\star) \right](t) = e^{-\abs{t}\tau} \indic\{\abs{\star}_{\infty} \leq R\} /(2R)^{p/2} $, we have 
\begin{align}
b_{\boldsymbol{m}}^2(t)&= \frac{1}{i^{\abs{\mt{m}}_1}}
\left(e^{-\tau \abs{t}}c^P_{\mt{m}}(t)	 
+ \gamma   \indic\{ \boldsymbol{m} =  \widetilde{\mt{N}}(q)\}\left(\frac{Rx_0\abs{t}}{2\pi} \right)^{p/2} \lambda(t)\right).\label{eq:b2m}
\end{align}	
Because $\left( \psi_{\boldsymbol{m}}^{Rx_0t}\right)_{\mt{m}\in\N_0^p}$ is an orthonormal basis, we have
\begin{align}\label{eq:FT}
\forall t\ne0,\ \sum_{\boldsymbol{m}\in\N_0^p}\abs{b_{\boldsymbol{m}}^2(t)}^2 & \leq 2  \left(e^{-2\tau \abs{t}} +\gamma ^2  \left( \frac{Rx_0\abs{t}}{2\pi} \right)^{p}  \lambda(t)^2\right) .
\end{align}
The first part of the first condition in $\mathcal{H}^{q,\phi,\omega}_{w,W}(l) $ holds by \eqref{eq:cond31} and because, by \eqref{eq:FT}, 
\begin{align*}
\sum_{\boldsymbol{m}\in\N_0^p}\int_{\R} \phi(t)^2 \left | b^2_{\boldsymbol{m}}(t) \right |^2dt 
& \leq   4\left( \int_{0}^{\infty}  \frac{\phi(t)^2}{e^{2\tau t}}dt+ \gamma ^2 \left(\frac{Rx_0}{2\pi}\right)^p  \int_{U/2}^{U}  \phi(t)^2  t^{p}\lambda^2(t) dt  \right).
\end{align*} 
The second part of the first condition holds by \eqref{eq:cond32} and because, by \eqref{eq:b2m} and Lemma \ref{lem:cst_f_0}, 
for all $ \tau \geq  \left(3e^{1/2} Rx_0/8\right)\vee (1/2)$ and $N\geq N(Rx_0U)$, 
\begin{align*}
\sum_{\boldsymbol{m}\in\N_0^p} |\boldsymbol{m}|_{q}^{2\sigma}  \int_{\R} \left | b_{\boldsymbol{m}}^2(t) \right |^2dt
&  \leq 2 \int_{\R}e^{-2 \tau\abs{t}} \sum_{\boldsymbol{m}\in\N_0^p}  |\boldsymbol{m}|_{q}^{2\sigma}\left(c^P_{\mt{m}}(t)\right)^2 dt \\
& \quad + 2\gamma ^2 \left(\frac{Rx_0}{2\pi}\right)^p\left|\widetilde{\mt{N}}(q)\right|_{q}^{2\sigma} \int_{\R}  \lambda^2(t) \abs{t}^{p} dt      \\
&  \leq 2\left(  \frac{C_{12}(\sigma,p)}{\tau p^{2\sigma/q}} +  \gamma ^2\frac{2U p^{2\sigma/q} N^{2\sigma}}{p+1}\left(\frac{Rx_0U}{2\pi}\right)^p  \right).
\end{align*}
\noindent\textbf{Step 2.} \eqref{ii} holds  with $4r(n)^2 = \gamma ^2 \left(Rx_0/(2\pi)\right)^p  \int_{U/2}^U t^p \lambda(t)^2 dt/\pi $.\\ 
\noindent\textbf{Step 3.} 
 (ii) page 97 in \cite{Tsybakov} yields $\chi_2(\mathbb{P}_{2,n}, \mathbb{P}_{1,n}) = \left( 1+ \chi_2\left(\mathbb{P}_2, \mathbb{P}_1\right)\right)^n - 1$ so
\begin{align*}
\chi_2(\mathbb{P}_{2,n}, \mathbb{P}_{1,n})
&= n\int_0^{\chi_2\left(\mathbb{P}_2, \mathbb{P}_1\right)}(1+u)^{n-1}du\le n\chi_2\left(\mathbb{P}_2, \mathbb{P}_1\right)\exp\left((n-1)\chi_2\left(\mathbb{P}_2, \mathbb{P}_1\right)\right).
\end{align*}
Thus, if $ \chi_2\left(\mathbb{P}_{2}, \mathbb{P}_1\right)  \leq 1/n $, we obtain $\chi_2(\mathbb{P}_{2,n}, \mathbb{P}_{1,n}) \leq e n  \chi_2\left(\mathbb{P}_2, \mathbb{P}_1\right)$. We have
\begin{align*}
\chi_2\left(\mathbb{P}_2, \mathbb{P}_1\right)
= \int_{\mathbb{S}_{\mt{X}}} \int_{\R}  \frac{f_{\mt{X}}(\mt{x})\left(f^1_{Y|\mt{X}}(y|\mt{x}) - f^2_{Y|\mt{X}}(y|\mt{x})\right)^2}{f^1_{Y|\mt{X}}(y|\mt{x})} d\mt{x}dy   
\end{align*}
and, 
for all $(y,\mt{x}) \in \R\times\mathbb{S}_{\mt{X}}$ such that $\mt{x}^{\underline{\mt{1}}}\ne 0$, 
\begin{align*}
f^1_{Y|\mt{X}}(y|\mt{x})
&= \frac{1}{\pi \tau (2R)^{p/2} \abs{\mt{x}^{\underline{\mt{1}}}}}  \int_{\R^p} \frac{ \prod_{k=1}^p \indic\{\abs{\mt{u}_k}\leq  \abs{\mt{x}_k}R\} }{\left(\left(y -  \sum_{k=1}^p \mt{u}_k\right)/\tau\right)^2 + 1}  d\mt{u}  \\
&\geq \frac{(2R)^{p/2} }{\pi \tau}  \inf_{\abs{u} \leq \abs{\mt{x}}_1R} \frac{1}{((y - u)/\tau)^2 + 1}.
\end{align*} 
This yields, using  $\mathbb{S}_{\mt{X}}= [-x_0,x_0]^p$ and Parseval's identity, 
\begin{align*}
&\chi_2\left(\mathbb{P}_2, \mathbb{P}_1\right)\\
&\leq \frac{\pi \tau  C_{\mt{X}}}{(2R)^{p/2} }\int_{[-x_0,x_0]^p} \int_{\R}  \left( \frac{2y^2}{\tau^2} + \frac{2(\abs{\mt{x}}_1R)^2}{\tau^2} + 1\right) \left(f^1_{Y|\mt{X}}(y|\mt{x}) - f^2_{Y|\mt{X}}(y|\mt{x})\right)^2 d\mt{x}dy   \\
&= \frac{ C_{\mt{X}}x_0^p\gamma ^2 }{\tau(2R)^{p/2} }   \int_{[-1,1]^p} \int_{\R}  \left|\partial_t \mathcal{F}\left[F\right](t, x_0t\mt{x}) \right|^2   +  \left((x_0pR)^2 + \frac{\tau^2}{2}\right)  \left|\mathcal{F}\left[F\right](t, x_0t\mt{x}) \right|^2  d\mt{x} dt.
\end{align*}
Lemmas \ref{lem:I12} and \ref{upper_bound} yield $\chi_2\left(\mathbb{P}_2, \mathbb{P}_1\right) \leq   C_{18}(U,x_0,R,\tau)  \gamma ^2N^2 \left(e Rx_0U/(4N)\right)^{2k_qN}$,
\begin{align*}
C_{18}(U,x_0,R,\tau):=& \frac{C_{\mt{X}}}{\tau} 
\left(\frac{R^{3/2}x_0Ue^3}{9\sqrt{2}}\right)^pC_{17}(Rx_0U,p,U) \left(\frac{eRx_0U}{4N(Rx_0U)}\right)^{\frac{2(p-1)N(Rx_0U)}{q}}\\
&+\frac{C_{\mt{X}}}{\tau} 
\left(\frac{R^{3/2}x_0Ue^3}{9\sqrt{2}}\right)^pU\frac{2(x_0pR)^2 +\tau^2}{2N(Rx_0U)^{2}} \left(\frac{eRx_0U}{4N(Rx_0U)}\right)^{\frac{2(p-1)N(Rx_0U)}{q}}.
\end{align*}
As a result, \eqref{iii} is satisfied if
\begin{equation} \label{eq:lower44}
n  \gamma ^2 N^2   \exp\left(-2k_qN\ln\left(\frac{4N}{e Rx_0U}\right)\right) \le\frac{\xi^2}{e C_{18}(U,x_0,R,\tau)}.
\end{equation}
\noindent \textbf{Step 4.} 
We take $U= 4/(eRx_0)$, $N = 2\left\lceil\underline{N}\right\rceil+1$ for $\underline{N}$ going to infinity with $n$,  and $\tau$ such that
$  \int_{0}^{\infty}  \phi(t)^2e^{-2\tau t}dt \bigvee C_{12}(\sigma,p)/(2\tau p^{2\sigma/q})\le
\pi l^2/4$.
Thus $N(Rx_0U)$ is universal and $ \underline{N} \geq N(Rx_0U)$ and $N\le(9/2)\underline{N}$ for $n$ large enough. \eqref{eq:pos}, \eqref{eq:cond31}-\eqref{eq:cond32} (by the pigeonhole principle), and \eqref{eq:lower44} hold for $n$ large enough if 
\begin{align}
& \gamma \underline{N}^{2+k_q/2}\leq    \frac{1}{\left(\tau+1/\tau\right)(9/2)^{1+(p+k_q)/2}C_{8}(4/e, p, U)U},\label{eq:posb}\\
&\gamma \underline{N}^{\sigma}
\leq  \frac{l}{2p^{\sigma/q}(9/2)2^{\sigma}}\sqrt{\frac{\pi(p+1)}{U}} \left(\frac{e\pi}{2}\right)^{p/2},
\label{eq:cond3b}\\
& \gamma  \leq\frac{l}{\phi(U)} \sqrt{\frac{\pi(p+1)}{ U }} \left(\frac{e\pi}{2}\right)^{p/2},\label{eq:cond3b1}\\
&n\gamma^2\underline{N}^2 \exp\left(-4k_q\underline{N}\ln\left(\underline{N}\right)\right)\le\frac{\xi^2}{(9/2)^2 e C_{18}(U,x_0,R,\tau)}\label{eq:lower44b},
\end{align}
and $\gamma$ goes to 0 with $n$. 
Taking 
$ \gamma  = C_\gamma\underline{N}^{-(2+k_q/2)\vee\sigma}/\left(\left(C_8(4/e,p,U)U\right)\wedge \sqrt{U}\right)$ for a small enough $C_\gamma$ depending on $l$, $\phi$, $\sigma$, $p$, and $q$, \eqref{eq:posb}-\eqref{eq:cond3b} hold because $Rx_0$, hence $U$ is fixed. Then, with  $\underline{N}=3\ln(n)/\left(8k_q\ln_2(n)\right)$,  \eqref{eq:lower44b} becomes, for $n$ large enough, 
\begin{align*}&\frac{C_{\gamma}^2 }{\sqrt{n}} \exp\left(  \frac{3\ln(n)\ln(8k_q\ln_2(n)/3)}{4\ln_2(n)} -\left(\left(2+\frac{k_q}{2}\right)\vee \sigma -1\right)\ln\left( \frac{3\ln(n)}{8k_q\ln_2(n)}\right)\right)\\
&\le\frac{\xi^2\left(C_8(4/e,p,U)U\right)^2\wedge U}{8 e C_{18}(U,x_0,R,\tau)}.
\end{align*}
Moreover, we have $ r(n)^2 = \underline{N}^{-2((2+k_q/2)\vee\sigma)}  C_{\gamma}^2  \left(Rx_0/(2\pi)\right)^p \int_{U/2}^U t^p \lambda(t)^2 dt  /(4\pi).\square$\\
All other steps 2 are the same as for (T\ref{theo:lower}.1\ref{it:lowerindic}).\\ 

\noindent {\bf Proof of (T\ref{theo:lower}.2\ref{it:lowerch}).} Denote by $E:=L^2(\R)\times L^2(\R)$. Equip $E$ with $\braket{\mt{g},\mt{h}}_{E}^2 = \braket{\mt{g}_{1},\mt{h}_{1}}_{L^2(\R)}^2 +\braket{\mt{g}_{2},\mt{h}_{2}}_{L^2(\R)}^2$.
Denote by 
$\mathbb{P}_{j,n}^{ \mt{m}}$ the law of $\left(\left(\mathfrak{R}\left(Z^j_{\mt{m}}(t)\right)\right)_{ \ t\in\R} , \left(\mathfrak{I}\left(Z^j_{\mt{m}}(t)\right)\right)_{ \ t\in\R} \right) $ in $E$ and by
$ \mathbb{P}_{j,n}$ the law on $\ell_2\left(E\right) $ 
of the sequence indexed by $\mt{m}\in \N_0^p$. The latter can be defined as a function of 
$f_{j,n}$ or  $ \left(b^j_{\mt{m}}(t)\right)_{\mt{m}\in\N_0^p, \ t\in\R} $, for $j=1,2$. 
Take $f_{1,n}:=0$ and  $ f_{2,n} $ like \eqref{eq:hyp} with
$\widetilde{\mt{N}}(1):= (N,\underline{\mt{0}}^{\top})^{\top}\in\N_0^p$.   
\eqref{eq:b2m} yields, for all $ \mt{m} \in\N_0^p$, $ b_{\mt{m}}^2(t) =i^{-\abs{\mt{m}}_1} \gamma  \indic\{ \boldsymbol{m} =  \widetilde{\mt{N}}(q)\} \left( Rx_0\abs{t}/2\pi \right)^{p/2}  \lambda(t) $. 
By independence, we have, for $j=1,2$, $\mathbb{P}_{j,n} = \bigotimes_{ \mt{m}\in \N_0^p} \mathbb{P}_{j,n}^{ \mt{m}}$.\\
\noindent \textbf{Step 1.} 
Using \eqref{ee} and \eqref{eq:H_N_IPP}, we have
$f_{2,n}\in L^2\left(w \otimes \cosh\left(\cdot/R\right)^{\otimes p}\right) $ and 
$f_{2,n}\in \mathcal{H}^{q,\phi,\omega}_{w,W}(l)$ if 
\begin{align}\label{eq:lower2}
& \left(\frac{Rx_0U}{2\pi}\right)^p \frac{2U\gamma ^2 }{p+1}\left(\phi(U)  \vee N^{\sigma}\right)^2\leq  \pi \l^2.
\end{align}
\noindent \textbf{Step 3.} Let $\xi<\sqrt{2}$, $ G_{\widetilde{\mt{N}}(q)}^W=\left( \mathfrak{R}\left(\sigma_{\widetilde{\mt{N}}(q)}^{W,x_0\cdot} b_{\widetilde{\mt{N}}(q)}^2(\cdot) \right),   \mathfrak{I}\left(\sigma_{\widetilde{\mt{N}}(q)}^{W,x_0\cdot} b_{\widetilde{\mt{N}}(q)}^2(\cdot) \right)  \right)^{\top} $, 
$\mathcal{Q}$ the covariance operator of $\mathbb{P}_{1,n}^{\widetilde{\mt{N}}(q)}$ on $E$, and, for all $\mt{h}\in E$, $\mathcal{L}[\mt{h}]:=(\sigma/\sqrt{n})\left(\int_{0}^{\cdot} \mt{h}_1(s)ds,
\int_{0}^{\cdot}   \mt{h}_2(s)ds\right)^{\top}$. 
The reproducing kernel Hilbert space  $H_{\mathbb{P}_{1,n}^{\widetilde{\mt{N}}(q)}}$ of $\mathbb{P}_{1,n}^{\widetilde{\mt{N}}(q)}$ on $E$ is the image of $\mathcal{Q}^{1/2}$ with the scalar product of the image structure.  By Corollary B.3 in \cite{da2014stochastic} and $\mathcal{Q}=\mathcal{L}\mathcal{L}^*$, it is 
the image of $\mathcal{L}$ with the norm $\|\mt{f}\|_{\mathbb{P}_{1,n}^{\widetilde{\mt{N}}(q)}}^2= (n/\sigma^2)\left(\left\|\mt{h}_1\right\|_2^2+\left\|\mt{h}_2\right\|_2^2\right)$ for
$\mt{f}=\mathcal{L}[\mt{h}]$ and derived scalar product. By (2.12) 
in \cite{da2014stochastic}, the scalar product is also defined when one function belongs to $H_{\mathbb{P}_{1,n}^{\widetilde{\mt{N}}(q)}}$ for $\mathbb{P}_{1,n}^{\widetilde{\mt{N}}(q)}\ a.e.$ other function in $E$.  
By the Cameron-Martin formula (Proposition 2.26 in \cite{da2014stochastic}), 
\begin{align*} &\frac{\mathrm{d}\mathbb{P}_{2,n}^{\widetilde{\mt{N}}(q)}}{\mathrm{d}\mathbb{P}_{1,n}^{\widetilde{\mt{N}}(q)}}(y)
= \exp\left(     \braket{ y, \frac{\sqrt{n}}{\sigma} \mathcal{L}\left[G_{\widetilde{\mt{N}}(q)}^W\right]   }_{\mathbb{P}_{1,n}^{\widetilde{\mt{N}}(q)}}  - \frac{1}{2} \left| \frac{\sqrt{n}}{\sigma} \mathcal{L}\left[G_{\widetilde{\mt{N}}(q)}^W\right]  \right|^2_{\mathbb{P}_{1,n}^{\widetilde{\mt{N}}(q)}}  \right)\quad \mathbb{P}_{1,n}^{\widetilde{\mt{N}}(q)}\ a.s.,
\end{align*}
and, because $K(\mathbb{P}_{2,n}, \mathbb{P}_{1,n}) 
=\int_{E} \ln \left( \mathrm{d}\mathbb{P}_{2,n}^{\widetilde{\mt{N}}(q)}/\mathrm{d}\mathbb{P}_{1,n}^{\widetilde{\mt{N}}(q)}(y)\right)d\mathbb{P}_{2,n}^{\widetilde{\mt{N}}(q)}(y)$, we have
\begin{align*} 
K(\mathbb{P}_{2,n}, \mathbb{P}_{1,n}) &  =   \E\left[ \braket{ Z^2_{\widetilde{\mt{N}}(q)},  \frac{\sqrt{n}}{\sigma} \mathcal{L}\left[G_{\widetilde{\mt{N}}(q)}^W\right]   }_{\mathbb{P}_{1,n}^{\widetilde{\mt{N}}(q)}}   \right] - \frac{n}{2\sigma^2} \int_{\R}  \left|\sigma_{\widetilde{\mt{N}}(q)}^{W,x_0s} b_{\widetilde{\mt{N}}(q)}^2(s) \right|^2       ds .
\end{align*}
Because 
\begin{align*} 
&  \braket{ Z^2_{\widetilde{\mt{N}}(q)},  \frac{\sqrt{n}}{\sigma} \mathcal{L}\left[G_{\widetilde{\mt{N}}(q)}^W\right]   }_{\mathbb{P}_{1,n}^{\widetilde{\mt{N}}(q)}}  = \left|  \frac{\sqrt{n}}{\sigma} \mathcal{L}\left[G_{\widetilde{\mt{N}}(q)}^W\right]   \right|^2_{\mathbb{P}_{1,n}^{\widetilde{\mt{N}}(q)}} + \braket{ \left(\begin{array}{c}
	B_{\widetilde{\mt{N}}(q)}^{\mathfrak{R}}
	\\ 
	B_{\widetilde{\mt{N}}(q)}^{\mathfrak{I}}
	\end{array}  \right)       ,  \mathcal{L}\left[G_{\widetilde{\mt{N}}(q)}^W\right]   }_{\mathbb{P}_{1,n}^{\widetilde{\mt{N}}(q)}},
\end{align*}
and the second term in the right-hand side is a limit in quadratic mean of  mean zero Gaussian random variables, hence has mean zero (see the arguments page 41 in \cite{da2014stochastic}), we have 
\begin{align}\label{eq:K_low} 
K(\mathbb{P}_{2,n}, \mathbb{P}_{1,n}) =  \frac{n}{2\sigma^2} \int_{\R}  \left|\sigma_{\widetilde{\mt{N}}(q)}^{W,x_0t} b_{\widetilde{\mt{N}}(q)}^2(t) \right|^2       dt .
\end{align}
By Proposition \ref{pdebut} \eqref{pdebutii}, we have 
$K(\mathbb{P}_{2,n}, \mathbb{P}_{1,n}) =  \gamma ^2nR^p   \int_{\R} \left( \sigma_{\widetilde{\mt{N}}(q)}^{\cosh,Rx_0t}  \right)^2  \left(Rx_0\abs{t}/(2\pi)\right)^{p} \lambda(t)^2 dt/(2\sigma^2)
=\gamma ^2n R^p \int_{\R} \rho_{\widetilde{\mt{N}}(q)}^{\cosh,Rx_0t} \lambda(t)^2 dt/(2\sigma^2)$ and, 
by Theorem 7 in \cite{Note} (there is difference of normalisation for $\mathcal{Q}_{t}$ by a factor $1/(2\pi)$), for all $ U/2 \leq \abs{t} \leq	 U $ and $ Rx_0U <1$, 
$  \rho_{\widetilde{\mt{N}}(q)}^{\cosh,Rx_0t}  \leq 	\left(Rx_0 U	e/\left(\pi	(1-(Rx_0 U)^2)\right)\right)^p\exp\left( 2k_q N \ln\left(Rx_0U\right) \right).$
Thus Lemma \ref{lem:lower} 
\eqref{iii} holds if
\begin{equation} \label{eq:lower4_0}
n \gamma ^2\exp\left(-2k_qN \ln\left(\frac{1}{Rx_0U}\right) \right)\le \xi^2\left(\frac{\pi	(1-(Rx_0 U)^2)}{R^2x_0 U	e}\right)^p \frac{2\sigma^2}{U}.
\end{equation}
\noindent \textbf{Step 4.} 
Let $U= 1/(eRx_0)$, $N =\left \lceil \underline{N} \right\rceil$,  
$ \gamma = C_{\gamma} \xi\sigma\sqrt{2eRx_0}\left(\pi	\left(1-e^{-2}\right)/R\right)^{p/2}/\underline{N}^{\overline{\sigma}\vee\sigma}$ for $\overline{\sigma}>1/2$ and $C_{\gamma}=l \left(2eR/\left(1-e^{-2}\right)\right)^{p/2}\sqrt{(p+1)e\pi Rx_0}/\left(2^{\sigma+1}  \xi\sigma\sqrt{2eRx_0}\right)$, and $\underline{N}= \ln(n/\ln(n))/(2k_q)$. 
\eqref{eq:lower4_0} holds if 
$n C_{\gamma}^2   \exp\left(-2k_q\underline{N}-2\left(\overline{\sigma}\vee \sigma\right)\ln\left(\underline{N}\right)\right)=   C_{\gamma}^2 (2k_q)^{2\left(\overline{\sigma}\vee \sigma\right)} \ln(n)^{1- 2(\overline{\sigma}\vee \sigma)}
\ln_2(n)^{2(\overline{\sigma}\vee \sigma)}\le 1$, so \eqref{eq:lower2} and \eqref{eq:lower4_0}  
hold for $n$ large enough. 
\hfill $\square$ \vspace{0.3cm} 

\noindent \textbf{Proof of (T\ref{theo:lower}.1\ref{it:expindic}).} 
\noindent \textbf{Step 1.} 
By the proof of (T.1\ref{it:lowerch}),  
$ f_{2,n}\in L^2\left(w \otimes W_{[-R,R]}^{\otimes p}\right)$ and $f_{2,n}\in \mathcal{H}^{q,\phi,\omega}_{w,W}(l)$ if 
\begin{equation}\label{eq:lower2_1}\left(\frac{Rx_0U}{2\pi}\right)^p \frac{2U\gamma^2}{p+1}\left(\phi(U)\bigvee \exp\left(\kappa N\ln\left(N+1\right)\right)\right)^2\le \pi l^2.\end{equation}
\noindent \textbf{Step 3.} Let $\xi<\sqrt{2}$ and $8/(eRx_0U) \geq 1$. 
By Lemma \ref{upper_bound}, we have, 
for all $ U/2 \leq \abs{t} \leq U $, 
$\left(\sigma_{\widetilde{\mt{N}}(q)}^{W_{[-1,1]},Rx_0t}\right)^2 \leq \left(2\pi e^3/9\right)^p \exp\left(-2k_qN\ln\left(4 (N+3/2)/(eRx_0U)\right) \right)$ and, by \eqref{eq:K_low} and Proposition \ref{pdebut} \eqref{pdebutii}, 
Lemma \ref{lem:lower} \eqref{iii} holds if
\begin{equation} \label{eq:lower4}
n  \gamma ^2\exp\left(-2k_qN\ln\left( \frac{4 (N+3/2)}{eRx_0U}\right)\right)\le \frac{(p+1)\xi^2\sigma^2}{U}  \left(\frac{9}{R^2Ux_0 e^3}\right)^p.
\end{equation}
\noindent \textbf{Step 4.}  Let $U= 4/(eRx_0)$,  $\gamma = \widetilde{C}_{\gamma}\exp\left(-\kappa N\ln\left(N+1\right)\right)$, $\widetilde{C}_{\gamma}=l (\pi e/2)^{p/2}\sqrt{(p+1)e\pi Rx_0/8}$,   
$ N =\left \lceil \underline{N} \right\rceil$,  
$2(\kappa+k_q) \underline{N} \ln\left(\underline{N}+1\right)=\ln\left(C_{\gamma}^2n\right)$, $C_{\gamma}=l(2\pi Re^3/9)^{p/2}\sqrt{\pi/2}/(\xi \sigma)$.  Under such a choice, 
\eqref{eq:lower2_1} and \eqref{eq:lower4} hold for $n$ large enough. Moreover, we have  
$ r(n)=C_r\exp\left(-\kappa N\ln\left(N+1\right)\right)$, where $C_r=\widetilde{C}_{\gamma}\left(Rx_0/(2\pi)\right)^{p/2}  \sqrt{\int_{U/2}^U \abs{t}^p \lambda(t)^2 dt /(4 \pi)}$, $\underline{N}\ln\left(\underline{N}\right)\le N\ln(N+1)\le \left(\underline{N}+1\right)\ln\left(\underline{N}+2\right)$
\begin{align*}
N\ln(N+1)&\le \underline{N}\ln\left(\underline{N}+1\right)+\ln\left(\underline{N}+1\right)+1+o(1)
\\&
= \ln\left(\left(C_{\gamma}^2n\right)^{1/(2\kappa+2 k_q)}\right)+(1+o(1))\ln_2\left(\left(C_{\gamma}^2n\right)^{1/(2\kappa+2 k_q)}\right),
\end{align*}
indeed, using iteratively the definition of $\underline{N}$, 
$\ln\left(\underline{N}+1\right)=\ln\left(\underline{N}\right)+(1+o(1))/\underline{N}$ so $\ln\left(\underline{N}+1\right)=\ln\left(\underline{N}\right)\left(1+(1+o(1))/\ln\left(\left(C_{\gamma}^2n\right)^{1/(2\kappa+2 k_q)}\right)\right)$ and 
\begin{align*}\ln\left(\underline{N}\right)&=\ln_2\left(\left(C_{\gamma}^2n\right)^{1/(2\kappa+2 k_q)}\right)-\ln_2\left(\underline{N}+1\right)\\
&=\ln_2\left(\left(C_{\gamma}^2n\right)^{1/(2\kappa+2 k_q)}\right)-\ln_2\left(\underline{N}\right)+(1+o(1))/\ln\left(\left(C_{\gamma}^2n\right)^{1/(2\kappa+2 k_q)}\right)\\
&=\ln_2\left(\left(C_{\gamma}^2n\right)^{1/(2\kappa+2 k_q)}\right)-(1+o(1))\ln_3\left(\left(C_{\gamma}^2n\right)^{1/(2\kappa+2 k_q)}\right)
\end{align*}
so
$1/\ln\left(\left(C_{\gamma}^2n\right)^{(1+o(1))\kappa/(2\kappa+2 k_q)}\right)\le r(n)\left(C_{\gamma}^2n\right)^{\kappa/(2\kappa+2 k_q)}/C_r\le 1$. \\
\noindent {\bf Proof of (T.\ref{theo:lower}.2\ref{it:expch}).} 
Let $U= 2/(eRx_0)$,  $\gamma = \widetilde{C}_{\gamma}\exp\left(-\kappa N\right)$, $\widetilde{C}_{\gamma}=l (\pi e)^{p/2}\sqrt{(p+1)e\pi Rx_0/2}$, $ N =\left \lceil \underline{N} \right\rceil$,  
$\underline{N} = \ln(n)/(2\kappa+2k_q)$, 
$C_{\gamma}=l\left(4\pi e/\left(\pi	\left(1-e^{-2}\right)\right)\right)^{p/2}\sqrt{\pi/2}/(\xi \sigma)$.  Under such a choice, 
$4\left(Rx_0U/(2\pi)\right)^p (2U\gamma^2/(p+1))\left(\phi(U)\bigvee \exp\left(\kappa N\right)\right)^2\le \pi l^2$
and \eqref{eq:lower4_0} hold for $n$ large enough, hence steps 1 and 3. By Step 2, we have  
$r(n)=C_r\exp\left(-\kappa N\right)\ge C_r\exp\left(-\kappa \underline{N}\right)/e$. 
$\square $

\subsection{Upper bounds}
We use, for all $\epsilon>0$, $N\in\N_0^{\R}$, $N_0\in\N_0$, $T>0$, $\widetilde{F}^{q,N,T,0}_1$ and $ \widetilde{f}_{\alpha,\mt{\beta}}^{q,N, T,\epsilon}$ which are defined like  $ \widehat{F}^{q,N,T,0}_1$ and $ \widehat{f}_{\alpha,\mt{\beta}}^{q,N,T,\epsilon} $ replacing $ \widehat{c}_{\mt{m}}$ by $\widetilde{c}_{\mt{m}}$ and $\widehat{F}_1^{q,N,T,0} $ by  $ \widetilde{F}^{q,N,T,0}_1$, $\widetilde{c}_{\mt{m}}:= \sum_{j=1}^n e^{i\star Y_j}  \overline{g_{\mt{m}}^{W,x_0\star}}\left(\mt{X}_j/x_0\right)\indic\{\mt{X}_j\in\mathcal{X}\}/(nx_0^p  f_{\mt{X}| \mathcal{X}}(\mt{X}_j))$, $Z_{n_0}:=\sup_{f_{\mt{X}|\mathcal{X}} \in \mathcal{E}} \left\| \Delta_f f_{\mt{X}|\mathcal{X}}\right\|_{L^{\infty}(\mathcal{X})}^2$, $\Delta_f:=1/\widehat{f}_{\mt{X}|\mathcal{X}}^{\delta}- 1/f_{\mt{X}|\mathcal{X}}$, $L := (2\pi)^p\left\| \mathcal{F}_{1\mathrm{st}}\left[f_{\alpha,\mt{\beta}}\right](\star,\cdot_2) \right\|^2_{L^2(\R^p)}$, $\widetilde{\omega}^{q,W,c}_{N(\star)}: = \sup_{\abs{\mt{m}}_{q}\leq N(\star)}1/\rho_{\mt{m}}^{W,c}$,  $\Delta_{\mt{m}}:=  \sum_{j=1}^nZ_j^{\mt{m},\star}/n$, $Z_j^{\mt{m},\star}: = e^{i\star Y_j}\Delta_f(\mt{X}_j)
\overline{g_{\mt{m}}^{W,x_0\star}\left(\mt{X}_j/x_0\right)} \indic\{\mt{X}_j\in\mathcal{X}\}/x_0^p$,\\
$S_0^N(\star,\cdot_2):=\sum_{\abs{\mt{m}}_{q}\leq N(\star)} g_{\mt{m}}^{W,x_0\star}(\cdot_2) \Delta_{\mt{m}} (\star)$,
$S_1^N(\star,\cdot_2):=\sum_{\abs{\mt{m}}_{q}\leq N(\star)} g_{\mt{m}}^{W,x_0\star}(\cdot_2) \E\left[\Delta_{\mt{m}} (\star)   
\right]$,\\
 $S_2^N(\star,\cdot_2):= \sum_{\abs{\mt{m}}_{q}\leq N(\star)} g_{\mt{m}}^{W,x_0\star}(\cdot_2) \left(\Delta_{\mt{m}} (\star) - \E\left[\Delta_{\mt{m}} (\star) \right] \right)$, \\
\begin{align*}
K_1:=&  \left\| \indic\left\{ \epsilon \leq \abs{\star} \right\}  \left(\widehat{F}_1^{q,N,T,0} -  \mathcal{F}_{1\mathrm{st}}\left[f_{\alpha,\mt{\beta}}\right]\right)(\star, \cdot_2)\right\|^2_{L^2\left(1\otimes W^{\otimes p}\right)},\\
K_2 := &\left\| \indic\left\{ \abs{\star} < \epsilon \right\}  \left(  \mathcal{I}_{\underline{a},\epsilon}\left[ \widehat{F}_1^{q,N,T,0} \right] -\mathcal{F}_{1\mathrm{st}}\left[f_{\alpha,\mt{\beta}}\right]\right) (\star, \cdot_2)  \right\|^2_{L^2\left(1\otimes W^{\otimes p}\right)}\\
R_1(\star,\cdot_2):=&\indic\{\epsilon \le |\star| 	\} \left(\widetilde{F}_1^{q,N, T, 0}- F_1^{q,N,T,0}\right)(\star,\cdot_2),\\
R_2(\star,\cdot_2):=&\indic\{\epsilon \le |\star| 	\}\left(\widehat{F}_1^{q,N, T, 0}- \widetilde{F}_1^{q,N,T,0}\right)(\star,\cdot_2),\\
R_3(\star,\cdot_2):=&\indic\{\epsilon \le |\star| 	\} \left(F_1^{q,N, T, 0}- F_1^{q,\infty,T,0}\right)(\star,\cdot_2),\\
R_4(\star,\cdot_2):=&\indic\{\epsilon \le |\star| 	\}   \left( F_1^{q,\infty,T,0}-  \mathcal{F}_{1\mathrm{st}}\left[f_{\alpha,\mt{\beta}}\right]\right)(\star,\cdot_2), \\
 \mathcal{R}^W_{n_0,\sup}:=&\sup_{f_{\alpha,\mt{\beta}} \in \mathcal{H}^{q,\phi,\omega}_{w,W}(l, M)\cap\mathcal{D},\ f_{\mt{X}|\mathcal{X}} \in \mathcal{E}} \mathcal{R}_{n_0}^{W}\left(\widehat{f}_{\alpha,\mt{\beta}}^{q,N,T,\epsilon},f_{\alpha,\mt{\beta}}\right),\\
\Delta^{W}_0(\star,N_0,n,z) :=&\frac{2}{\pi(2\pi)^p} \frac{c_{\mt{X}}|\star|^p}{n} \nu_q^{W}(x_0\star,N_0) \\
&\quad     + \frac{2z}{\pi(2\pi)^p}\left(L(\star) + \frac{ c_{\mt{X}} (N_0+1)^p\abs{\star}^p}{n} \right) \widetilde{\omega}^{q,W,x_0\star}_{N_0}.
\end{align*}
\begin{lemma}\label{lem:def_est}
	For all $ \mt{m} \in \N_0^p $, we have 
	$\E\left[ \widetilde{c}_{\mt{m}}(t)\right] =  c_{\mt{m}}(t)$ and $\E\left[\abs{ \widetilde{c}_{\mt{m}}(t) - c_{\mt{m}}(t)}^2\right] \leq c_{\mt{X}}/(nx_0^p)$.
\end{lemma}
\noindent {\bf Proof.} This comes from 
\begin{align*} \E\left[ \widetilde{c}_{\mt{m}}(t)\right] & = \frac{1}{x_0^p}\E\left[\frac{ e^{it Y}}{f_{\mt{X}|\mathcal{X}}(\mt{X})}  \overline{g_{\mt{m}}^{W,x_0t}}\left(\frac{\mt{X}}{x_0}\right)\indic\{\mt{X}\in\mathcal{X}\}\right] \\
&=  \frac{1}{x_0^p} \int_{\mathcal{X}}    \E\left[ e^{it\alpha+ it\mt{\beta}^{\top}\mt{x}}\right]  \overline{g_{\mt{m}}^{W,x_0t}}\left(\frac{\mt{x}}{x_0}\right)d\mt{x},\\
\E\left[\abs{ \widetilde{c}_{\mt{m}}(t) - c_{\mt{m}}(t)}^2\right] & \leq   \frac{1}{n x_0^{2p} }\E\left[ \left.\left| \frac{ e^{it Y}}{f_{\mt{X}|\mathcal{X}}(\mt{X})}  \overline{g_{\mt{m}}^{W,x_0t}}\left(\frac{\mt{X}}{x_0}\right) \right|^2\right|\mt{X}\in\mathcal{X}\right]  
\\& \leq    \frac{1}{nx_0^{2p}} \int_{\mathcal{X}}  \frac{1 }{f_{\mt{X}|\mathcal{X}}(\mt{x})} \left|  \overline{g_{\mt{m}}^{W,x_0t}}\left(\frac{\mt{x}}{x_0}\right) \right|^2 d\mt{x}\\ 
&\leq   \frac{c_{\mt{X}}}{n x_0^{p}} \  \int_{[-1,1]^p}  \left| \overline{g_{\mt{m}}^{W,x_0t}}\left(\mt{u}\right) \right|^2 d\mt{u}.\quad\square
\end{align*} 

\begin{lemma}\label{lem:Z_n0}
	If $ \widehat{f}_{\mt{X}|\mathcal{X}}$  satisfies 
	(H1.\ref{E4}) then $Z_{n_0}= O_{p}\left(v(n_0,\mathcal{E})/{\delta(n_0)}\right)$.
\end{lemma}
\noindent {\bf Proof.} 
For all $n_0$ large enough so that $ \sqrt{\delta(n_0)} c_{\mt{X}} \leq 1$ and $x\in \mathcal{X}$, we have 
\begin{align*}
\left| \left(\widehat{f}_{\mt{X}|\mathcal{X}}^{\delta}-f_{\mt{X}|\mathcal{X}}\right)(x)\right|\leq &\left| \left(\widehat{f}_{\mt{X}|\mathcal{X}}-f_{\mt{X}|\mathcal{X}}\right)(x)\right| \indic\left\{ \widehat{f}_{\mt{X}|\mathcal{X}}(x) \geq \sqrt{\delta(n_0)} \right\} \\
&+  \left| \sqrt{\delta(n_0)}-f_{\mt{X}|\mathcal{X}}(x)\right| \indic\left\{ \widehat{f}_{\mt{X}|\mathcal{X}}(x) <  \sqrt{\delta(n_0)}  \right\}\\
\leq& \left| \left(\widehat{f}_{\mt{X}|\mathcal{X}}-f_{\mt{X}|\mathcal{X}}\right)(x)\right|\ (\text{using} \ \sqrt{\delta(n_0)} c_{\mt{X}} \leq 1)
\end{align*}
and $\delta(n_0)Z_{n_0}\le\sup_{f_{\mt{X}|\mathcal{X}} \in \mathcal{E}}\left\|\widehat{f}_{\mt{X}|\mathcal{X}}-f_{\mt{X}|\mathcal{X}}\right\|_{L^{\infty}(\mathcal{X})}^2$. We conclude by (H1.\ref{E4}). $ \square $\vspace{0.3cm}

In the remaining, $\mathcal{E}$ is a class of densities, $f_{\mt{X}|\mathcal{X}} \in \mathcal{E}$, $\eta,l, M>0$, and $f_{\alpha,\mt{\beta}} \in \mathcal{H}^{q,\phi,\omega}_{w,W}(l, M)\cap\mathcal{D}$. By Lemma  \ref{lem:Z_n0}, there exists $M_{\mathcal{E},\eta}$ such that, for all $ n_0\in\N $, $\mathbb{P}\left(E\left(\mathcal{G}_{n_0},\mathcal{E},\eta\right)\right)\geq 1-\eta$, where $E\left(\mathcal{G}_{n_0},\mathcal{E},\eta\right):=\left\{ Z_{n_0} \leq M_{\mathcal{E},\eta}v(n_0,\mathcal{E})/\delta(n_0) \right\}$. We work on $E\left(\mathcal{G}_{n_0},\mathcal{E},\eta\right)$. 
\vspace{0.3cm}

\noindent {\bf Proof of theorems \ref{theo:compact} and \ref{theo:non_compact}.}
The proof consists in three parts.\\
In Part 1 
we show, for $W= W_{[-R,R]}$ and $W= \text{cosh}(\cdot/R)$, 
\begin{equation}
\mathcal{R}^W_{n_0,\sup}\le C\Bigg(
\int_{\epsilon\le \abs{t}\le T} \Delta_0^{W}(t,N,n,Z_{n_0})  dt +4l^2 \left( \sup_{t\in \R} \frac{1}{\omega_{N(t)+1}^2}+ \frac{1}{\phi(T)^2} \right)\Bigg) + CM^2\widetilde{w}(\underline{a}). \label{eq:START}
\end{equation}
In Part 2 we take $W= W_{[-R,R]}$ and, particularising \eqref{eq:START} to the different smoothness cases, obtain (T2.\ref{t_comp_1}),  (T\ref{theo:compact}.\ref{poly00}), (T\ref{theo:compact}.\ref{poly01}), and (T\ref{theo:compact}.\ref{t_comp_3}) in Theorem \ref{theo:compact}. 
In Part 3 we proceed similarly for the weight $W= \text{cosh}(\cdot/R)$ and prove (T\ref{theo:non_compact}.\ref{t_noncomp_1}) and (T\ref{theo:non_compact}.2) in Theorem \ref{theo:non_compact}.\\
We  use 
$\theta := 7e\pi/(Rx_0)$, $\theta_0:=\pi/(4Rx_0)$, $\theta_1 := 7e^2\pi/(2Rx_0)$, $Q_q:= 2^{k_q}\left((p/2)^p/(p!q) + \indic\{q=\infty\}\right)$, 
for all $k, l\geq 0$, $N\geq 1$, $f_{\alpha,\mt{\beta}} \in \mathcal{H}^{q,\phi,\omega}_{w, W}(l,M)$, 
\begin{equation}\label{eq:Nk}
\left(N+l\right)^{k} \leq ((l+1)N)^{k}, \ \int_{\epsilon\leq \abs{t}\leq T}L(t)dt \leq (2\pi)^{p+1}l^2. 
\end{equation}
\noindent {\bf Part 1.} 
The Plancherel and Chasles identities yield $ \left\|  \widehat{f}_{\alpha,\mt{\beta}}^{q,N, T, \epsilon}   -f_{\alpha,\mt{\beta}} \right\|^2_{L^2\left(1\otimes W^{\otimes p}\right)} \leq ( K_1 + K_2)/(2\pi) $. By the Jensen inequality, we have $K_1 \leq  4  \sum_{j=1}^4 \left\| R_j\right\|_{L^2\left(1\otimes W^{\otimes p}\right)}^2$ and, using \eqref{eq:interpol} for the first display and Lemma \ref{lem:interp_error} for the second, 
\begin{align}
&K_1+K_{2} \notag \\
&\leq K_1+ \int_{\R^p} 2 (1+ C_0(\underline{a}\epsilon))  \left\|\mathcal{F}_{1\mathrm{st}}\left[f_{\alpha,\mt{\beta}}\right](\star,\mt{b})-\mathcal{P}_{\underline{a}}\left[\mathcal{F}_{1\mathrm{st}}\left[f_{\alpha,\mt{\beta}}\right](\cdot,\mt{b})\right](\star)   \right\|_{L^2(\R)}^2   W^{\otimes p}(\mt{b}) d\mt{b} \notag  \\ & \quad   +  \int_{\R^p} 2C_0(\underline{a}\epsilon)\left\| \indic\left\{\abs{\star}\geq  \epsilon \right\}  \left( \widehat{F}_1^{q,N,T,0} - \mathcal{F}_{1\mathrm{st}}\left[f_{\alpha,\mt{\beta}}\right]\right)(\star,\mt{b})\right\|^2_{L^2(\R)} W^{\otimes p}(\mt{b}) d\mt{b}   \notag  \\
&\leq  K_1+4\pi (1+  C_0(\underline{a}\epsilon))\widetilde{w}(\underline{a}) \int_{\R^p} \norm{f_{\alpha,\mt{\beta}}(\cdot,\mt{b})}^2_{L^2(w	)}  W^{\otimes p}(\mt{b}) d\mt{b}   + 2C_0(\underline{a}\epsilon) K_1 \notag \\
&\leq  (1+ 2C_0(\underline{a}\epsilon)) K_1+4\pi(1+ C_0(\underline{a}\epsilon))  M^2\widetilde{w}(\underline{a})\le C\left(K_1+2\pi M^2\widetilde{w}(\underline{a})\right)\label{eq:K1K2}.
\end{align}

\noindent Using successively Proposition \ref{sec:upper:extension} and lemmas \ref{lem:def_est} and \ref{lem:lowerbound_weight10}, we have 
\begin{align}
\E\left[\left\| R_{1} \right\|_{L^2\left(1\otimes W^{\otimes p}\right)}^2 \right] 
= &  \int_{\epsilon\le |t|\le T}   \sum_{\abs{\mt{m}}_{q}\leq N(t)} \frac{\E\left[ \left|\widetilde{c}_{\mt{m}}(t) -c_{\mt{m}}(t)\right|^2   \right]}{ \left(\sigma_{\mt{m}}^{W,x_0t}\right)^2} dt  \notag \\
\leq&  \frac{c_{\mt{X}}}{(2\pi)^{p}n} \int_{\epsilon\le |t|\le T}   \abs{t}^p \nu_q^{W}(x_0t,N(t))  dt,\label{eq:R_0_w}
\end{align} 
also 
\begin{align}
\left\| R_2 \right\|_{L^2\left(1\otimes W^{\otimes p}\right)}^2 
&  \le  \int_{\epsilon\le |t|\le T} \left(\frac{x_0|t|}{2\pi}\right)^p \widetilde{\omega}^{q,W,x_0t}_{N(t)} \left\|S_0^N(t,\cdot_2)\right\|^2_{L^2([-1,1]^p)}  dt,\notag\\
\E\left[\left\|S_0^N(t,\cdot_2)\right\|^2_{L^2([-1,1]^p)}\right]&=\left\|S_1^N(t,\cdot_2)\right\|^2_{L^2([-1,1]^p)}+ \E\left[\left\|S_2^N(t,\cdot_2)\right\|^2_{L^2([-1,1]^p)} \right],\notag
\end{align}
\begin{align}
&\left\|S_1^N(t,\cdot_2)\right\|^2_{L^2([-1,1]^p)} \notag \\
&= \left\| \sum_{\abs{\mt{m}}_{q}\leq N(t)} g_{\mt{m}}^{W,x_0t} 
\braket{\mathcal{F}\left[f_{Y|\mt{X}=x_0\cdot_2}\right](t)\left(\Delta_f f_{\mt{X}|\mathcal{X}}\right)\left(x_0\cdot_2\right) , g_{\mt{m}}^{W,x_0t}}_{L^2([-1,1]^p)} \right\|^2_{L^2([-1,1]^p)} \notag \\
& \leq   \left\| 
\mathcal{F}\left[f_{Y|\mt{X}=x_0\cdot_2}\right](t) \left(
\Delta_f f_{\mt{X}|\mathcal{X}}
\right)\left(x_0\cdot_2\right)  \right\|^2_{L^2([-1,1]^p)}  \notag \\ 
&\leq   Z_{n_0}\left\| \mathcal{F}\left[f_{\alpha,\mt{\beta}}\right](t,x_0t\cdot_2) \right\|^2_{L^2([-1,1]^p)} 
\leq   Z_{n_0} \left(\frac{2\pi}{x_0|t|}\right)^p\left\| \mathcal{F}_{1\mathrm{st}}\left[f_{\alpha,\mt{\beta}}\right](t,\cdot_2) \right\|^2_{L^2(\R^p)}\label{eq:P3},
\end{align}
and, by independence 
and $\sum_{\abs{\mt{m}}_{q}\leq N}1 ={{N+p}\choose{p}}\indic\{q=1\}+(N+1)^p\indic\{q=\infty\}
\le (N+1)^p$, 
\begin{align}
\E\left[\left\|S_2^N(t,\cdot_2)\right\|^2_{L^2([-1,1]^p)} 
\right] 
= &     \sum_{\abs{\mt{m}}_{q}\leq N(t)} \frac{1}{n}\E\left[\left|Z_j^{\mt{m},t}  - \E\left[Z_j^{\mt{m},t} 
\right]\right|^2  
\right]\notag \\
\leq &     \sum_{\abs{\mt{m}}_{q}\leq N(t)} \frac{Z_{n_0}}{nx_0^{2p}}    \int_{\mathcal{X}}  \frac{1 }{f_{\mt{X}|\mathcal{X}}(\mt{x})} \left|  \overline{g_{\mt{m}}^{W,x_0t}}\left(\frac{\mt{x}}{x_0}\right) \right|^2 d\mt{x} \notag\\
\leq&\frac{(N(t)+1)^pc_{\mt{X}}Z_{n_0}}{nx_0^p}.\label{eq:P2}
\end{align} 
Collecting \eqref{eq:P3} and \eqref{eq:P2}, we obtain
\begin{align}
\E\left[\left\|R_2 \right\|_{L^2\left(1\otimes W^{\otimes p}\right)}^2  
\right]  
\leq   \frac{Z_{n_0}}{(2\pi)^p} \int_{\epsilon\leq \abs{t}\leq T }  \left(L(t) + \frac{c_{\mt{X}} (N(t)+1)^{p}  \abs{t}^p}{  n}\right) \widetilde{\omega}^{q,W,x_0t}_{N(t)}dt. \label{eq:R2_temp}
\end{align}
\noindent By Lemma \ref{lem:bc} and Proposition \ref{sec:upper:extension}, we have
\begin{align}
\left\| R_3\right\|_{L^2\left(1\otimes W^{\otimes p}\right)}^2 
\leq  \int_{\R}   \sum_{k> N(t)} \sum_{\abs{\mt{m}}_{q}=k} \left | b_{\boldsymbol{m}}(t) \right |^2  dt   
&\leq  \int_{\R}   \sum_{k> N(t)} \frac{\omega_{k}^2\theta_{q,k}^2(t) }{\omega_{N(t)+1}^2} 
dt \notag \\
&\leq \sup_{t\in \R} \frac{2\pi l^2}{\omega_{N(t)+1}^2}\label{eq:R_2_w}
\end{align}
\noindent 
and, by Proposition \ref{sec:upper:extension},
\begin{align}
\left\|R_4 \right\|^2_{L^2\left(1\otimes W^{\otimes p}\right)}
\leq    \sum_{k\in\N_0}   \int_{\abs{t} \geq T}     \sum_{\abs{\mt{m}}_{q}=k} \left | b_{\boldsymbol{m}}(t) \right |^2   dt 
&\leq      \sum_{k\in\N_0}\int_{\R}  \frac{\phi^2(\abs{t})}{\phi^2(T)} \theta_{q,k}^2(t)  dt \notag \\
&\leq   \frac{2\pi l^2}{ \phi^2(T)} . \label{eq:R3_w}
\end{align}
Thus we have \eqref{eq:START}.


\noindent{\bf Part 2.} We consider now all smoothness cases when $q\in\{1,\infty\}$.  Let  $t\neq0$ and $z>0$.  \eqref{eq:lower1} yields
\begin{equation}\label{eq:tildeW}
\widetilde{\omega}_N^{q,W,x_0t}\leq 2^p\left(1\bigvee  \frac{\theta (N+1)}{\abs{t}} \right)^{2k_qN}.
\end{equation} 
This yields, for all $N\geq 1$,
$\Delta_0^{W_{[-R,R]}}(t,N,n,z) \leq \Delta^{W_{[-R,R]}}(t,N,n,z)$, where 
\begin{align*}
\Delta^{W_{[-R,R]}}(\star,N,n,z) := & \left( 1 \bigvee \frac{\theta (N+1)}{\abs{\star}}\right)^{2k_qN} \frac{2}{\pi^{p+1}} \left(  \frac{Q_q c_{\mt{X}}  N^p\abs{\star}^p }{n}+   z \left(L(t)+ \frac{c_{\mt{X}} (N+1)^{p}\abs{\star}^p}{ n} \right)\right).
\end{align*}
Let $n_e$ be large enough to ensure $\underline{N} \geq (p+1)/k_q$. Using $N\leq \underline{N}$, $\epsilon \leq \theta  \leq \theta(\underline{N}+1) $, 
\begin{align}\notag 
\int_{\epsilon}^T \left( 1 \bigvee \frac{\theta (\underline{N}+1)}{t}\right)^{2k_q\underline{N}}t^p dt &=   \left(\theta (\underline{N}+1)\right)^{2k_q\underline{N}}  \int_{\epsilon}^{\theta(\underline{N}+1)} t^{p-2k_q\underline{N}}   dt +  \indic\{\theta (\underline{N}+1) \leq T\}\int_{\theta (\underline{N}+1)}^T t^p dt \\
&\leq    \frac{\epsilon^{p+1}}{2k_q\underline{N}-p-1} \left(\frac{\theta (\underline{N}+1)}{\epsilon}\right)^{2k_q\underline{N}} +  \frac{T^{p+1}}{p+1}, \label{eq:int_t0}
\end{align}
\begin{align}
\int_{\epsilon\leq \abs{t}\leq T} \left( 1 \bigvee \frac{\theta (\underline{N}+1)}{\abs{t}}\right)^{2k_q\underline{N}} L(t) dt  \leq    (2\pi)^{p+1}l^2 \left(\frac{\theta (\underline{N}+1)}{\epsilon}\right)^{2k_q\underline{N}} \quad \text{(by} \ \eqref{eq:Nk}), \label{eq:int_t1}
\end{align}
$n_e /n\leq 1$, and $n_e   v(n_0,\mathcal{E})/\delta(n_0)  \leq 1$, 
we have
\begin{align*}
&\int_{\epsilon\leq \abs{t}\leq T } \Delta^{W_{[-R,R]}}(t,N,n,Z_{n_0})dt \\
&\leq \frac{4c_{\mt{X}}\underline{N}^p }{\pi^{p+1}n_e}  \left(\frac{\epsilon^{p+1}}{k_q\underline{N}}\left(\frac{\theta (\underline{N}+1)}{\epsilon}\right)^{2k_q\underline{N}} +  \frac{T^{p+1}}{p+1}\right)\left(  Q_q +\frac{M_{\mathcal{E},\eta}2^p}{n} \right)   +    \frac{M_{\mathcal{E},\eta}2^{p+2}l^2}{n_e} \left(\frac{\theta (\underline{N}+1)}{\epsilon}\right)^{2k_q\underline{N}}  \\ 
&\leq \frac{\tau_0\underline{N}^{p-1}}{e^{2k_q}n_e} \left(\frac{\theta (\underline{N}+1)}{\epsilon}\right)^{2k_q\underline{N}}+ \frac{\tau_1 \underline{N}^{p-1} T^{p+1}}{n_e} ,  \\
& \tau_0:=\frac{e^{2k_q}4c_{\mt{X}}\theta^{p+1}}{\pi^{p+1}k_q} \left( Q_q +e^2M_{\mathcal{E},\eta}2^p\right)+  M_{\mathcal{E},\eta}2^{p+2}l^2,\ \tau_1 :=\frac{4c_{\mt{X}}}{\pi^{p+1}(p+1)}\left(Q_q +M_{\mathcal{E},\eta}2^p\right).
\end{align*} 
\eqref{eq:START}, $(\underline{N}+1)^{2k_q\underline{N}}\le e^{2k_q}\underline{N}^{2k_q\underline{N}}$, $\theta/\epsilon=K_{\underline{a}}(1)$,  $N+1\geq \underline{N}$, and the definition of $\underline{a}$, yield
\begin{align}\label{eq:START_2}
\mathcal{R}^W_{n_0,\sup} 
\le C\Bigg(&  \tau_0\underline{N}^{p-1}\frac{\left(\underline{N}K_{\underline{a}}(1)\right)^{2k_q\underline{N}}}{n_e} + \tau_1 \underline{N}^{p-1}\frac{T^{p+1}}{n_e} 
+\frac{8l^2+ M^2\indic\{w\neq W_{[-\underline{a}, \underline{a}]}\} }{\omega_{\underline{N}}^2} \Bigg).
\end{align}
The choices of $\underline{N}$ are such that the first and third terms 
have the same and largest order.\\ 
\noindent{\bf Proof of (T2.\ref{t_comp_1}).}  Let $n_e\ge e^e$ be large enough so that $(\ln(n_e)/\tau_2)^{\sigma ( (p+1)/s +2) + p-1} \leq n_e$,  where $\tau_2:=4k_q(2\sigma/\mu +1)\mathcal{W}(1/( 4k_q(2\sigma/\mu+1)))$. 
We have
\begin{align*}
\notag	2k_q\underline{N}\ln\left(\underline{N}w^{I}\left(\omega_{\underline{N}}^2\right)\right) + \ln\left(\omega_{\underline{N}}^2\right) + (p-1)\ln(\underline{N})
&= 2k_q\left(\frac{2\sigma}{\mu} +1\right)\underline{N}\ln\left(\underline{N}\right) + 2\sigma\ln\left(\underline{N}\right) + (p-1)\ln(\underline{N})\\
&\geq 2k_q\left(\frac{2\sigma}{\mu} +1\right)\underline{N}\ln\left(\underline{N}\right)
\end{align*}
and, for all $x\ge1/e$, $\mathcal{W}\left(x\ln(x)\right)=\ln(x)$. Using as well the definition of $\mathcal{W}$, this yields 
\begin{equation}\label{eq:ln1}
\underline{N} \leq \frac{\ln(n_e)}{ 4k_q(2\sigma/\mu +1)\mathcal{W}(\ln(n_e)/( 4k_q(2\sigma/\mu+1)))} \leq \frac{\ln(n_e)}{\tau_2}.
\end{equation}
Using \eqref{eq:ln1}, we have
$ T^{p+1} \underline{N}^{p-1}/n_e = \underline{N}^{\sigma(p+1)/s + p-1}/n_e \leq (\ln(n_e)/\tau_2)^{\sigma(p+1)/s+ p-1}/n_e \le (\ln(n_e)/\tau_2)^{-2\sigma} \le \underline{N}^{-2\sigma}=\omega_{\underline{N}}^{-2}$. Using the definition of $\underline{N}$ and \eqref{eq:START_2}, we obtain
\begin{align}
\mathcal{R}^W_{n_0,\sup} \le \frac{C}{\underline{N}^{2\sigma}}\Bigg(& \tau_0 + \tau_1
+8l^2+ M^2\Bigg).\label{eq_compact1}
\end{align}
We also have
\begin{equation*}\label{eq:ln}
\ln(n_e)=2k_q\underline{N}\ln\left(\underline{N}w^{I}\left(\omega_{\underline{N}}^2\right)\right) + \ln\left(\omega_{\underline{N}}^2\right) + (p-1)\ln(\underline{N}) \leq \left(2k_q\left(\frac{2\sigma}{\mu}+1\right)+\frac{2\sigma+p-1}{p+1}\right)\underline{N}\ln(\underline{N}),
\end{equation*}
hence $ \underline{N}\ln(\underline{N}) \ge \ln(n_e)/\tau_3$, $\tau_3:=2k_q\left(2\sigma/\mu+1\right)+(2\sigma+p-1)/(p+1)$. Similarly to \eqref{eq:ln1} and using for the second inequality, for all $x >0$, $\mathcal{W}(x) \leq  \ln(x +1 )$ (see Theorem 2.3 in \cite{hoorfar2007approximation}), we have
\begin{align}
&\underline{N}\geq  \frac{ \ln(n_e)}{\tau_3 \mathcal{W}\left(\ln(n_e)/\tau_3\right)}\geq   \frac{ \ln(n_e)}{\tau_3 \ln\left(\ln(n_e) + \tau_3\right)}
\ge\frac{ \ln(n_e)}{\tau_3\ln_2(n_e)  \left(1 + \ln(1+\tau_3/e) \right)}. \label{eq:W_lam}
\end{align}
This and \eqref{eq_compact1} yield the result.\\
\noindent{\bf Proof of (T\ref{theo:compact}.\ref{poly00}).}  Let $\tau_4:= \kappa(p+1)/\left(2s\left(k_q(\nu+1)+\kappa\right)\right)$. We have  
\begin{align}
2k_q\underline{N}\ln\left(\underline{N}w^{I}\left(\omega_{\underline{N}}^2\right)\right) + \ln\left(\omega_{\underline{N}}^2\right) 
\ge  2 \left(k_q(\nu+1)+\kappa\right)\underline{N} \ln\left(\underline{N} \right),
\label{eq:final_comp1}
\end{align}
hence 
$
\underline{N} \ln\left(\underline{N} \right) \leq \left( \ln(n_e) - (p-1)\ln(\underline{N})\right)/\left( 2 \left(k_q(\nu+1)+\kappa\right)\right)$ 
and $$\frac{T^{p+1}\underline{N}^{p-1}}{n_e} = \frac{e^{\kappa(p+1) \underline{N}\ln(\underline{N}+1)/s}\underline{N}^{p-1}}{n_e} \leq e^{\kappa(p+1)/s}n_e^{\tau_4 -1}\underline{N}^{(p-1)\left(1-\tau_4\right)}.$$
Because $s\ge\kappa(p+1)/(2k_q(1+\nu))$, we have $\tau_4 -1\leq  -\kappa/\left(\kappa + k_q(1+\nu) \right)$ and
\begin{align}
\mathcal{R}^W_{n_0,\sup} 
\le C\Bigg(
\frac{\tau_1e^{\kappa(p+1)/s}}{n_e^{\kappa/\left(\kappa +k_q(1+\nu)\right)}}\underline{N}^{(p-1)(1 - \tau_4)}+  \frac{\tau_0+  8l^2+ M^2\indic\{w\neq W_{[-\underline{a}, \underline{a}]}\}}{e^{2\kappa\underline{N}\ln(\underline{N}+1)}}\Bigg) .\label{eq_compact2}
\end{align}
Using again $\tau_4 -1\leq  -\kappa/\left(\kappa + k_q(1+\nu) \right)$, $\tau_5:=2(k_q(\nu +1) + \kappa)\ln(2)$, and
\begin{align}
\ln\left(n_e\right)-(p-1)\ln\left(\underline{N}\right)= 2k_q\underline{N}\ln\left(\underline{N}w^{I}\left(\omega_{\underline{N}}^2\right)\right) + \ln\left(\omega_{\underline{N}}^2\right) 
& \leq \frac{\tau_5}{\ln(2)}\underline{N}\ln\left(\underline{N}+1\right),
\end{align}
we obtain $e^{2\kappa\underline{N}\ln(\underline{N}+1)} \geq n_e^{\kappa/\left(\kappa +k_q(1+\nu)\right)}/\underline{N}^{\kappa(p-1)/ (k_q(\nu+1)+\kappa)}\geq n_e^{\kappa/\left(\kappa +k_q(1+\nu)\right)}/\underline{N}^{(p-1)(1-\tau_4)}$.  
We conclude 
because, by \eqref{eq:final_comp1}, $ \underline{N} \leq  \ln(n_e)/\tau_5$.\\
\noindent{\bf Proof of (T\ref{theo:compact}.\ref{poly01}).} It is derived from \eqref{eq_compact2} with $w= W_{[-\underline{a}, \underline{a}]}$ and $\nu = 0$.\\ 
\noindent{\bf Proof of (T\ref{theo:compact}.\ref{t_comp_3}).} By 
$\ln(n_e)\geq 2k_q\underline{N}\ln\left(\underline{N}w^{I}\left(\omega_{\underline{N}}^2\right)\right) + \ln\left(\omega_{\underline{N}}^2\right) \ge  \ln\left(\omega_{\underline{N}}^2\right)$, 
we have
\begin{align}
\label{eq:final_comp11}
(\underline{N}\ln(\underline{N}+1))^r \leq \frac{\ln(n_e)}{2\kappa},
\end{align}
hence, using the value of $T$ and $\underline{N}\le\underline{N}\ln\left(\underline{N}+1\right)/\ln(p+2)$,  
$$\frac{T^{p+1}\underline{N}^{p-1}}{n_e} = \frac{\kappa^{p+1} (\underline{N}\ln(\underline{N}+1))^{r(p+1)}\underline{N}^{p-1}}{\gamma^{p+1} n_e} \leq  \frac{ \ln(n_e)^{p+1+(p-1)/r}}{\kappa^{(p-1)/r} 2^{p+1+(p-1)/r}\gamma^{p+1} \ln(p+2)^{p-1} n_e}.$$
Moreover, because $\ln(n_e)^{p+1+(p-1)/r}$ is smaller than $ \varphi(n_e)$ by definition, 
\begin{align}\label{eq:end1}
\mathcal{R}^W_{n_0,\sup}  \le C\Bigg(&  \frac{\tau_1 \varphi(n_e) }{\kappa^{(p-1)/r}2^{p+1+(p-1)/r}\gamma^{p+1}\ln(p+2)^{p-1}n_e}
+\frac{\tau_0 + 8l^2+ M^2}{e^{2\kappa(\underline{N}\ln(1+\underline{N}))^r}} \Bigg).
\end{align}
We also have 
\begin{align}
2k_q\underline{N}\ln\left(\underline{N}w^{I}\left(\omega_{\underline{N}}^2\right)\right) + \ln\left(\omega_{\underline{N}}^2\right) +(p-1)\ln(\underline{N})
& \leq  2\kappa \left(\underline{N}\ln\left(\underline{N}+1\right)\right)^r \left( 1 + h(\underline{N}) \right),\label{eq:final_comp10}
\end{align}
where $h:= \left(k_q(1+\nu)\cdot + p-1\right)\ln\left(\cdot\right)/(\kappa\left(\cdot\ln\left(\cdot+1\right)\right)^{r}) $.
This yields, for $n_e$ large enough, 
\begin{align}
\exp\left( 2\kappa \left(\underline{N}\ln(\underline{N}+1)\right)^r\right) & \geq \exp\left( \frac{\ln(n_e)}{1 + h(\underline{N})}\right) 
= n_e \exp\left( \sum_{k=1}^{\infty} (-1)^k h(\underline{N})^{k}\ln(n_e)\right) .\label{eq:expinq}
\end{align}
By \eqref{eq:final_comp10}, 
we have $\underline{N}\ln(\underline{N}+1) \geq \ln(n_e)^{1/r}/d_0^{1/r}$. We obtain, by \eqref{eq:final_comp11} for the second inequality,
\begin{align*}
& h\left(\underline{N}\right)\leq \frac{(k_q+1)p/(p+1)+k_q-1 + k_q\nu}{\kappa (\underline{N}\ln(\underline{N}+1))^{r-1}} \leq \frac{((k_q+1)p/(p+1)+k_q-1 + k_q\nu)d_0^{1-1/r}}{\kappa \ln(n_e)^{1-1/r}} ,\\
& h\left(\underline{N}\right) \geq \frac{k_q(1+\nu)}{\kappa (1+1/((p+1)\ln(p+1)))^{r}(\underline{N}\ln(\underline{N}))^{r-1}}\geq \frac{k_q(1+\nu)(2\kappa)^{1-1/r}}{\kappa (1+1/((p+1)\ln(p+1)))^{r}\ln(n_e)^{1-1/r}},
\end{align*}
and we conclude using that, for $n_e$ large enough so that the remainder below is smaller in absolute value than a converging geometric series,  
\begin{equation*}\label{eq:phi1}
\exp\left(\sum_{k=1}^{\infty} (-1)^k h(\underline{N})^{k}\ln(n_e) \right)\geq \exp\left(\sum_{k=1}^{k_0}  (-1)^k  d_k\ln(n_e)^{(1/r-1)k+1}+ O(1)\right).\quad\square
\end{equation*}

\noindent{\bf Part 3.} Let $q\in\{1,\infty\}$. Let $t\neq0$ and $z>0$. By \eqref{eq:rho},  \eqref{eq:rho1}, and Proposition \ref{pdebut} \eqref{pdebutii},   we have, for $q\in{1,\infty}$ and $|t| \leq \theta_0$
\begin{align}
\widetilde{\omega}^{q,W,x_0t}_{N} \leq& \left( \frac{e\pi}{2} \right)^{2p} \exp\left( 2k_q\ln\left(\frac{\theta_1}{|t|}\right)N \right)  \indic\left\{ |t| \leq\theta_0 \right\} \label{eq:tildech}\\ 
&+  2^p \exp\left( \frac{2k_q\theta_0 (N+k_q')}{\abs{t}}\right)\indic\left\{ |t| > \theta_0\right\} .\notag
\end{align}
For all $N \geq 1$, we have
$\Delta_0^{\cosh(\cdot/R)}(t,N,n,z) \leq \Delta^{\cosh(\cdot/R)}(t,N,n,z)$, where
\begin{align}
&\Delta^{\cosh(\cdot/R)}(\star,N,n,z)\notag\\
&:=\frac{2}{\pi}\left(\frac{\pi e^2}{8}\right)^p\left(\frac{2^{p/q}Q_q c_{\mt{X}}N^{(p-1)/q}\abs{\star}^p}{n}+z\left(L(\star)+    \frac{c_{\mt{X}} (N+1)^{p}\abs{\star}^p}{ n}\right) 
\right)\left(\frac{\theta_1}{|\star|}\right)^{2k_qN}\indic\left\{ |\star| \le  \theta_0\right\}\notag\\
&\quad+\frac{2}{\pi^{p+1}}\left(\frac{2^{p/q}Q_q c_{\mt{X}}N^{(p-1)/q}\abs{\star}^{p+k_q}}{\left(4\theta_0/e\right)^{k_q}n}+z\left(L(\star)+\frac{c_{\mt{X}} (N+1)^{p}\abs{\star}^p}{n}\right)\right)\exp\left(\frac{2k_q\theta_0 (N+k_q')}{\abs{\star}}\right)\indic\left\{ |\star| >  \theta_0\right\}.\notag
\end{align}
Let $n_e\ge e^e$ be large enough so $\underline{N} \geq (p+2)/(2k_q)$. We have 
\begin{align}
\int_{\theta_0}^{T} t^{p+k_q}  e^{2k_q\theta_0(\underline{N}+k_q')/t} dt &= \int_{1/T}^{1/\theta_0}\frac{ e^{2k_q\theta_0(\underline{N}+k_q')u}}{u^{p+2+k_q}} du \leq \frac{\theta_0^{p+1+k_q} e^{2k_q(\underline{N}+k_q')}}{2k_q(\underline{N}+k_q')}\le \frac{\theta_0^{p+1+k_q} e^{2k_q(\underline{N}+k_q')}}{p+2(1+k_qk_q')}.\label{eq:ch0}
\end{align}
Then, using $N\leq \underline{N}$, $n_e /n\leq 1$, $n_e   v(n_0,\mathcal{E})/\delta(n_0)  \leq 1$, and $\epsilon \leq \theta_0 $, for the first display,
\begin{align}\label{eq:ch1} 
\int_{\epsilon}^{\theta_0} t^p \left(\frac{\theta_1}{t}\right)^{2k_q\underline{N}} dt 
&\leq    \frac{\epsilon^{p+1}}{2k_q\underline{N}-p-1} \left(\frac{\theta_1}{\epsilon}\right)^{2k_q\underline{N}}, 
\end{align}
\begin{align*}\notag 
&\int_{\epsilon \leq \abs{t}\leq \theta_0}  \left(\frac{\theta_1}{\abs{t}}\right)^{2k_q\underline{N}} L(t) dt +  \int_{\theta_0 \leq \abs{t}\leq T}  e^{2k_q\theta_0(\underline{N}+k_q')/\abs{t}} L(t) dt\\
& \leq    (2\pi)^{p+1}l^2 \left(  \left(\frac{\theta_1}{\epsilon}\right)^{2k_q\underline{N}} \indic\{\epsilon < \theta_0\}\bigvee  e^{2k_q(\underline{N}+k_q')} \right),
\end{align*}
and \eqref{eq:ch0} for the second display, we obtain  
\begin{align}&\int_{\epsilon\leq \abs{t}\leq T }\Delta^{\cosh(\cdot/R)}(t,N,n,Z_{n_0}) \notag\\
&\leq \frac{2^{2+p/q}Q_q c_{\mt{X}} \underline{N}^{(p-1)/q}}{\pi n_e} \left(\frac{\pi e^2}{8} \right)^p \int_{\epsilon}^{\theta_0} t^p \left(\frac{\theta_1}{t}\right)^{2k_q\underline{N}} dt \notag \\
&+\frac{M_{\mathcal{E},\eta}}{n_e}\left( \left(\frac{\pi e^2}{8} \right)^p  \bigvee \frac{1}{\pi^p} \right)   2^{p+2}\pi^{p}l^2  \left(  \left(\frac{\theta_1}{\epsilon}\right)^{2k_q\underline{N}} \indic\{\epsilon < \theta_0\}\bigvee  e^{2k_q(\underline{N}+k_q')} \right) \notag\\
&+ \frac{4 c_{\mt{X}}}{\pi n_e } \left(\frac{\pi e^2}{8} \right)^p \frac{2^{p} M_{\mathcal{E},\eta}  \underline{N}^{p}}{n} \int_{\epsilon}^{\theta_0} t^p \left(\frac{\theta_1}{t}\right)^{2k_q\underline{N}} dt +  \frac{2^{p/q}Q_q c_{\mt{X}}\underline{N}^{(p-1)/q}e^{k_q} }{\pi^{p+1}4^{k_q-1}\theta_0^{k_q} n_e} \int_{\theta_0}^{T} t^{p+k_q}  e^{2k_q\theta_0(\underline{N}+k_q')/t} dt \notag\\
& + \frac{4 c_{\mt{X}}}{\pi^{p+1} n_e }  \frac{2^{p} M_{\mathcal{E},\eta}  \underline{N}^{p}}{n}\int_{\theta_0}^{T} t^{p}  e^{2k_q\theta_0(\underline{N}+k_q')/t} dt \notag\\
& \leq G_1\left(\frac{\underline{N}^{k_q}}{n}\right) \frac{ \underline{N}^{(p-1)/q}}{n_e} \left(\frac{\theta_1}{\epsilon}\right)^{2k_q\underline{N}} \indic\{\epsilon < \theta_0\}  + G_2\left(\frac{\underline{N}^{k_q}}{n}\right)\frac{\underline{N}^{(p-1)/q}}{n_e}  e^{2k_q(\underline{N}+k_q')} 
, \label{eq:cosh_up}\\
&G:= 4\left(\frac{2k_q}{p+2}\right)^{(p-1)/q} M_{\mathcal{E},\eta}  \left(\frac{\pi e^2}{4} \right)^p  \pi^{p}l^2, \ G_1 :=  \frac{4 c_{\mt{X}}\theta_0^{p+1}}{\pi}\left( 2^{p/q}Q_q  + 2^{p}M_{\mathcal{E},\eta} \cdot \right)\left(\frac{\pi e^2}{8} \right)^p + G,  \notag\\
& G_2:= \frac{4 c_{\mt{X}}\theta_0^{p+1+k_q}}{(p+2(1+k_qk_q'))\pi^{p+1}}\left(2^{p/q}Q_q\left(\frac{e}{4\theta_0}\right)^{k_q} + 2^{p} M_{\mathcal{E},\eta} \cdot \right) +G. \notag 
\end{align} 

\noindent{\bf Proof of (T\ref{theo:non_compact}.\ref{t_noncomp_1}).} We have, using $K_{\underline{a}}(e) =w^I\left(\omega_{\underline{N}}^2\right)$, 
\begin{align*}
\ln(n_e) - \frac{p-1}{q}\ln(\underline{N}) =2k_q\underline{N}\ln\left(w^{I}\left(\omega_{\underline{N}}^2\right)\right) + \ln\left(\omega_{\underline{N}}^2\right) 
&= \frac{4k_q\sigma}{\mu} \underline{N}\ln\left(\underline{N}\right) + 2\sigma\ln\left(\underline{N}\right) \geq\frac{4k_q\sigma}{\mu} \underline{N}\ln\left(\underline{N}\right),  
\end{align*}
hence, for $n_e$ large enough so that $\ln(\underline{N})\geq \mu/(k_q\sigma)$, 
\begin{equation}\label{eq:mu_4}
\underline{N} \leq \frac{\ln(n_e)}{4k_q\sigma\ln(\underline{N})/\mu}\leq \frac{\ln(n_e)}{4k_q} .
\end{equation}
Thus, using \eqref{elog},  we have $\underline{N}/n\leq 1/(4k_qe)$. This yields, using \eqref{eq:cosh_up}, \eqref{eq:START}, $\theta_1/\epsilon=w^I\left(\omega_{\underline{N}}^2\right)$,  $N+1\geq \underline{N}$, and the definition of $\underline{a}$,
\begin{align*}
\mathcal{R}^W_{n_0,\sup} 
\le C\Bigg(&  G_1\left(\frac{k_q^{k_q-1}}{4e^{k_q}}\right)\underline{N}^{p-1}\frac{\left(w^I\left(\omega_{\underline{N}}^2\right)\right)^{2k_q\underline{N}}}{n_e} + G_2\left(\frac{k_q^{k_q-1}}{4e^{k_q}}\right)\underline{N}^{p-1}  \frac{e^{2k_q(\underline{N}+k_q')}}{n_e}
+\frac{8l^2+ M^2 }{\omega_{\underline{N}}^2} \Bigg).
\end{align*}
By \eqref{eq:mu_4},  we obtain $\omega_{\underline{N}}^2\underline{N}^{(p-1)/q}e^{2k_q(\underline{N}+k_q')}/n_e \leq  \ln(n_e)^{2\sigma+(p-1)/q}e^{2k_qk_q'}/(4^{2\sigma+p-1}\sqrt{n_e})$. Thus, using \eqref{elog}, we have $\underline{N}^{(p-1)/q}e^{2k_q(\underline{N}+k_q')}/n_e \leq e^{2k_qk_q'}((2\sigma+(p-1)/q)/(2e))^{2\sigma+(p-1)/q}/\omega_{\underline{N}}^2$ and using the definition of $\underline{N}$, 
\begin{align}
\mathcal{R}^W_{n_0,\sup} \le \frac{C}{\underline{N}^{2\sigma}}\Bigg(&G_1\left(\frac{k_q^{k_q-1}}{4e^{k_q}}\right)+ G_2\left(\frac{k_q^{k_q-1}}{4e^{k_q}}\right)e^{2k_qk_q'}\left(\frac{2\sigma+(p-1)/q}{2e}\right)^{2\sigma+(p-1)/q}+8l^2+ M^2\Bigg).\label{eq_non_compact1}
\end{align}
We also have
\begin{equation*}
\ln(n_e)=2k_q\underline{N}\ln\left(w^{I}\left(\omega_{\underline{N}}^2\right)\right) + \ln\left(\omega_{\underline{N}}^2\right) + \frac{p-1}{q}\ln(\underline{N})\leq\left( 2\left(\frac{2k_q}{\mu}+1\right)\sigma + \frac{2(p-1)}{q(p+2)}\right)\underline{N}\ln(\underline{N}),
\end{equation*}
hence $ \underline{N}\ln(\underline{N}) \geq \ln(n_e)/\tau_6$, $\tau_6:=2\left(2k_q/\mu+1\right)\sigma + 2(p-1)/(q(p+2))$. Similarly to \eqref{eq:W_lam}, we have $\underline{N}\geq \ln(n_e)/ (\tau_6 \left(1 + \ln(1+\tau_6/e) \right)\ln_2(n_e) )$, which yields the result with \eqref{eq_non_compact1}.\\ 
\noindent{\bf Proof of (T\ref{theo:non_compact}.2).} 
Because $K_{\underline{a}}(e)=1$ then $ 2(k_q+\kappa)\underline{N} + (p-1)\ln(\underline{N}) \geq 2(k_q + \kappa)\underline{N} $, we obtain $\underline{N} \leq \ln(n_e)/(2(k_q+\kappa))$. Thus using $n\geq n_e$ and \eqref{elog}, we have $G_2(\underline{N}^{k_q}/n)\leq G_2(k_q^{k_q}/(2(k_q+\kappa)e)^{k_q})$. Using \eqref{eq:START},  $  w=W_{\mathcal{A}} $, \eqref{eq:cosh_up},  $\epsilon =\theta_0$, yield
\begin{align*}
\mathcal{R}^W_{n_0,\sup} 
\le C\Bigg(& G_2\left(\frac{k_q^{k_q}}{2(k_q+\kappa)e^{k_q}}\right)\underline{N}^{(p-1)/q} \frac{e^{2k_q(\underline{N}+k_q')}}{n_e} +\frac{8l^2}{\omega_{\underline{N}}^2} \Bigg).
\end{align*}
We conclude using the definition of $\underline{N}$, which yields  $\underline{N}^{(p-1)/q}e^{2k_q\underline{N}}/n_e= n_e^{-\kappa/(k_q+\kappa)} \underline{N}^{\kappa (p-1)/(q(1+\kappa))}$
and $\omega_{\underline{N}}^{-2}=e^{-2\kappa \underline{N}}=n_e^{-\kappa/\left(k_q+ \kappa\right)}\underline{N}^{\kappa (p-1)/(q(k_q+\kappa))}$. \hfill $\square$

\subsection{Auxiliary lemmas and proof of Theorem \ref{cor:adp_rate}}

The proof of Theorem \ref{cor:adp_rate} uses several auxiliary lemmas. Lemmas \ref{adp2} and \ref{adp1} are  particularly important.\\

Let $\mathcal{N}_n$ be the set of functions $N\in \N_0^{\R}$ such that, for all $t\in \R\setminus(-\epsilon,\epsilon)$, $N(t) \in \{0,\dots,N_{\max,q}^{W}\}$ and $\Pi(n, Z_{n_0},T_{\max},N^W_{\max,q})$. 
Let, for all $t\neq0$ and $N\in \N_0^{\R}$,  $\Delta_{\mt{m}}:=\widehat{c}_{\mt{m}} - \widetilde{c}_{\mt{m}}$, $\widetilde{\Delta}_{\mt{m}}:=\widetilde{c}_{\mt{m}}- c_{\mt{m}}$,  
\begin{align}&\Xi\left(t, N\right):= \sum_{\abs{\mt{m}}_{q}> N} \abs{\frac{c_{\mt{m}}(t)}{\sigma^{W,x_0t}_{\mt{m}}}}^2,\ 
S_{1}\left(t, N\right):=\sum_{\abs{\mt{m}}_{q} \leq  N } \abs{ \frac{\E\left[\Delta_{\mt{m}}(t)
		\right]}{\sigma^{W,x_0t}_{\mt{m}}}}^2, \notag\\
&S_{2}\left(t, N\right) := \sum_{\abs{\mt{m}}_{q} \leq N} \abs{ \frac{\Delta_{\mt{m}}(t) - \E\left[ \Delta_{\mt{m}}(t)
		\right]}{\sigma^{W,x_0t}_{\mt{m}}}}^2, \  S_{3}\left(t,N\right) 
:= \sum_{\abs{\mt{m}}_{q} \leq N}\abs{\frac{\widetilde{\Delta}_{\mt{m}}(t)}{\sigma^{W,x_0t}_{\mt{m}}}}^2,\notag\\
& K_n(t):= H_W(t)  \left(N_{\max,q}^{W}+\frac{1}{2}\right)^{p}, \  L:=  \frac{1}{42} \sqrt{\frac{2x_0^p}{c_{\mt{X}}}}, \notag\\
&\Psi_{0,n}(t):= \exp\left(- \frac{p_n}{6} \right) + \frac{294c_{\mt{X}}  K_n^2(t)}{x_0^p n} \exp\left(-  \frac{L \sqrt{p_nn} }{ K_n(t) } \right),\notag\\
& \widetilde{B}\left(\widehat{N}\right) := \E\left[ \sup_{T' \in \mathcal{T}_n}\int_{\epsilon\leq \abs{t} \leq T'} \left( S_{3}\left(t, \widehat{N}(t)\right)  - \frac{\Sigma\left(t,\widehat{N}(t)\right)}{2(2+\sqrt{5})} \right)_+  dt 
\right],\notag\\
&\Pi(n, Z_{n_0},T_{\max},N^W_{\max,q}) :=  Z_{n_0} \int_{\epsilon\leq \abs{t} \leq T_{\max}}  \Psi_n(t)  dt +  \Pi_1\left(n,T_{\max},N^W_{\max,q}\right), \label{eq:pi} \\
&\Pi_1(n,T_{\max},N^W_{\max,q}) : =  \frac{96\left(1+2\sqrt{5}\right)c_{\mt{X}}  K_{\max}}{(2\pi)^p n}  \int_{\epsilon}^{T_{\max}}   \left(N_{\max,q}^{W}+1\right) t^p \nu_q^{W}\left(x_0t, N_{\max,q}^{W}\right)    \Psi_{0,n}(t)  dt, \notag\\ 
&\Psi_n :=  \left( 2+\frac{1}{\sqrt{5}} \right) \left( \left(\frac{2\pi}{x_0\abs{\star}}\right)^p  \widetilde{\omega}^{q,W,x_0t}_{N_{\max,q}^W} \left\|\mathcal{F}_{1\mathrm{st}}\left[f_{\alpha,\mt{\beta}}\right](\star,\cdot_2)\right\|^2_{L^2(\R^{p})}+ \left(\frac{ \abs{\star}}{2\pi}\right)^p\frac{c_{\mt{X}}\nu_q^{W}\left(x_0\star, N_{\max,q}^{W}\right)}{n} \right), \notag\\
&\widetilde{\Delta}^{W}_0(\star,N,n,z) :=\frac{c_{\mt{X}}|\star|^p}{n}\left(1 +2(1+2p_n)(1+c_1) \right) \nu_q^{W}(x_0\star,N ) 
+ z\left(L(\star) + \frac{ c_{\mt{X}} (N+1)^p\abs{\star}^p}{n} \right) \widetilde{\omega}^{q,W,x_0\star}_{N}, \label{eq:tildedelta}
\end{align}
where $H_W(t)$ is defined in Proposition \ref{eq:g_m_control}. 
For all $t\in[-T,T]\setminus[-\epsilon,\epsilon] $ and $N\in \N_0$, using \eqref{Young} with $c=\sqrt{5}$, we have
\begin{align} \label{eq_decom_adp} 
&\mathcal{L}_{q}^{W}(t,N)
\le  \Xi(t, N)  + \left(1+\frac{2}{\sqrt{5}}\right)\left( S_{1}(t,N) 
+ S_{2}(t,N)  \right) + (1+2\sqrt{5}) S_{3}(t,N).
\end{align}

\begin{lemma}\label{lem:adp} 
	\noindent For all $q\in \{1,\infty\}$, $0<\epsilon<1< T<T_{\max}=2^{K_{\max}}$, $t \in [-T,T]\setminus(-\epsilon,\epsilon)$, 
	and $N\in\{0,\dots,N_{\max,q}^{W}\}$, we have
	\begin{align}
	& \E\left[ S_1\left(t,\widehat{N}(t)\right) \right] \leq  Z_{n_0}\left(\frac{2\pi}{x_0|t|}\right)^p   \widetilde{\omega}^{q,W,x_0t}_{N_{\max,q}^{W}}  \left\|\mathcal{F}_{1\mathrm{st}}\left[f_{\alpha,\mt{\beta}}\right](t,\cdot_2) \right\|^2_{L^2(\R^{p})} \label{lem:eq1},\\
	&\E\left[  S_2\left(t,\widehat{N}(t)\right) \right]\leq   Z_{n_0}  \frac{ c_{\mt{X}} }{n}\left( \frac{|t|}{2\pi}\right)^{p} \nu_q^{W}(x_0t,N_{\max,q}^{W})  \label{lem:eq2},\\
	&\E\left[ \left( S_{3}(t,N)-\frac{\Sigma(t,N)}{2(2+\sqrt{5})}\right)_+\right]\leq48\frac{c_{\mt{X}}}{n}\left(\frac{\abs{t}}{2\pi}\right)^p\nu_q^{W}(x_0t,N) \Psi_{0,n}(t).\label{lem:eq3}
	\end{align} 
\end{lemma}

\noindent {\bf Proof of Lemma \ref{lem:adp}.}  
Let the parameters in the for all statement be given. 
\eqref{lem:eq1} follows from 
\begin{align*}
S_{1}\left(t,\widehat{N}(t)\right)&\leq \widetilde{\omega}^{q,W,x_0t}_{\widehat{N}(t)}  \left\|
\mathcal{F}\left[f_{Y|\mt{X}=x_0 \cdot}\right](t)\left(\Delta_f f_{\mt{X}|\mathcal{X}}
\right)\left(x_0\cdot\right)\right\|^2_{L^2([-1,1]^p)} \\
&\le Z_{n_0} \left(\frac{2\pi}{x_0|t|}\right)^p \widetilde{\omega}^{q,W,x_0t}_{N_{\max,q}^{W}} \left\|\mathcal{F}_{1\mathrm{st}}\left[f_{\alpha,\mt{\beta}}\right](t,\cdot_2) \right\|^2_{L^2(\R^{p})} \quad (\text{by} \  \eqref{eq:P3}).
\end{align*}
By  Lemma \ref{lem:lowerbound_weight10}, 
$\sum_{\abs{\mt{m}}_{q} \leq  N } \left(\sigma_{\mt{m}}^{W,x_0t}\right)^{-2}
\leq |x_0t|^p\nu_q^{W}(x_0t,N)/(2\pi)^{p}$, so we obtain \eqref{lem:eq2} by the following sequence of inequalities, which uses \eqref{eq:P2} for the second display,
\begin{align*}
\E\left[  S_2\left(t,\widehat{N}(t)\right)
\right]\leq &  \sum_{\abs{\mt{m}}_{q} \leq  N_{\max,q}^{W}}  \frac{ \E\left[\abs{\Delta_{\mt{m}}(t) - \E\left[ \Delta_{\mt{m}}(t) 
		\right]}^2 \right]}{\left(\sigma_{\mt{m}}^{W,x_0t}\right)^2}  \leq \frac{c_{\mt{X}} Z_{n_0}|t|^p\nu_q^{W}\left(x_0t,N_{\max,q}^{W}\right)}{(2\pi)^{p}n}. 
\end{align*}
To prove \eqref{lem:eq3}, we use
\begin{align}
&S_{3}(t,N)   =   \int_{\R^p}  \left| \widetilde{F}_1^{q,N,T,0}(t,\mt{b}) - F_1^{q,N,T,0}(t,\mt{b})  \right|^2  W^{\otimes p}(\mt{b})d\mt{b} 
= \sup_{u\in\mathcal{U}} \left| \nu_n^t(u)\right|^2,\notag\\
\notag \nu_n^t(u)
&:= \left\langle  \widetilde{F}_1^{q,N,T,0}(t,\cdot_2) - F_1^{q,N,T,0}(t,\cdot_2), u(\cdot)  \right\rangle_{L^2(W^{\otimes p})} 
= \frac{1}{n} \sum_{j=1}^n \left( f^t_u(Y_j,X_j)-\E\left[f^t_u(Y_j,X_j)\right]\right),\\
\notag  f^t_u (\star, \cdot) &:= \indic\left\{\cdot\in\mathcal{X}\right\} \frac{ e^{it\star }}{x_0^p f_{\mt{X}|\mathcal{X}}(\cdot)}\int_{\R^p} \sum_{\abs{\mt{m}}_{q}\leq N}  \overline{g ^{W,x_0t}_{\mt{m}}}\left(\frac{\cdot}{x_0}\right)\frac{1}{\sigma_{\mt{m}}^{W,x_0t}}  \varphi^{W,x_0t}_{\mt{m}}(\mt{b}) \overline{u}(\mt{b})W^{\otimes p}(\mt{b}) d\mt{b},
\end{align}
and 
$\mathcal{U}$ is a countable dense set of measurable functions of $\left\{u:\ \left\|u \right\|_{L^2\left(W^{\otimes p}\right)}=1\right\}$ 
and check the conditions of the Talagrand inequality given in Lemma \ref{Talagrand}  with $\eta =  p_n$ and  $\Lambda(p_n)=1$. 
For all $u\in \mathcal{U}$, the Cauchy-Schwarz inequality 
yield
\begin{align*}
\left\| f^t_u \right\|_{L^{\infty}( \R\times \mathcal{X})} & \leq c_{\mt{X}}\left(\frac{\abs{t}}{2\pi x_0}\right)^{p/2}\left\|\left( \sum_{\abs{\mt{m}}_{q}\leq N} 
\frac{\left|\overline{g ^{W,x_0t}_{\mt{m}}}\left(\cdot/x_0\right)\right|^2}{\rho_{\mt{m}}^{W,x_0t}}\int_{\R^p}  \left| \varphi^{W,x_0t}_{\mt{m}}(\mt{b})\right|^2  W^{\otimes p}(\mt{b}) d\mt{b}\right)^{1/2}  \right\|_{L^{\infty}(  \mathcal{X})}\\ 
&
\leq c_{\mt{X}} K_n(t)\left(\frac{\abs{t}}{2\pi x_0}\right)^{p/2} \sqrt{\nu_q^{W}(x_0t,N)}.
\end{align*}
By the Cauchy-Schwarz inequality and the computation leading to \eqref{eq:R_0_w}, we have
\begin{align*} 
\E\left[  \sup_{u \in \mathcal{U}}\left| \nu_n^t(u)\right|\right]^2  \leq 
\E\left[  \sup_{u \in \mathcal{U}}\left| \nu_n^t(u)\right|^2   \right]& \leq 
\E\left[\left\|  \widetilde{F}_1^{q,N,T,0}(t,\cdot_2) - F_1^{q,N,T,0}(t,\cdot_2) \right\|^2_{L^2(W^{\otimes p})} \right]\\
& \leq \frac{c_{\mt{X}}}{n}\left(\frac{\abs{t}}{2\pi}\right)^p\nu_q^{W}(x_0t,N)=\frac{\Sigma(t,N)}{8(2+\sqrt{5})(1+2p_n)}.
\end{align*}
Finally, by the Cauchy-Schwarz inequality and Proposition \ref{eq:g_m_control} for the second display and Lemma \ref{lem:lowerbound_weight10} for the third display, we have
\begin{align*}
\text{Var}\left(\mathfrak{R}(f^t_u(Y_j,X_j))\right) \vee \text{Var}\left(\mathfrak{I}(f^t_u(Y_j,X_j))\right)
& \leq \int_{\R\times \mathcal{X}} \left|  f^t_u (y,\mt{x}) \right|^2 f_{Y,\mt{X}}(y,\mt{x}) dy d\mt{x} \\
&\leq c_{\mt{X}} \left(\frac{\abs{t}}{2\pi}\right)^p \nu_q^{W}(x_0t,N).\quad\square
\end{align*}
\begin{lemma}\label{adp2}
	For all $\epsilon>0$, $q\in \{1,\infty\}$, and $ T \in \mathcal{T}_n $, we have
	\begin{align*}
	\mathcal{R}_{n_0}^{W}\left(\widehat{f}_{\alpha,\mt{\beta}}^{q,\widehat{N},\widehat{T},\epsilon},f_{\alpha,\mt{\beta}}\right) \leq&    \frac{C \left(2+\sqrt{5}\right)^2}{2\pi}\int_{\epsilon\leq \abs{t}}
	\E\left[\mathcal{L}_{q}^{W}\left(t,\widehat{N}(t), T\right) \right] + \frac{\E\left[  \indic\{\abs{t}\leq T\} \Sigma\left(t, \widehat{N}(t)\right)\right]}{2+\sqrt{5}} dt 
	\\ &  
	+ \frac{C2(2+\sqrt{5})^2 }{\pi}\Pi(n, Z_{n_0},T_{\max}, N^W_{\max,q})  + C M^2\widetilde{w}(\underline{a}).
	\end{align*} 
\end{lemma}

\noindent {\bf Proof of Lemma \ref{adp2}.} 
Let $\epsilon>0$, $q\in \{1,\infty\}$, and $T \in\mathcal{T}_n $.\\ Start from \eqref{eRisk}. Using, for all $T_1,T_2\ge\epsilon$, 
$R_{T_1}^{T_2} :=  \indic\left\{ \epsilon\le  \abs{\star} \right\}\left(\widehat{F}_1^{q,\widehat{N},T_1,0}- \widehat{F}_1^{q,\widehat{N},T_2\vee T_1,0}\right)(\star,\cdot_2)$ and  $ R^{T_1}:=  \indic\left\{ \epsilon\le\abs{\star} \right\}\left(\widehat{F}_1^{q,\widehat{N}, T_1,0}  -  \mathcal{F}_{\rm{1st}}\left[f_{\alpha,\mt{\beta}}\right]\right)(\star,\cdot_2)$, we have $R^{\widehat{T}}=R_{\widehat{T}}^{T}-R_{T}^{\widehat{T}}+R^{T}$ and $\left\| R^{\widehat{T}} \right\|^2_{L^2\left(W^{\otimes p}\right)}=\indic\left\{ \epsilon\le\abs{\star} \right\}\mathcal{L}_{q}^{W}\left(\star,\widehat{N}(\star), \widehat{T}\right)$.
Because 
\begin{equation}\label{eB2} B_2\left(T_1, \widehat{N}\right)  = \underset{ T'\in \mathcal{T}_n}{\max}\left(\int_{T_1\le \abs{t} \leq T_1\vee T'} \left\|R_{T_1}^{T'}(t,\cdot_2)  \right\|^2_{L^2\left(W^{\otimes p}\right)}  -\Sigma\left(t, \widehat{N}(t)\right)dt\right)_{+}, 
\end{equation}
we have $\E\left[  \left\|R_{T_1}^{T_2}\right\|^2_{L^2\left(1\otimes W^{\otimes p}\right)}     \right] 
\leq \E\left[B_2\left(T_1, \widehat{N}\right) \right] + \E\left[ \Sigma_2\left(T_2, \widehat{N}\right) \right]$ for possibly random $T_1$ and $T_2$ on $\mathcal{T}_n$. By \eqref{Young} with $c=\sqrt{5}$  and  \eqref{eq:choiceT}, we have
\begin{align*}
\notag &\E\left[ \left\| R^{\widehat{T}}\right\|^2_{L^2\left(1\otimes W^{\otimes p}\right)} \right]\leq   2(2+\sqrt{5}) \left( \E\left[B_2\left(T, \widehat{N}\right)   \right]+ \E\left[\Sigma_2\left(T, \widehat{N}\right)  \right]\right)  +  \left(1+\frac{2}{\sqrt{5}}\right)   \E\left[   \left\|R^T	 \right\|^2_{L^2\left(1\otimes W^{\otimes p}\right)}    \right].  
\end{align*}
Using, for all $T' \in \mathcal{T}_n$,  $ R_{T,1}^{T'} : = \widehat{F}_{1}^{q,\widehat{N},T \vee T',0} -  F_{1}^{q,\widehat{N},T \vee T', 0} $, $ R_{T,2}^{T'}  :=F_{1}^{q,\widehat{N},T,0} -  \widehat{F}_{1}^{q,\widehat{N},T,0}  $, and $ R_{T,3}^{T'} := F_{1}^{q,\widehat{N},T \vee T',0} -  F_{1}^{q,\widehat{N},T,0}  $, by \eqref{Young}, 
the objective function in \eqref{eB2} is 
smaller than
\begin{align*}
\int_{T\le \abs{t} \leq T\vee T'}\left( (2+\sqrt{5})\sum_{j=1}^2 \left\| R_{T,j}^{T'} (t,\cdot_2) \right\|^2_{ L^2\left(W^{\otimes p}\right)} +   \left(1+\frac{2}{\sqrt{5}}\right)  \left\| R_{T,3}^{T'}  (t,\cdot_2)\right\|^2_{ L^2\left(W^{\otimes p}\right)} -\Sigma\left(t, \widehat{N}(t)\right) \right)_{+}  dt . 
\end{align*}
Using that $ F_{1}^{q,\infty,\infty ,0} = \mathcal{F}_{1\mathrm{st}}\left[f_{\alpha,\beta}\right]$, we have, for all $t\in \R\setminus(-\epsilon,\epsilon)$, 
\begin{align*}
\left\| R_{T,3}^{T'} (t,\cdot_2) \right\|^2_{L^2\left(W^{\otimes p}\right) } & =  \indic\{T\leq \abs{t} \leq T \vee T'\} \sum_{ 0\leq   \abs{\mt{m}}_{q} \leq \widehat{N} } \left| \frac{c_{\mt{m}}(t)}{\sigma_{\mt{m}}^{W,x_0t} }\right|^2  \\
&\leq \left\| \left(F_{1}^{q,\widehat{N},T ,0} -  \mathcal{F}_{1\mathrm{st}}\left[f_{\alpha,\beta}\right]\right)(t,\cdot_2)  \right\|^2_{L^2\left(W^{\otimes p}\right)},
\end{align*} hence 
\begin{align*}  B_2\left(T, \widehat{N}\right)&
\leq    \underset{T' \in \mathcal{T}_n}{\max}\int_{T\le \abs{t} \leq T'} \left(2(2+\sqrt{5}) \left\|\left( \widehat{F}_{1}^{q,\widehat{N},T',0} -  F_{1}^{q,\widehat{N},T',0}\right)(t,\cdot_2) \right\|^2_{L^2\left(W^{\otimes p}\right)}    - 
\Sigma\left(t, \widehat{N}(t)\right)\right)_{+} dt  \\
& \quad+    \left(1+\frac{2}{\sqrt{5}}\right)  \int_{\epsilon \leq \abs{t}} \left\|R^T(t,\cdot_2) \right\|^2_{L^2\left(W^{\otimes p}\right)} dt.
\end{align*} 
Finally, we have
\begin{align*}
&\E\left[ \left\|R^{\widehat{T}}\right\|^2_{L^2\left(1\otimes W^{\otimes p}\right)}\right]\le 2(2+\sqrt{5})\E\left[ \Sigma_2\left(T, \widehat{N}\right)\right] +   (5+2\sqrt{5}) \left(1+\frac{2}{\sqrt{5}}\right)  \E\left[   \left\| R^T \right\|^2_{L^2\left(1\otimes W^{\otimes p}\right)}    \right]\\   
&+  4(2+\sqrt{5})^2  \E\left[ \underset{T' \in \mathcal{T}_n}{\max} \int_{T\le \abs{t} \leq T'} \left(\left\|\left( \widehat{F}_{1}^{q,\widehat{N},T',0} -  F_{1}^{q,\widehat{N},T',0}\right)(t,\cdot_2) \right\|^2_{L^2\left(W^{\otimes p}\right)}  -\frac{\Sigma\left(t, \widehat{N}(t)\right)}{2(2+\sqrt{5})}  \right)_{+}  \right]. 
\end{align*}
\noindent Using \eqref{eq_decom_adp}  and Lemma \ref{lem:adp}, we have  
\begin{align*}
&  \E\left[ \underset{T' \in \mathcal{T}_n}{\max}\int_{T\le \abs{t} \leq T'}\left( \left\|\left( \widehat{F}_{1}^{q,\widehat{N},T',0} -  F_{1}^{q,\widehat{N},T',0}\right)(t,\cdot_2) \right\|^2_{L^2\left(W^{\otimes p}\right)} -\frac{\Sigma\left(t,\widehat{N}(t)\right)}{2(2+\sqrt{5})} \right)_+ dt  \right]\notag  \\ 
&  \le	\E\left[ \underset{T' \in \mathcal{T}_n}{\max}\int_{\epsilon\leq \abs{t} \leq T'} \left(\sum_{\abs{\mt{m}}_{q}\leq \widehat{N}(t)} \left(\frac{\left| \widehat{c}_{\mt{m}}(t) - c_{\mt{m}}(t)\right|}{\sigma^{W,x_0t}_{\mt{m}}}\right)^2 - \frac{\Sigma\left(t,\widehat{N}(t)\right)}{2(2+\sqrt{5})}    \right)_+ dt     \right]   \notag \\
&  \leq  \left(1+2\sqrt{5}\right) \widetilde{B}\left(\widehat{N}\right) + Z_{n_0}  \int_{\epsilon\leq \abs{t} \leq T_{\max}}   \Psi_n(t)  dt .
\end{align*} 
Considering the first term of the last inequality and using \eqref{lem:eq3} for the second display yields 
\begin{align*}
\widetilde{B}\left(\widehat{N}\right) &\leq \E\left[  \sum_{T' \in \mathcal{T}_n} \int_{\epsilon\leq \abs{t} \leq T'}\left( S_{3}\left(t, \widehat{N}(t)\right) -  \frac{\Sigma\left(t,\widehat{N}(t)\right)}{2(2+\sqrt{5})} \right)_+  dt     \right]     \\
& \leq  \sum_{T' \in \mathcal{T}_n}  \int_{\epsilon}^{T'} \sum_{0 \leq N \leq N_{\max,q}^{W}} 96\frac{ c_{\mt{X}}   }{n}\left(\frac{t}{2\pi}\right)^p  \nu_q^{W}(x_0t,N)  \Psi_{0,n}(t)  dt   \\
& \leq  \frac{96c_{\mt{X}} K_{\max} }{ (2\pi)^pn}  \int_{\epsilon}^{T_{\max}}   \left(N_{\max,q}^{W} +1\right)t^p \nu_q^{W}\left(x_0t, N_{\max,q}^{W}\right)    \Psi_{0,n}(t)  dt.\quad\quad\quad\quad\quad\quad\quad\square
\end{align*} 

\begin{lemma}\label{adp1} For all $\epsilon >0 $, 
	$q\in \{1,\infty\}$, and $(T,N)\in\mathcal{T}_n\times \mathcal{N}_n$, 
	\begin{align*}
	&\int_{\epsilon\leq \abs{t}\leq T } \E\left[\mathcal{L}_{q}^{W}\left(t,\widehat{N}(t),T\right)\right]  + \frac{ \indic\{\abs{t}\leq T\}}{2+\sqrt{5}} \E\left[\Sigma\left(T,\widehat{N}(t)\right) 
	\right]dt \\& \leq  \left(2+\sqrt{5}\right)^2  \left(\int_{\epsilon\leq \abs{t}\leq T} 	
	\E\left[\mathcal{L}_{q}^{W}\left(t,N(t),T\right)\right] dt
	+ \frac{1+c_1}{2+\sqrt{5}}\E\left[\Sigma_2(T, N) 
	\right] \right) \\
	&\quad + \left(2+\sqrt{5}\right)^2 4\Pi(n,Z_{n_0},T_{\max},N_{\max,q}^W) .
	\end{align*}
\end{lemma}

\noindent{\bf Proof of Lemma} \ref{adp1}.
Let $ t\in [-T,T]\setminus(-\epsilon,\epsilon)$, $  N\in\{0,\dots,N_{\max,q}^{W}\} $, $T\in \mathcal{T}_n$. Using, for all $N_1,N_2\in \N_0$,  $\widetilde{R}_{N_1}^{N_2}(t,\cdot_2) :=  \left(\widehat{F}_1^{q,N_1, T,0} - \widehat{F}_1^{q,N_2\vee N_1, T,0}\right)(t,\cdot_2)$, 
we have $ \mathcal{L}_q^W\left(t,\widehat{N}(t),T\right) =  \widetilde{R}_{\widehat{N}(t)}^{N} - \widetilde{R}_{N}^{\widehat{N}(t)} + \mathcal{L}_q^W\left(t,N,T\right)$. Using \eqref{Young} yields
\begin{align*}	
\E\left[\mathcal{L}_{q}^{W}\left(t,\widehat{N}(t),T\right)\right]   \leq & (2+\sqrt{5})\left( \E\left[ \norm{ \widetilde{R}_{\widehat{N}(t)}^{N} (t,\cdot_2)}^2_{L^2\left(W^{\otimes p}\right)}    \right] + \E\left[ \norm{ \widetilde{R}_{N}^{\widehat{N}(t)} (t,\cdot_2)}^2_{L^2\left(W^{\otimes p}\right)}    \right] \right) \\
&+ \left(1+ \frac{2}{\sqrt{5}}\right)\E\left[ \mathcal{L}_{q}^{W}\left(t,N,T\right) \right] .
\end{align*} 
Because 
$B_1\left(t,N\right) = \underset{N'\in\N_0: \ N'\leq N_{\max,q}^{W}  }{ \max}\left( \sum_{N\leq \abs{\mt{m}}_{q}\leq N'\vee N	}\left(\abs{\widehat{c}_{\mt{m}}(t)}/\sigma_{\mt{m}}^{W,x_0t}\right)^2  - \Sigma\left(t, N'\right)\right)_{+}$, we have $\E\left[  \left\|\widetilde{R}_{N_1}^{N_2}(t,\cdot_2) \right\|^2_{L^2\left( W^{\otimes p}\right)}     \right] \leq \E\left[B_1\left(t, N_1\right) \right] + \E\left[ \Sigma\left(t, N_2\right) \right]$ for possibly random $N_1$ and $N_2$. Using $c_1 \geq 1+1/(2+\sqrt{5})^2$ and 
\eqref{eq:choiceN} yield
\begin{align*}
&\E\left[\mathcal{L}_{q}^{W}\left(t,\widehat{N}(t),T\right)\right] + \frac{1}{2+\sqrt{5}}\E\left[\Sigma\left(t,\widehat{N}(t)\right)\right]\\
&\leq    (2+\sqrt{5})  \left( 2 \E\left[B_1\left(t, N\right)   \right]+ (1+c_1)\E\left[\Sigma\left(t, N\right) \right]\right)  +  \left(1+\frac{2}{\sqrt{5}}\right) \E\left[ \mathcal{L}_{q}^{W}\left(t,N,T\right)  \right].
\end{align*}
By \eqref{Young} and, for all $N'\in \mathcal{N}_n$,  $ \widetilde{R}_{N,1}^{N'}(t, \cdot_2) := \left(\widehat{F}_1^{q,N\vee N',T,0} -  F_1^{q,N\vee N',T,0}\right) (t, \cdot_2)$, 
$ \widetilde{R}_{N,2}^{N'}(t, \cdot_2):=\left(F_{1}^{q,N,T,0} -  \widehat{F}_{1}^{q,N,T,0}\right)(t, \cdot_2) $, and $ \widetilde{R}_{N,3}^{N'}(t, \cdot_2): = \left(F_{1}^{q,N\vee N',T,0} -  F_1^{q,N,T,0}\right)(t, \cdot_2)$, we have that $B_1(t, N)$ is lower or equal than
\begin{align*}
 \underset{\substack{0\leq N' \leq N_{\max,q}^{W} \\  N'\in \N_0} }{\max}\left( (2+\sqrt{5}) \sum_{j=1}^{2}\left\|\widetilde{R}_{N,j}^{N'}(t,\cdot_2)\right\|^2_{ L^2\left(W^{\otimes p}\right)}  +  \left(1+\frac{2}{\sqrt{5}}\right)\left\|\widetilde{R}_{N,3}^{N'}(t,\cdot_2)\right\|^2_{ L^2\left(W^{\otimes p}\right)}  -\Sigma(t, N') \right)_{+} .
\end{align*}
Using $ F_{1}^{q,\infty,T ,0}(t, \cdot) = \mathcal{F}_{1\mathrm{st}}\left[f_{\alpha,\beta}\right](t, \cdot)$, we have
\begin{align*}
\left\| \widetilde{R}_{N,3}^{N'}(t,\cdot_2) \right\|^2_{L^2\left(W^{\otimes p}\right) } & =  \sum_{ N < \abs{\mt{m}}_{q} \leq N\vee N' } \left| \frac{c_{\mt{m}}(t)}{\sigma_{\mt{m}}^{W,x_0t} }\right|^2\leq \left\| \left(F_{1}^{q,N,T ,0} - \mathcal{F}_{1\mathrm{st}}\left[f_{\alpha,\beta}\right]\right)(t,\cdot_2)  \right\|^2_{L^2\left(W^{\otimes p}\right)},\\
B_1(t, N) 
& \leq  \underset{\substack{0\leq N' \leq N_{\max,q}^{W} \\  N'\in \N_0} }{\max} \left( 2(2+\sqrt{5}) \left\| \left(F_1^{q,N',T,0} -  \widehat{F}_1^{q,N',T,0} \right)(t,\cdot_2) \right\|^2_{L^2\left(W^{\otimes p}\right)}   - \Sigma(t, N') \right)_{+} \\
& \quad+  \left(1 + \frac{2}{\sqrt{5}}\right) \left\| \left(F_1^{q,N,T,0}  -  \mathcal{F}_{\rm{1st}}\left[f_{\alpha,\mt{\beta}}\right]\right)(t,\cdot_2)  \right\|^2_{L^2\left(W^{\otimes p}\right)}  .
\end{align*}
Finally, 
we have
\begin{align*}
&\E\left[\mathcal{L}_{q}^{W}\left(t,\widehat{N}(t),T\right)\right]  + \frac{1}{2+\sqrt{5}}\E\left[\Sigma\left(t,\widehat{N}(t)\right)\right] \\
&\leq  4(2+\sqrt{5})^2 \E\left[ \underset{0\leq N' \leq N_{\max,q}^{W}}{\max} \left( \left\| \left(\widehat{F}_1^{q,N',T,0} -  F_{1}^{q,N',T,0}\right)(t,\cdot_2) \right\|^2_{L^2\left(W^{\otimes p}\right)}  - \frac{\Sigma(t, N')}{2(2+\sqrt{5})}\right)_+   \right]  \\ 
& \hspace{0.2cm}  + (2+\sqrt{5})(1+c_1)\E\left[ \Sigma\left(t, N\right) \right] + \left(2+\sqrt{5}\right)^2\E\left[ \mathcal{L}_{q}^{W}\left(t,N,T\right)  \right].
\end{align*}
\noindent Using \eqref{eq_decom_adp} for the second display and Lemma \ref{lem:adp} for the third, we obtain  
\begin{align*}
&\E\left[ \underset{0\leq N' \leq N_{\max,q}^{W}}{\max} \left\| \left(\widehat{F}_1^{q,N',T,0} -  F_{1}^{q,N',T,0}\right)(t,\cdot_2) \right\|^2_{L^2\left(W^{\otimes p}\right)}  - \frac{\Sigma(t, N')}{2(2+\sqrt{5})}  \right]\notag  \\ 
& = \E\left[ \underset{0\leq N' \leq N_{\max,q}^{W}}{\max} \left(\sum_{\abs{\mt{m}}_{q}\leq N'} \left(\frac{\left|\widehat{c}_{\mt{m}}(t) - c_{\mt{m}}(t)\right| }{\sigma^{W,x_0t}_{\mt{m}}}\right)^2 -  \frac{\Sigma(t,N')}{2(2+\sqrt{5})} \right)_+   \right]  \notag \\
&\le  \left(1+2\sqrt{5}\right)\E\left[ \underset{0\leq N' \leq N_{\max,q}^{W}}{\max} \left(S_3(N',t)-  \frac{\Sigma(t,N')}{2(2+\sqrt{5})}  \right)_+   \right]    \notag  \\ 
&\quad + \left(2+\frac{1}{\sqrt{5}}\right) \E\left[ \underset{0\leq N' \leq N_{\max,q}^{W}}{\max}  \left( S_1(N',t) + S_2(N',t) \right)  \right] \notag \\
&\leq    \left(N_{\max,q}^{W}+1\right) \left(1+2\sqrt{5}\right)48\frac{c_{\mt{X}}  }{n}\left(\frac{\abs{t} }{2\pi}\right)^p  \nu_q^{W}\left(x_0t, N_{\max,q}^{W}\right)  \Psi_{0,n}(t) +  Z_{n_0} \Psi_{n}(t).\notag 
\end{align*} 
Hence the result.\hfill $\square$
\vspace{0.3cm}

\noindent Hereafter, let $(n,n_0) \in \N^2$ such that $ v(n_0,\mathcal{E})/\delta(n_0) \leq n^{-2}\ln(n)^{-p} $, $n \geq e^{7e^2/(2\pi)}$ large enough so that  $N_{\max,q}^{W}\geq(p+1)/k_q$. Using $\theta_1/\ln(n) \leq 2\theta_0$, let
$(\theta \indic\{W=W_{[-R,R]}\} + \theta_1 \indic\{W=\cosh(\cdot/R)\} )/\ln(n) \leq \epsilon\leq \theta \indic\{W=W_{[-R,R]}\} + 2\theta_0 \indic\{W=\cosh(\cdot/R)\}  $. Using the definition of $N_{\max,q}^{W}$ yields $N_{\max,q}^{W}\leq  \underline{N_{\max,q}^{W}}$ and $\epsilon \leq \theta$ yields $\underline{N_{\max,q}^{W_{[-R,R]}}}\ln\left(\underline{N_{\max,q}^{W_{[-R,R]}}}\right)\le \ln(n)/(2k_q)$. Then, using that, for all $x\ge1/e$, $\mathcal{W}\left(x\ln(x)\right)=\ln(x)$, and the definition of $\mathcal{W}$ for the bound on $N_{\max,q}^{W_{[-R,R]}}$ and else the definition of $N_{\max,q}^{\cosh(\cdot/R)}$, we have, for all  $t\neq0$,
\begin{align}\label{eq:N_max_indic}
&\quad N_{\max,q}^{W}\leq \frac{\ln(n)}{\tau_7}, \ \tau_7 :=2k_q\mathcal{W}\left(\frac{7e^2}{4\pi k_q}\right) \indic\{W=W_{[-R,R]}\} + 2k_q \indic\{W=\cosh(\cdot/R)\}. 
\end{align}
Using, for all $|t|\geq \epsilon$ and $N\ge1$, $2k_qN\ln(\theta(N+1)/\abs{t})\leq 2k_qN\ln(\theta N/\epsilon) + 2k_q$, we have $\left(\theta \left(N_{\max,q}^{W_{[-R,R]}}+1\right)/\abs{t}\right)^{2k_qN_{\max,q}^{W_{[-R,R]}}} \leq e^{2k_q}\left(\theta \underline{N_{\max,q}^{W_{[-R,R]}}}/\epsilon\right)^{2k_q\underline{N_{\max,q}^{W_{[-R,R]}}}}$. \eqref{eq:tildeW} and the definition of $N_{\max,q}^{W_{[-R,R]}}$ yield
\begin{align}
&\forall |t|\ge\epsilon,\ \widetilde{\omega}^{q,W_{[-R,R]},x_0t}_{N^W_{\max,q}} \leq 2^pe^{2k_q}n. \label{eq:control_N}
\end{align}

\begin{lemma}\label{lem:N_max_control}
	For all $q\in \{1,\infty\}$,  $\epsilon\leq T_{\max} \leq n^{\zeta_0}$, and $W\in\{W_{[-R,R]},\cosh(\cdot/R)\}$, we have
	\begin{align*}
	&\int_{\epsilon}^{T_{\max}} t^p\nu_q^{W}\left(x_0t, N_{\max,q}^{W}\right)   dt\leq  A_3^{W,q} \ln(n)^{a_0}n,\\
	& a_0 := p\indic\{W=W_{[-R,R]}\}+\frac{p-1}{q}\indic\{W=\cosh(\cdot/R)\},\\
	& A_3^{W_{[-R,R]}, q}:= \frac{2^p Q_q}{\tau_7^p(p+1)}\left(\frac{\theta^{p+1}e^{2k_q}}{k_q} + e^{(1-\zeta_0(p+1))7e^2/(2\pi)}\right), \\
	& A_3^{\cosh(\cdot/R),q}:= \frac{2^{1/q}Q_q}{k_q^{(p-1)/q}} \left(\left(\frac{e\pi}{2}\right)^{2p}  \indic\{\epsilon<\theta_0\} \frac{\theta_0^{p+1}}{p+1} + 2^{p}\left(\frac{e}{4}\right)^{k_q} \frac{ \theta_0^{p+1}  e^{2k_qk_q'}}{p+2(1+k_qk_q')}\right).
	\end{align*}
\end{lemma}
\noindent{\bf Proof.}
Let  $ W=W_{[-R,R]} $. 
Using the definition of $\nu^{W_{[-R,R]}}_q$, $N_{\max,q}^{W}\geq 1$, and \eqref{eq:Nk} for the first inequality, \eqref{eq:int_t0} and \eqref{eq:N_max_indic} for the second, and the arguments which yield \eqref{eq:control_N} for the last inequality, the result follows from 
\begin{align*}
\int_{\epsilon}^{T_{\max}} t^p\nu_q^{W}\left(x_0t, N_{\max,q}^{W}\right)   dt &\leq 2^p Q_q \left(\underline{N_{\max,q}^{W}}\right)^p \int_{\epsilon}^{T_{\max}} t^p \left(1 \bigvee \left(\frac{\theta \left(\underline{N_{\max,q}^{W}}+1\right)}{\abs{t}}\right)\right)^{2k_q\underline{N_{\max,q}^{W}}} dt\\
& \leq  \frac{2^p Q_q \ln(n)^p}{\tau_7^p}\left(\frac{\epsilon^{p+1}}{k_q\underline{N_{\max,q}^{W}}}\left(\frac{\theta(\underline{N_{\max,q}^{W}}+1)}{\epsilon}\right)^{2k_q\underline{N_{\max,q}^{W}}}  + \frac{T_{\max}^{p+1}}{p+1}\right)\\
& \leq  \frac{2^p Q_q \ln(n)^p}{\tau_7^p(p+1)}\left(\frac{\theta^{p+1}e^{2k_q}}{k_q}+ \frac{1}{n^{(p+1)\zeta_0-1}}\right) n. 
\end{align*}
Let now  $ W=\cosh(\cdot/R) $. Using the definition of $\nu_q^{W}$, we have for all $N\geq 1$ and $t\neq0$, 
\begin{align*}
\nu_q^{\cosh(\cdot/R)}\left(x_0t, N\right)=& 2^{p/q}Q_q N^{(p-1)/q} \left(\frac{e\pi}{2}\right)^{2p}  \left(\frac{\theta_1}{|t|}\right)^{2k_qN} \indic\{|t|\leq \theta_0\}\\
& + 2^{p(1+1/q)}Q_q N^{(p-1)/q}  \left(\frac{e|t|}{4\theta_0}\right)^{k_q} \exp\left(2\theta_0 k_q(N+k_q') \right)  \indic\{|t|> \theta_0\}.
\end{align*}
Because of  \eqref{eq:N_max_indic}, we have, for $t\neq0$, $ e^{2k_qN_{\max,q}^{\cosh(\cdot/R)}} \le n$. By definition of $\underline{N_{\max,q}^{\cosh(\cdot/R)}}$, when $\epsilon<\theta_0$, we also have, for $|t|\leq \theta_0$, $ (\theta_1/\epsilon)^{2k_qN_{\max,q}^{\cosh(\cdot/R)}} \leq n$. Then, using \eqref{eq:N_max_indic} for the first display and using \eqref{eq:ch0} and  \eqref{eq:ch1} for the second, the result follows from
\begin{align*}
\int_{\epsilon}^{T_{\max}} t^p\nu_q^{W}\left(x_0t, N_{\max,q}^{W}\right)   dt 
\leq& 2^{p/q}Q_q \left(\frac{\ln(n)}{2k_q}\right)^{(p-1)/q} \left(\frac{e\pi}{2}\right)^{2p}  \indic\{\epsilon<\theta_0\} \int_{\epsilon}^{\theta_0}t^{p}\left(\frac{\theta_1}{|t|}\right)^{2k_q\underline{N_{\max,q}^{W}}}dt\\
&+ Q_q \left(\frac{\ln(n)}{2k_q}\right)^{(p-1)/q} \frac{2^{p(1+1/q)}e^{k_q}}{(4\theta_0)^{k_q}}  \int_{\theta_0}^{T_{\max}}t^{p+k_q}e^{2\theta_0 k_q\left(\underline{N_{\max,q}^{W}}+k_q'\right)}dt\\
\leq &2^{p/q} Q_q \left(\frac{\ln(n)}{2k_q}\right)^{(p-1)/q} \left(\frac{e\pi}{2}\right)^{2p}  \indic\{\epsilon<\theta_0\} \frac{\theta_0^{p+1}}{p+1} \left(\frac{\theta_1}{\epsilon}\right)^{2k_q\underline{N_{\max,q}^{W}}}\\
&+ Q_q \left(\frac{\ln(n)}{2k_q}\right)^{(p-1)/q}  \frac{2^{p(1+1/q)}e^{k_q}}{(4\theta_0)^{k_q}}  \frac{\theta_0^{p+1+k_q}e^{2k_qk_q'}}{p+2(1+k_qk_q')}e^{2k_q \underline{N_{\max,q}^{W}}}\\
\leq & A_3^{\cosh(\cdot/R),q}\ln(n)^{(p-1)/q}n. \quad \square
\end{align*}

\begin{lemma}\label{lem:Z_control}
	For all 
	$q\in \{1,\infty\}$, $W\in\{W_{[-R,R]},\cosh(\cdot/R)\}$, and $(\epsilon\vee 1)\leq T_{\max} \leq n^{\zeta_0}$,  
	we have
	$ \Pi(n, Z_{n_0},T_{\max},N_{\max,q}) \leq ( A_0  + A_1)/n$, where 
	\begin{align*}
	&A_0 :=\frac{M_{\mathcal{E},\eta}}{n}\left(2+\frac{1}{\sqrt{5}}\right)\left( \left(\frac{4\pi^2}{\theta x_0}\right)^p b_0(2\pi)l^2 + \frac{2c_{\mt{X}}A_3^{W,q}}{(2\pi)^pe^{7e^2/(2\pi)}}\right),\\
	& b_0:=2^pe^{2k_q}\indic\{W=W_{[-R,R]}\} + \left(\frac{e\pi}{2}\right)^{2p}\indic\{W=\cosh(\cdot/R)\},\\
	&A_1:= \frac{96  \left(1+2\sqrt{5}\right) c_{\mt{X}} \zeta_0 A_3^{W,q}}{(2\pi)^{p} k_q \ln(2)\left(1/\tau_7+\pi/(7e^2)\right)^p }    \left(\left(\frac{a_0+2}{e\zeta_0}\right)^{a_0+2}+\frac{e^{1/b_1}294c_{\mt{X}} a_1^{2}}{x_0^{p} }  \left(\frac{a_0+2(p+1)}{e}\right)^{a_0+2(p+1)}\right),\\
	& a_1 := \left(\frac{1}{\tau_7}+\frac{\pi}{7e^2}\right)^p (H_0\indic\{W=W_{[-R,R]}\}    +H_1\indic\{W=\cosh(\cdot/R)\}  )^p \left( 1+x_0^2\right)^p, \\
	&  b_1:= \frac{L\sqrt{p_2}  }{  a_1} \left(\frac{ e  (1-4(1+\indic\{W=\cosh(\cdot/R)\})p\zeta_0)}{2p+3}\right)^{p+3/2}.
	\end{align*}
\end{lemma}

\noindent \textbf{Proof.} 
Let us show
\begin{align}
Z_{n_0} \int_{\epsilon\leq \abs{t} \leq T_{\max}}  \Psi_n(t)  dt \leq  \frac{A_0  }{n}\quad\text{and}  \quad  \Pi_1(n,T_{\max},N_{\max,q})  \leq \frac{A_1}{n}.\label{eq:ZC_1} 
\end{align}
Let $W=W_{[-R,R]}$. 
Using for the second display $ v(n_0,\mathcal{E})/\delta(n_0) \leq n^{-2}\ln(n)^{-p} $, \eqref{eq:control_N}, $\epsilon \geq \theta/\ln(n)$, \eqref{eq:Nk}, and Lemma \ref{lem:N_max_control}, we obtain
\begin{align*}
Z_{n_0}\int_{\epsilon\leq \abs{t} \leq T_{\max}} \Psi_n(t)dt
= &  Z_{n_0} \left( 2+\frac{1}{\sqrt{5}} \right)  
\left(\frac{2\pi}{x_0}\right)^p  \int_{\epsilon\leq \abs{t} \leq T_{\max}}  \widetilde{\omega}^{q,W,x_0t}_{N_{\max,q}^{W}} \abs{t}^{-p}\left\|\mathcal{F}_{1\mathrm{st}}\left[f_{\alpha,\mt{\beta}}\right](t,\cdot_2) \right\|^2_{L^2(\R^{p})}dt\\ & + Z_{n_0} \left( 2+\frac{1}{\sqrt{5}} \right)\frac{c_{\mt{X}}}{(2\pi)^pn}\int_{\epsilon\leq \abs{t} \leq T_{\max}} |t|^p\nu_q^{W}(x_0t,N_{\max,q}^{W}) dt \quad\\
\leq &\frac{M_{\mathcal{E},\eta}}{n^{2}\ln(n)^p}\left(2+\frac{1}{\sqrt{5}}\right)\left( \left(\frac{4\pi}{\theta x_0}\right)^p e^{2k_q} n \ln(n)^p  (2\pi)^{p+1}l^2+ \frac{2c_{\mt{X}}A_3^{W,q}\ln(n)^p}{(2\pi)^p}\right).
\end{align*}
Using $n\geq e^{7e^2/(2\pi)}$ and \eqref{elog} yield the first inequality in \eqref{eq:ZC_1}.
Similarly, by definition of $N^{\cosh(\cdot/R)}_{\max,q}$ and \eqref{eq:tildech}, we have $ \widetilde{\omega}^{q,\cosh(\cdot/R),x_0t}_{N^W_{\max,q}} \leq (e\pi/2)^{2p}n$. This and \eqref{eq:N_max_indic} yield the first inequality in \eqref{eq:ZC_1} for the other instances of $W$ and $q$. \\
By  \eqref{eq:N_max_indic}, we have
\begin{align}\label{eq:Kn_1}
\quad K_n(t) \leq \left(\frac{\ln(n)}{\tau_7}+\frac{1}{2}\right)^p T_{\max}^{2p} H_0^p \left(\frac{1}{T_{\max}^2}+x_0^2\right)^p \leq a_1  \ln(n)^pT_{\max}^{2p} .
\end{align}
We obtain, using $T_{\max}\leq n^{\zeta_0}$ for the third inequality and \eqref{elog} for the fourth, 
\begin{align*}
\frac{L\sqrt{p_nn}}{ K_n(t)}& \geq  \frac{L \sqrt{p_2\ln(n)n}}{ a_1 \ln(n)^{p} T_{\max}^{2p	}}\geq  \frac{L \sqrt{p_2} n^{(1 - 4\zeta_0p) /2}}{ a_1 \ln(n)^{p-1/2} }
\geq b_1 \ln(n)^2. 
\end{align*}
Using \eqref{eq:N_max_indic} for the first inequality,  Lemma \ref{lem:N_max_control}  for the second, and using the definition of $p_n$, $6(1+\zeta_0)\ln(n)>3$, \eqref{eq:Kn_1}, and $T_{\max}^{4p} \leq n^{4p\zeta_0}$, for the third, we have
\begin{align*}
&\Pi_1(n,T_{\max},N_{\max,q})\\
& \leq \frac{ 96  \left(1+2\sqrt{5}\right)  c_{\mt{X}}\zeta_0   \ln(n)^2 }{ (2\pi)^p k_q \ln(2) n} \left(\frac{1}{\tau_7}+\frac{\pi}{7e^2}\right)^p  \int_{\epsilon}^{T_{\max}} t^p\nu_q^{W}\left(x_0t, N_{\max,q}^{W}\right)   dt \sup_{\epsilon\le t\le T_{\max}}\Psi_{0,n}(t)  \\
& \leq \frac{96  \left(1+2\sqrt{5}\right) c_{\mt{X}} \zeta_0 A_3^{W,q}}{(2\pi)^{p} k_q \ln(2)n} \left(\frac{1}{\tau_7}+\frac{\pi}{7e^2}\right)^p   \sup_{\epsilon\le t \le T_{\max}} \ln(n)^{p+2} n \Psi_{0,n}(t) \\
& \leq \frac{96  \left(1+2\sqrt{5}\right) c_{\mt{X}} \zeta_0 A_3^{W,q}}{(2\pi)^{p} k_q \ln(2)n} \left(\frac{1}{\tau_7}+\frac{\pi}{7e^2}\right)^p  \left(\frac{\ln(n)^{p+2}}{n^{\zeta_0}}+\sup_{n>0}\left( e^{-b_1\ln(n)^2}n^2\right) \frac{294c_{\mt{X}} a_1^{2}}{x_0^{p} }\frac{\ln(n)^{3p+2}}{n^{2(1-2p\zeta_0)}}\right) .\ 
\end{align*}
Thus,  \eqref{elog}, $1-2p\zeta_0>1/2$, $\sup_{x>0}\left( e^{-b_1\ln(x)^2}x^2\right)=e^{1/b_1}$ yield the second inequality in \eqref{eq:ZC_1}.
We obtain similarly the bounds for the other instances of $W$ and $q$.
\hfill $\square$\vspace{0.3cm}

\noindent {\bf Proof of Theorem \ref{cor:adp_rate}.} Let $n,n_0$ such that $ v(n_0,\mathcal{E})/\delta(n_0) \leq n^{-2}\ln(n)^{-p} $,  
$T \in \mathcal{T}_n$, and $N \in \mathcal{N}_n$. The proof  of Theorem \ref{cor:adp_rate} has two parts. First, we prove 
\begin{align} \notag 
&\mathcal{R}_{n_0}^{W}\left(\widehat{f}_{\alpha,\mt{\beta}}^{q,\widehat{N},\widehat{T},\epsilon},f_{\alpha,\mt{\beta}}\right)\\
&\leq\frac{C2 \left(2+\sqrt{5}\right)^4}{\pi} \left(\int_{\epsilon\leq \abs{t}\leq T }  \frac{\widetilde{\Delta}^{W}_0(t,N(t),n,Z_{n_0})}{(2\pi)^p} dt + \sup_{t\in [-T,T]} \frac{ 2\pi l^2}{ \omega_{N(t)+1}^2} +\frac{ 2\pi l^2}{\phi(T)^2}\right) \notag
\\ & \quad + CM^2\widetilde{w}(\underline{a}) + \frac{2(2+\sqrt{5})^2 C\left(1+\left(2+\sqrt{5}\right)^2\right)(A_0+A_1)}{\pi n}, \label{eq:start_adpt0} 
\end{align}
where  $\widetilde{\Delta}^{W}_0$ is defined in \eqref{eq:pi} and \eqref{eq:tildedelta}.
Second, we particularise \eqref{eq:start_adpt0} to the different smoothness cases and prove (T\ref{cor:adp_rate}.\ref{a0}), (T\ref{cor:adp_rate}.\ref{a01}), (T\ref{cor:adp_rate}.\ref{b0}), and (T\ref{cor:adp_rate}.\ref{a02}).\\
\noindent {\bf Part 1.}  By Lemma \ref{adp2} and Lemma \ref{adp1}, 
we have 
\begin{align*}
&\mathcal{R}_{n_0}^{W}\left(\widehat{f}_{\alpha,\mt{\beta}}^{q,\widehat{N},\widehat{T},\epsilon},f_{\alpha,\mt{\beta}}\right)\\
&\leq \frac{C \left(2+\sqrt{5}\right)^4}{2\pi} \left(\int_{\epsilon\leq \abs{t} } \E\left[\mathcal{L}_{q}^{W}\left(t,N(t),T\right)\right] dt +  \frac{1+c_1}{2+\sqrt{5}}\E\left[\Sigma_2(T, N)\right] \right) \notag \\
& \quad + \frac{2(2+\sqrt{5})^2 C\left(1+\left(2+\sqrt{5}\right)^2\right)}{\pi}  \Pi(n, Z_{n_0},T_{\max},N_{\max,q})   + CM^2\widetilde{w}(\underline{a}). \notag 
\end{align*}
The definition of $\Sigma$,
\eqref{eq:R_0_w}, \eqref{eq:R2_temp}, \eqref{eq:R_2_w}, \eqref{eq:R3_w}, 
and Lemma \ref{lem:Z_control} 
yield \eqref{eq:start_adpt0}. 

\noindent {\bf Part 2.} 
\noindent We 
start from \eqref{eq:start_adpt0} 
and use $A_4 := 2(2+\sqrt{5})^2 C\left(1+\left(2+\sqrt{5}\right)^2\right)(A_0+A_1)/\pi$,  $T^*: = 2^{k^*}$,  $k^*: = \lfloor \ln(\phi^I(\omega_{\underline{N}^*}))/\ln(2) \rfloor$, $N^*(t) := \lfloor \underline{N}^*\rfloor$, where $\underline{N}^*$ is defined below, and 
$$\mathcal{R}^{W,\rm{adp}}_{n_0,\sup}:=\sup_{f_{\alpha,\mt{\beta}} \in \mathcal{H}^{q,\phi,\omega}_{w,W}(l, M)\cap\mathcal{D},\ f_{\mt{X}|\mathcal{X}} \in \mathcal{E}} \mathcal{R}_{n_0}^{W}\left(\widehat{f}_{\alpha,\mt{\beta}}^{q,\widehat{N},\widehat{T},\epsilon},f_{\alpha,\mt{\beta}}\right).$$
%
We have, for all $|t|\geq \epsilon$, $W\in\{W_{[-R,R]},\cosh(\cdot/R)\}$, and $N\geq 1$, $2\widetilde{\Delta}^{W}_0(t,N,n,Z_{n_0})/(\pi(2\pi)^p) \leq \widetilde{\Delta}^{W}(t,N,n,Z_{n_0})$ where  $\widetilde{\Delta}^{W}$ is defined like $\Delta^{W}$ replacing $Q_q$  by $Q_{q,n}:= Q_q\left(1 +2(1+2p_n)(1+c_1)  \right)$. 
Thus, by \eqref{eq:start_adpt0}, we obtain, for all $W\in\{W_{[-R,R]},\cosh(\cdot/R)\}$,
\begin{align}
\mathcal{R}^{W,\rm{adp}}_{n_0,\sup} \notag  \leq &   C \left(2+\sqrt{5}\right)^4\int_{\epsilon\leq \abs{t}\leq T }  \widetilde{\Delta}^{W}\left(t,N(t),n, \frac{M_{\mathcal{E},\eta} v(n_0,\mathcal{E})}{\delta(n_0)} \right)  dt \notag \\ 
&+ C\left(2+\sqrt{5}\right)^4 \left(\sup_{t\in [-T,T]} \frac{ 4 l^2}{ \omega_{N(t)+1}^2} +\frac{ 4 l^2}{\phi(T)^2} \right)   + C M^2\widetilde{w}(\underline{a})+ \frac{A_4}{n}.\label{eq_end1}\end{align}
\noindent {\bf Proof of (T\ref{cor:adp_rate}.\ref{a0})}.  Let $\underline{N}^*$ solution of 
\begin{equation}\label{eq:nstar}
2k_q\underline{N}^*\ln\left(\underline{N}^*\ln(n)\right) + \ln\left(\omega_{\underline{N}^*}^2\right)+(p-1)\ln\left(\underline{N}^*\right)+ \ln_2(n)= \ln(n),
\end{equation} 
$n\geq e^{7e^2/(2\pi)}$ large enough so $N^*\geq 1$, and $(\ln(n) /\tau_2')^{\sigma/s} \leq  n^{\zeta_0}/2$, where  $\tau_2':= 2k_q\mathcal{W}(e/(2k_q))$. By \eqref{eq:nstar} and the definition of $\underline{N_{\max,q}^{W}}$,  we have $ N^* \leq N_{\max,q}^{W}$ for all $t\in \R\setminus(-\epsilon,\epsilon)$, hence $N^*\in\mathcal{N}_n$. 
Also $T^*\in  \mathcal{T}_n$ because, by the arguments in the proof of (T\ref{theo:compact}.\ref{t_comp_1}), 
\begin{equation*}
\underline{N}^* \leq \frac{\ln(n)}{ 2k_q\mathcal{W}(\ln(n)/(2k_q))} \leq \frac{\ln(n)}{\tau_2'},
\end{equation*}
so  we have $ T^* \leq  (\ln(n) /\tau_2')^{\sigma/s} \leq  n^{\zeta_0}/2\leq T_{\max}$. \eqref{eq:int_t0}, \eqref{eq:int_t1}, and 
$p_n=6\left(1+\zeta_0\right)\ln(n)$, yield
\begin{align*}
\mathcal{R}^{W,\rm{adp}}_{n_0,\sup} &\leq   C\left(2+\sqrt{5}\right)^4 \ln(n)(N^*)^{p-1} \left( \frac{\tau_0'}{n}\left(N^*\ln(n)\right)^{2k_qN^*} + \frac{\tau_1'}{n} (T^*)^{p+1} \right) +\frac{ 8C\left(2+\sqrt{5}\right)^4 l^2}{ \omega_{N^*}^{2} } \\
&  \quad  + \frac{\theta C M^2}{\ln(n)} + \frac{A_4}{n},\\ 
\tau_0'&:=\frac{e^{2k_q}4c_{\mt{X}}\theta^{p+1}}{\pi^{p+1}k_q} \left( Q_q \left(\frac{1}{e}+ \left(\frac{1}{e}+6\left(1+\zeta_0\right)\right)2(1+c_1)\right)+eM_{\mathcal{E},\eta}2^p\right)+  \frac{M_{\mathcal{E},\eta}2^{p+2}l^2}{e},  \notag\\ 
\tau_1' &:=\frac{4c_{\mt{X}}}{\pi^{p+1}(p+1)}\left(Q_q\left(\frac{1}{e}+ \left(\frac{1}{e}+6\left(1+\zeta_0\right)\right)2(1+c_1)\right) +\frac{M_{\mathcal{E},\eta}2^p}{e}\right). \notag
\end{align*} 
The computation below gives lower bounds on $N^*\ln(N^*)$ and $N^*$: 
\begin{align*}
\ln(n)&=2k_qN^*\ln\left(N^*\ln(n)\right) + (2\sigma + p-1)\ln\left(N^*\right) + \ln_2(n)\\
&\leq 2\left( (2(k_q+\sigma)+p-1) N^*\ln(N^*) \bigvee  (2k_qN^*+1) \ln_2(n)\right)\\
&\leq 2\left( (2(k_q+\sigma)+p-1) N^*\ln(N^*) \bigvee  (2k_q+1)N^* \ln_2(n)\right).
\end{align*}
Using both and $\mathcal{W}(x) \leq  \ln(x +1 )$ for all $x >0$, we obtain $ N^* \geq \ln(n)/((\tau_8\bigvee (2k_q+1)) \ln_2(n) (1 +\ln(1+\tau_8/e))  $, where $\tau_8 :=2(2(k_q+\sigma) + p-1)$.
We conclude proceeding like for \eqref{eq_compact1}.\\
\noindent {\bf Proof of (T\ref{cor:adp_rate}.\ref{a01}).} 
Start from \eqref{eq_end1}, where,  because $ w= W_{\mathcal{A}}$, the term $ M^2\widetilde{w}(\underline{a}) $ is zero.
Let $\underline{N}^*$ solution of $2k_q\underline{N}^*\ln\left(\underline{N}^*\right) + \ln\left(\omega_{\underline{N}^*}^2\right)+(p-1)\ln\left(\underline{N}^*\right)+ \ln_2(n)= \ln(n)$. By definition of $\underline{N_{\max,q}^{W}}$ this yields $ N^*\leq N_{\max,q}^{W} $ hence $N^*\in \mathcal{N}_n$. Using arguments from the proof of (T\ref{theo:compact}.\ref{poly00}) we have $ T^* \leq n^{\kappa/(2(\kappa+k_q)s)}$ and, using $  s> 2p+1/2$, for $n$ large enough  $n^{\kappa/(2(\kappa+k_q)s)}\leq n^{1/(4p+1)} /2 \leq T_{\max} $, hence $T^*\in \mathcal{T}_n$. Thus, we obtain   
$$
\mathcal{R}^{W,\rm{adp}}_{n_0,\sup} 
\leq    C\left(2+\sqrt{5}\right)^4\ln(n)(N^*)^{p-1}\left( \frac{\tau_0'}{n}\left(N^*\right)^{2k_qN^*} + \frac{\tau_1'}{n} (T^*)^{p+1} \right) + \frac{ 8 C\left(2+\sqrt{5}\right)^4l^2}{ \omega_{N^*}^{2} } + \frac{A_4}{n}.
$$
This yields the result following the proof of (T\ref{theo:compact}.\ref{poly00}).\\ 
\noindent {\bf Proof of (T\ref{cor:adp_rate}.\ref{b0})}. 
Starting from \eqref{eq_end1}, the proof is similar to that of (T\ref{cor:adp_rate}.\ref{a0}) with elements from that of (T\ref{theo:non_compact}.\ref{t_noncomp_1}), using $\underline{N}^*$ solution of $ 2k_q\underline{N}\ln_2\left(n\right) + (p-1)\ln(\underline{N})/q + \ln\left(\omega_{\underline{N}^*}^2\right) = \ln(n)$.\\
\noindent {\bf Proof of (T\ref{cor:adp_rate}.\ref{a02}).} The proof is similar to that of (T\ref{cor:adp_rate}.\ref{a01}). Start from  \eqref{eq_end1}. Let $\underline{N}^*$ solution of $  2(k_q+\kappa)\underline{N}^* + (p-1)\ln(\underline{N}^*)/q + \ln_2(n)=\ln(n)$. Then, using 
the definition of $\underline{N_{\max,q}^{W}}$, which satisfies $ 2k_q\underline{N_{\max,q}^{W}} = \ln(n)$, we have $ N^* \in \mathcal{N}_{n}$. Using arguments from the proof of (T\ref{theo:non_compact}.\ref{t_noncomp_2}), we have $ T^* \leq n_e^{\kappa/(2s(\kappa + k_q))}$ and, using $  s>4p+1/2$, for $n$ large enough  $ n_e^{\kappa/(2s(\kappa + k_q))}\leq n^{1/(8p+1)}/2 \leq T_{\max}  $,  hence $T^*\in \mathcal{T}_n$. This yields the result using the proof of (T\ref{theo:non_compact}.\ref{t_noncomp_2}). \hfill $\square$ 

\renewcommand{\theequation}{B.\arabic{equation}}
\renewcommand{\thelemma}{B.\arabic{lemma}}
\renewcommand{\thecorollary}{B.\arabic{corollary}}
\renewcommand{\thedefinition}{B.\arabic{definition}}
\renewcommand{\theproposition}{B.\arabic{proposition}}
\renewcommand{\theremark}{B.\arabic{remark}}
\renewcommand{\thetheorem}{B.\arabic{theorem}}
\renewcommand{\theassumption}{B.\arabic{assumptio}}
\renewcommand{\thesubsection}{B.\arabic{subsection}}
\setcounter{equation}{0}  
\setcounter{lemma}{0}
\setcounter{corollary}{0}
\setcounter{proposition}{0}
\setcounter{remark}{0}
\setcounter{definition}{0}
\setcounter{lemma}{0}
\setcounter{theorem}{0}
\setcounter{assumption}{0}
\setcounter{subsection}{0}
\setcounter{footnote}{0}
\setcounter{figure}{0}

\section{Harmonic analysis}\label{s3}
\subsection{Preliminaries}
$P_m$ is the Legendre polynomial of degree $m$ with $\|P_m\|_{L^2([-1,1])}=1$. 

\begin{lemma}\label{lem:interp_error}
	For all $ f\in L^2_{w}(\R)$, $w$ even, nondecreasing on $[0,\infty)$, and $w(0),R>0$,  we have
	$     \bigl\| \mathcal{P}_{R}\left[ \mathcal{F}\left[f\right]\right] - \mathcal{F}\left[f\right]  \bigr\|^2_{L^2(\R)} \leq (2\pi/w(R)) \norm{f}_{L^2(w)}^2.$
	
\end{lemma}
\noindent {\bf Proof.}
The result uses the Plancherel identity and 
\begin{align*}
\bigl\|\mathcal{P}_{R}\left[ \mathcal{F}\left[f\right]\right] - \mathcal{F}\left[f\right]\bigr\|^2_{L^2(\R)} 
& = 2\pi \int_{\R} \indic\{|a|>R\} \left|f\left(  a\right)\right|^2da      \leq  \frac{2\pi}{w(R)} \int_{\R}  \abs{f\left(a\right)}^{2}  w(a)da.\quad\square
\end{align*}

\begin{proposition}\label{pdebut}
	For all weighting function $W$, $c\in\R$, $R>0$, and $m\in\N_0$, we have
	\begin{enumerate}[\textup{(}i\textup{)}]
		\item\label{pdebuti} $g_m^{W(\cdot/R),c}=g_m^{W,Rc}$ in $L^2([-1,1])$,
		\item\label{pdebutii} $\sigma_m^{W(\cdot/R),c}=\sigma_m^{W,Rc}\sqrt{R}$,
		\item\label{pdebutiii} $\varphi_m^{W(\cdot/R),c}=\varphi_m^{W,Rc}\left(\star/R\right)/\sqrt{R}\ a.e.$
	\end{enumerate} 
	
\end{proposition}
\noindent {\bf Proof.} \eqref{pdebuti} follows from 
$\mathcal{Q}_c^{W(\cdot/R)}=\mathcal{Q}_{Rc}^{W}$ and 
\eqref{pdebutii} from 
$\sigma_m^{W(\cdot/R),c}=\sqrt{2\pi\rho_{{m}}^{W(\cdot/R),c}/\abs{c}}=\sqrt{2\pi\rho_{{m}}^{W,Rc}/\abs{c}}
$ (by the argument yielding \eqref{pdebuti}).  
Now, using \eqref{pdebuti} in the first display and \eqref{pdebutii} in the last display, we have, for a.e. $t\in\R$, 
\begin{align*}
\sigma_m^{W,Rc}\varphi_m^{W,Rc}\left(\frac{t}{R}\right)
&=\mathcal{F}_{Rc}^*\left[g_m^{W(\cdot/R),c}\right]\left(\frac{t}{R}\right)\quad(\text{where} \ \mathcal{F}_{Rc}^*: \ L^{2}([-1,1]) \to L^{2}(W))\\
&=\mathcal{F}_{c}^*\left[g_m^{W(\cdot/R),c}\right](t) \quad \left(\text{where} \ \mathcal{F}_{c}^*: \ L^{2}([-1,1]) \to L^{2}(W(\cdot/R))\right)   \\
&=\sigma_m^{W(\cdot/R),c}\varphi_m^{W(\cdot/R),c}(t)
=\sigma_m^{W,Rc}\sqrt{R}\varphi_m^{W(\cdot/R),c}(t),
\end{align*}
hence \eqref{pdebutiii} when we divide by $\sigma_m^{W,Rc}$ which is nonzero.
\hfill$\square$

\begin{proposition}\label{eq:g_m_control}
	For all $ \mt{m}\in\N_0^p $, $R>0$, $W= W_{[-R,R]}$ or   $W=\cosh(\cdot/R)$, $t\neq 0$, we have 
	$\left\|g^{W,x_0t}_{\mt{m}}\right\|_{L^\infty([-1,1]^p)} \leq H_W(t) \prod_{j=1}^p  \sqrt{m_j+1/2}$, where 
	$H_{W_{[-R,R]}}(t)=H_0^p\left(1+(x_0|t|)^2\right)^p$, $H_0=2(1+1/\sqrt{3})$, 
	$H_{\cosh(\cdot/R)}(t)= H_1^p  (1 \vee (x_0|t|)^{4})^p$, $H_1 >0$.
\end{proposition}
\noindent {\bf Proof.} 
When $ W=W_{[-R,R]} $, this is (66) in \cite{bonami2016uniform} else 
this is (50) in \cite{Note}.\hfill$\square$

\begin{lemma}\label{lem:lowerbound_weight10}
	For all $q\in \{1,\infty\}$, $t\ne0$,  $R >0$, $N\in\N_0$,  in cases (N.1) and (N.2) of Section \ref{s4},  we have
	$\sum_{\abs{\mt{m}}_{q} \leq N} 1/\rho_{\mt{m}}^{W,t} \leq  \nu_q^{W}(t,N)$. 
\end{lemma}
\noindent {\bf Proof.} Let $R >0$. We use repeatedly, for all $x>0$ and $N\in \N_0$,  
\begin{align}
\sum_{k\leq N} \exp\left(k x \right)  &\le \frac{\exp\left((N+1/2)x\right)}{2\sinh\left(x/2\right)} \leq  
\frac{\exp\left((N+1/2)x\right)}{x} \quad(\text{because }\sinh(\abs{x}) \geq \abs{x}),\label{eq:up1}\\
&\le \frac{\exp\left(Nx\right)}{1-\exp(-x)},\label{eq:up2}
\end{align}
the cardinal of $\{\mt{m}\in\N_0^p:\ \abs{\mt{m}}_1 =k\}$ is $
{k+p-1\choose p-1}$, 
and $(k+p-1)!/k! \leq (k+p-1)^{p-1}$, and for all $m\in\N_0$, $\rho_{m}^{\cosh, Rt}=\rho_{m}^{\cosh(\cdot/R), t}$ and 
$\rho_{m}^{W_{[-1,1]}, Rt} = \rho_{m}^{W_{[-R,R]},t}$.\\
Start by case (N.2). Let $|t| > \pi/4$ and $q=1$.   By (11) in \cite{Note} (there $\mathcal{Q}_{t}$ differs by a factor $1/(2\pi)$), we have, for all $ m \in \N_0  $, 
\begin{align}
\rho_m^{\cosh,t} & \geq \frac{1}{2}\exp\left(-\frac{\pi(m+1)}{2\abs{t}}\right). \label{eq:rho}
\end{align}
The result is obtained from the above with \eqref{eq:up1} and
\begin{align}
\sum_{\abs{\mt{m}}_1\leq N} \frac{1}{ \rho_{\mt{m}}^{\cosh,t}} \leq &  2^{p} \sum_{k\leq N} \sum_{\abs{\mt{m}}_1 =k} \exp\left( \frac{\pi (\abs{\mt{m}}_1+p)}{2\abs{t}} \right)\notag\\
\leq  &  \frac{2^{p+1}(N+p-1)^{p-1}e\abs{t}}{\pi (p-1)!}   \exp\left( \frac{\pi (N+p)}{2\abs{t}} \right)\label{efin}.
\end{align}
Let $|t| \leq  \pi/4$ and $q=1$. By Theorem 1 in \cite{Note}, we have, for all $m \in \N_0$, 
\begin{align}
\rho_m^{\cosh,t} & \geq   \left(\frac{2}{e\pi}\right)^2 \exp\left( - 2\ln\left(\frac{7e^2\pi}{2|t|}\right)m \right). \label{eq:rho1}
\end{align}
The result is obtained from the above with \eqref{eq:up2} and
\begin{align}
\sum_{\abs{\mt{m}}_1\leq N} \frac{1}{ \rho_{\mt{m}}^{\cosh,t}} 
&\leq   \left(\frac{e\pi}{2} \right)^{2p} \sum_{k\leq N} \sum_{\abs{\mt{m}}_1 =k}  \exp\left(  2 \ln\left(\frac{7e^2\pi}{2|t|}\right)\abs{\mt{m}}_1 \right)  \notag\\
&\leq  \left(\frac{e\pi}{2}   \right)^{2p}  \frac{(N+p-1)^{p-1}}{(p-1)!}\exp\left( 2	 \ln\left(\frac{7e^2\pi}{2|t|}\right) N \right) \frac{1}{1-\left(1/(14e^2\pi)\right)^2} \label{efinb}.
\end{align}
The results for $q=\infty$ are obtained using \eqref{efin}  and \eqref{efinb} with $p=1$ and  
\begin{align}\label{eq:dimp}
\sum_{\abs{\mt{m}}_{\infty} \leq N}\frac{1}{ \rho_{\mt{m}}^{\cosh,t}}  & \leq \prod_{j=1}^p \left( \sum_{\mt{m}_j=0}^{N}  \frac{1}{\rho_{\mt{m}_j}^{\cosh,t}} \right).
\end{align}
Consider case (N.1). Let $t\neq0$. Because $14e \geq 1$ and by Lemma \ref{lem:lowerbound}, we have, for all $m\in\N_0$,
\begin{align}\label{eq:lower1}
\rho_m^{W_{[-1,1]},t} \geq \frac{1}{2}\left( \frac{\abs{t}}{7e\pi(m+1)} \bigwedge 1   \right)^{2m} .
\end{align}
When $q=1$, the result follows from the following sequence of inequalities
\begin{align*}
\sum_{\abs{\mt{m}}_1\leq N} \frac{1}{ \rho_{\mt{m}}^{W_{[-1,1]},t}} \leq & 2^p \sum_{k\leq N} \sum_{\abs{\mt{m}}_1 =k} \prod_{j=1}^p\exp\left( 2\mt{m}_j \ln\left(\frac{7e\pi(\mt{m}_j+1)}{\abs{t}}\bigvee 1\right)   \right)\\
\leq & \frac{2^{p}(N+p-1)^{p-1}(N+1) }{(p-1)!} \exp\left( 2 N \ln\left(\frac{7e\pi(N+1)}{\abs{t}} \bigvee 1\right)   \right).  
\end{align*}
When $q=\infty$, we obtain the result using the above with $p=1$ and \eqref{eq:dimp}.\hfill$\square$\vspace{0.3cm}

The proof of the next lemma is straightforward. 
\begin{lemma}\label{lem:bc}
	Let $f_{\alpha,\mt{\beta}}\in L^2\left(w\otimes W^{\otimes p}\right)$. 
	For all $ \mt{m}\in \N_0^p,\ t\neq0$, we have $c_{\mt{m}}=\sigma_{\mt{m}}^{W,x_0t} b_{\mt{m}}\ a.e.$
\end{lemma}


\subsection{Properties of the PSWF and eigenvalues}\label{app:lower1} 
\begin{lemma}\label{upper_bound}
	For all $ c \ne0$ and $ m\in\N_0$, we have $\left|\mu_{m}^c\right| \leq
	\sqrt{2\pi} e^{3/2} \left( e\abs{c}/\left(4 (m+3/2)\right)\right)^m/3$. 
\end{lemma}
\noindent {\bf Proof.} Let $c\ne0$ and $m\in\N_0$. By (69) in \cite{Xiao}, 6.1.18 in \cite{abramowitz1965handbook}, 
(7) in \cite{gamma}, (1.3) in \cite{inq_gamma2}, 
we obtain
\begin{align*} 
\left|\mu_{m}^c\right|
&\leq \frac{\sqrt{\pi}\abs{c}^m(m!)^2}{(2m)!\Gamma(m+3/2)}\\ 
&\leq \frac{\pi \abs{c}^m}{4^m\Gamma\left(m+3/2\right)} \frac{\Gamma(m +1)}{\Gamma\left(m + 1/2\right)}\\
&\leq \frac{\pi \abs{c}^m}{4^m\Gamma(m+3/2) } (m+1)^{1/2}\leq \frac{\sqrt{\pi e^3}(e\abs{c})^m(m+1)^{1/2}}{4^m\sqrt{2}(m+3/2)^{m+1}}
\end{align*}  
and conclude using $\sup_{x\ge0} (x+1)^{1/2}/(x+3/2)\le2/3$.\hfill$\square$

\begin{lemma}\label{lem:lowerbound}
	For all $c \ne0$ and $m\in\N_0$, we have
	\begin{equation*}
	\rho_{m}^{W_{[-1,1]},c} \geq \frac{1}{2}\left(\indic\left\{m  \leq \frac{2\abs{c}}{\pi} -1\right\}+\frac{c}{7e\pi(m+1)}\indic\left\{m>\frac{2\abs{c}}{\pi} -1\right\}\right)^{2m}. 
	\end{equation*} 
\end{lemma}
\noindent{\bf Proof.} When $m\ge 2\abs{c}/\pi -1$, the result follows from the fact that, by Proposition 5.1 in \cite{bonami2018} and the Tur\'an-Nazarov inequality (see \cite{nazarov2000complete} page 240),
$ \rho_{m}^{W_{[-1,1]},c} \geq \left(c/\left(7e\pi (m+1)\right)\right)^{2m}/2.$
For all $m  \leq   2\abs{c}/\pi -1$, the result follows from Remark 5.2 in \cite{bonami2018}  and that, for all $m\in\N_0$,  $ c\in (0,\infty) \mapsto \rho_{m}^{c} $ is nondecreasing (by the arguments in the proof of  Lemma 1 in \cite{Note}).\hfill $\square$\vspace{.3cm}

We now use $\Pi(c) :=3c^2\exp\left(2 c^2/\sqrt{3}\right)/16$, $H(c) := \sqrt{2\Pi(c)}\vee 2$, $l(c): =  \left(1+4c^2/3^{3/2}\right)\left(1+2c^23^{3/2}\right)$, if $N\ge H(c)$ then $N\geq c$ (for all $c\geq 2$, $N\geq c \sqrt{3 \exp(	8/\sqrt{3})/16}> c$ else $N\ge H(c)\ge 2>c$), $f(x):= | x |/(1-x^2)$, $g(x):=\abs{x} /(1-x)^2$, $h(x):= | x |/(1-\abs{x})$, $c_f:=4/3$, $c_g := 4$, $c_h := 2$,
\begin{align}
&\forall x\in[-1/2,1/2],\ f(x) \leq c_f\abs{x},  \  g(x) \leq c_g\abs{x},\  h(x) \leq c_h\abs{x};\label{ubfg}\\ 
&2\sum_{k \equiv N[2],\ 0<k < N} 2k+1 =N(N-1).\label{sumod}
\end{align}
\eqref{sumod} is obtained because for all $N$  even the sum is 
$2\sum_{p=1}^{N/2-1}4p +1$ and else $2\sum_{p=0}^{(N-1)/2-1}4p+3$. 
\begin{lemma}\label{lem:ratio}
	For all $c\ne0$ and $ m \geq 2$, 
	we have
	$\abs{\mu_m^c/\mu_{m-2}^c}\leq \Pi(c)/m^2$.
\end{lemma}
\noindent {\bf Proof.}  
Let $c>0$ and $m\in\N_0$ (for $ c<0 $, we use $\mu_m^c=\overline{\mu_m^{-c}}$). By Theorem 8.1 in \cite{Osipov}, we have 
$$
\abs{\mu_m^c} = \frac{ \sqrt{\pi} c^m (m!)^2  }{(2m)! \Gamma(m + 3/2)} e^{F_m(c)},\ 
F_m(c)= \int_{0}^c \left( \frac{2  \left(\psi_m^t(1)\right)^2 -1}{2t} -\frac{m}{t} \right) dt. 
$$
Moreover, by (65) in \cite{bonami2016uniform}, for all $t>0$,
\begin{equation*}
\left(\sqrt{m+\frac{1}{2}} - \frac{t^2}{\sqrt{3}\sqrt{m+1/2}}\right)^2  \leq \left(\psi_m^t(1)\right)^2 \leq \left(\sqrt{m+\frac{1}{2}} + \frac{t^2}{\sqrt{3}\sqrt{m+1/2}}\right)^2
\end{equation*} 
which yields, if $m \geq 2$,
\begin{align}
\left(\psi_m^t(1)\right)^2 - \left(\psi_{m-2}^t(1)\right)^2 
&\leq   \left(\sqrt{m+\frac{1}{2}} + \frac{t^2}{\sqrt{3}\sqrt{m+1/2}}\right)^2 - 
\left(\sqrt{(m -2)+\frac{1}{2}} - \frac{t^2}{\sqrt{3}\sqrt{(m -2)+1/2}}\right)^2\notag \\
&= 2 + \frac{4t^2}{\sqrt{3}} + \frac{t^4}{3}\left(\frac{1}{m+1/2} - \frac{1}{m-3/2} \right)
\leq  2 + \frac{4t^2}{\sqrt{3}} .\label{eq:FCdiff}
\end{align}
Using $\sup_{x\ge 2}x^3(x-1)/\left((x^2 - 1/4)(x - 1/2)(x-3/2)\right)\le3$ and \eqref{eq:FCdiff}, for all $m \geq 2$,
\begin{align*}
\abs{ \frac{\mu_m^c}{\mu_{m-2}^c}} & = \frac{c^2}{16}\frac{m(m-1)}{(m^2 - 1/4)(m - 1/2)(m-3/2)} \exp\left(F_m(c)- F_{m-2}(c) \right) \\
& \leq \frac{3c^2}{16m^2}\exp\left( \int_{0}^c \left( \frac{ \left(\psi_m^t(1) \right)^2 - \left(\psi_{m-2}^t(1) \right)^2 }{t} -\frac{2}{t} \right) dt  \right) 
\leq \frac{3c^2}{16m^2} \exp\left(  \frac{2c^2}{\sqrt{3}}  \right).\quad\square
\end{align*}

\begin{lemma}\label{rem:psi}
	For all $ c\neq 0 $ and $k\in \N$, we have $\left(\psi_k^c(1)\right)^2\leq  \left(k+1/2\right)\left(1+2c^2/3^{3/2} \right)^2$ and  $\|\psi_k^c\|_{L^{\infty}([-1,1])}^2\le\left(k+1/2\right)\left(1+4c^2/3^{3/2} \right)^2$. For all $ c\neq 0 $ and $k\ge c$, we have $\|\psi_k^c\|_{L^{\infty}([-1,1])}^2\le k+1/2$. We also have $\|\psi_0^c\|_{L^{\infty}([-1,1])}^2\le 2|c|/\pi$.
\end{lemma}
\noindent {\bf Proof.} 
The first assertion follows from (65) in \cite{bonami2016uniform}. For the second, we use (66) in \cite{bonami2016uniform} in the first display,  22.14.7 and 22.2.10 in \cite{abramowitz1965handbook}, hence $\|P_k\|_{L^{\infty}([-1,1])} \leq \sqrt{k+1/2}$, in the second inequality, 
\begin{align*}
\|\psi_k^c\|_{L^{\infty}([-1,1])}&\leq \|P_k\|_{L^{\infty}([-1,1])} +  \frac{c^2}{\sqrt{3(k+1/2)}}\left(1+\frac{\sqrt{3/2}}{\sqrt{k+1/2}}\right)\\
&\leq  \sqrt{k+1/2} \left(1+\frac{c^2}{\sqrt{3}(k+1/2) }\left(1+\frac{\sqrt{3/2}}{\sqrt{k+1/2}}\right)\right)
\leq  \sqrt{k+1/2}\left(1+\frac{4c^2}{3^{3/2}} \right). 
\end{align*}
The third uses (3.4) and (3.125) in \cite{Osipov}. We obtain the last by the proof of Proposition 1 in \cite{karoui2008new} which yields $\|\psi_0^c\|_{L^{\infty}([-1,1])}^2\le 2/(\mu_0^c)^2$ and Lemma \ref{lem:lowerbound}. For all $ c<0 $, we use 
$\psi_m^{-c}=\psi_m^{c}$.\hfill $ \square  $

\begin{lemma}\label{lem:first}
	For all $c\ne0$ and $ N\ge H(c)$, 
	we have
	\begin{align*}
	&\left\|\frac{\partial  \psi_N^{c}}{\partial c}\right\|_{L^{\infty}([-1,1])} \leq \frac{2c_f\left( C_1(c)+ C_2(c) \right)C_3(c)\Pi(c)}{\abs{c}}\sqrt{N},\\
	&C_1(c):=\frac{2H(c)+9}{(H(c)+2)^2}
	,\  C_2(c):= \frac{2\abs{c}}{\pi H(c)(H(c)-1)} + \frac{l(c)}{4},\  C_3(c):= \sqrt{1+\frac{1}{2H(c)}}.
	\end{align*}
\end{lemma}
\noindent {\bf Proof.} Take $c\ne0$, $N\ge H(c)$, and $w\in[-1,1]$. Theorem 7.11 in \cite{Osipov} yields
\begin{equation}\label{eq:first_decomp}
\frac{\partial  \psi_N^{c}}{\partial c}(w) =\frac{2\psi_N^{c}(1)}{\abs{c}}  \sum_{ k \equiv N[2],\ k \neq N}\frac{\mu_N^{c}\mu_k^{c} }{\left(\mu_N^{c}\right)^2 - \left(\mu_k^{c}\right)^2}  \psi_k^{c}(1)\psi_k^{c}(w). 
\end{equation}
\noindent Using $\mu_k^c/\mu_N^c\in\R$ if $k\equiv N[2]$ and Lemma \ref{rem:psi}, we obtain
\begin{align*}
&\abs{\frac{\partial  \psi_N^{c}}{\partial c}(w )} \leq \frac{\sqrt{4N+2} }{\abs{c}}\mathcal{C}(f,N,c),\\
&\mathcal{C}(f,N,c):= f\left( \frac{\mu_N^c}{\mu_0^c}\right)\frac{2|c|\indic\{N\equiv 0 [2]\}}{\pi}+\sum_{\underset{ k \equiv N[2]}{0<k < N} } f\left( \frac{\mu_N^c}{\mu_k^c}\right)l(c)\left(k + \frac{1}{2}\right)
+\sum_{\underset{ k \equiv N[2]}{k > N} } f\left( \frac{\mu_k^c}{\mu_N^c}\right)\left(k + \frac{1}{2}\right).\notag
\end{align*}
Lemma \ref{lem:ratio} yields, if $k\equiv N[2]$,
\begin{equation}\label{b1}
\left|\frac{\mu_N^c}{\mu_k^c}\right|\le\left|\frac{\mu_N^c}{\mu_{N-2}^c}\right|\le \frac{\Pi(c)}{N^2}\le \frac{1}{2}\ \textrm{if }k<N\  \textrm{and }
\left|\frac{\mu_k^c}{\mu_N^c}\right|\le \left(\frac{\sqrt{\Pi(c)}}{N+2}\right)^{k-N}\le\frac{1}{2}\ \textrm{if } k>N.
\end{equation} 
Using \eqref{sumod}, \eqref{ubfg}, 
\eqref{b1}, and $\sum_{k\in\N}k2^{-k}  = 2$ in the third display,  
the result follows from
\begin{align}
\mathcal{C}(f,N,c)
&\leq c_f \left(\left(\frac{2\abs{c}}{\pi} + \frac{l(c) N(N-1)}{4} \right) \frac{\Pi(c)}{N^2} +\sum_{k \equiv N[2],\ k > N}\frac{k+1/2}{ 2^{(k-N)/2}}\left(\frac{\sqrt{2\Pi(c)}}{N+2}\right)^{k-N}\right) \notag\\
&\leq c_f\Pi(c)\left(\frac{2\abs{c}}{\pi H(c)(H(c)-1)} + \frac{l(c)}{4}     
+\frac{2}{(N+2)^2}\sum_{l \equiv 0[2],\ l \ge2 }\left(l+N+\frac12\right)\frac{1}{ 2^{l/2}}\right)\notag \\
&\leq  c_f\Pi(c)\left(C_2(c)
+\frac{2}{(N+2)^2}\left(N+\frac{9}{2}\right)\right)\leq c_f\Pi(c)\left(C_1(c)+C_2(c)  \right). \quad\quad\quad\quad\quad\quad\square\label{eq:lem4}
\end{align}
\begin{lemma}\label{lem:second}
	For all $ c\ne0 $ and $N\geq H(c)$, 
	we have
	\begin{align*}
	&\left\|\frac{\partial^2  \psi_N^{c}}{\partial c^2 }\right\|_{L^{\infty}([-1,1])}  \leq \frac{\Pi(c)C_3( c)}{c^2}\left(C_4(c) N^{5/2} +  C_5(c) N^{3/2}   +  C_6(c) \sqrt{N} +  C_7(c)\right),\\
	&C_4(c): = c_g\left(C_2(c) - C_1(c)\right),\ C_7(c): =
	\frac{c_g}{(H(c)+2)^{1/2}} 
	\left( 85 + \frac{246}{H(c)+2}\right),\\
	&C_5(c):=   8\left(c_f\left(C_1( c) +C_2( c) \right)C_3( c)\right)^2  \Pi(c)+\left(c_g +4 c_f\right)C_2( c)+\left(8 c_f-c_g  \right)C_1( c)+2c_g,\\
	&C_6(c) :=8c_h c_f ( C_1( c)+C_2( c))^2\Pi(c)+  (C_1(c)+C_2( c))\left( c^2c_g + 4 c_f\right)+19c_g.
	\end{align*}
\end{lemma}
\noindent {\bf Proof.} For all $ c<0 $, 
$\mu_m^c=\overline{\mu_m^{-c}}$ and $\psi_m^{-c}=\psi_m^{c}$, hence we only consider $c>0$. Using 
$c\in (0,\infty)\mapsto \psi_{N}^c(x)$ is analytic (see \cite{Fuchs} page 320) and (7.99) in \cite{Osipov}, we have by differentiating 
\begin{align}
\mu_N^c \psi_N^c(x)&= \int_{-1}^1 e^{icxt}\psi_N^c(t) dt:\label{eq:diff1}\\
\mu_N^c \frac{\partial\psi_N^c}{\partial x}(x) &= \int_{-1}^{1} ict e^{ic x t}  \psi_N^c(t) dt\label{eq:first_psi2},\\ 
\mu_N^c \frac{\partial^2\psi_N^c}{\partial x^2}(x) &= -\int_{-1}^{1} (ct)^2 e^{ic x t}  \psi_N^c(t) dt\label{eq:first_psi3},\\
\left(\frac{\partial^2 \mu_N^c}{\partial c^2} \psi_N^c + 2 \frac{\partial \mu_N^c}{\partial c} \frac{\partial \psi_N^c}{\partial c} + \mu_N^c \frac{\partial^2 \psi_N^c}{\partial c^2}\right)(x) 
&= \int_{-1}^1  e^{icxt}\left(\frac{\partial^2 \psi_N^{c}}{\partial c^2}(t) + 2 ixt\frac{\partial \psi_N^c}{\partial c}(t) -(xt)^2 \psi_N^c(t)\right)  dt\label{eq:first_psi}.
\end{align}
Multiplying 
\eqref{eq:first_psi}
by $ \psi_k^{c}(x)$, integrating, 
and using 
\eqref{eq:diff1}-\eqref{eq:first_psi3}, we obtain, for all $k\neq N$, 
\begin{align*}
&2 \frac{\partial \mu^c_N}{\partial c} \int_{-1}^{1} \frac{\partial \psi_N^c}{\partial c}(x) \psi_{k}^c(x) dx+
\mu_N^c \int_{-1}^{1} \frac{\partial^2 \psi_N^c}{\partial c^2}(x)  \psi_k^c(x) dx\\
&=\mu_{k}^c\int_{-1}^{1} \frac{\partial^2 \psi_N^c}{\partial c^2}(x)  \psi_k^c(x) dx+ 2 \frac{\mu_k^c}{c} \int_{-1}^{1} x  \frac{\partial \psi_N^c}{\partial c}(x)
\frac{\partial\psi_{k}^c}{\partial x}(x)dx+ \frac{\mu_{k}^c}{c^2} \int_{-1}^{1} x^2 \psi_N^c(x) \frac{\partial^2\psi_{k}^c}{\partial x^2}(x) dx.
\end{align*}
Recombining and using that, for all $k\ne N$, $ \mu_{k}^c\ne\mu_N^c$ (see (3.45) in \cite{Osipov}), we obtain
\begin{align*}
&\left(\mu_N^c - \mu_{k}^c\right) \int_{-1}^{1} \frac{\partial^2 \psi_N^c}{\partial c^2}(x)  \psi_k^c(x) dx \\
&=
2  \frac{\mu_k^c}{c} \int_{-1}^{1} x  \frac{\partial \psi_N^c}{\partial c}(x) \frac{\partial\psi_{k}^c}{\partial x}(x) dx 
+
\frac{\mu_{k}^c}{c^2} \int_{-1}^{1} x^2 \psi_N^c(x)\frac{\partial^2\psi_{k}^c}{\partial x^2}(x) dx -2 \frac{\partial \mu^c_N}{\partial c} \int_{-1}^{1} \frac{\partial \psi_N^c}{\partial c}(x) \psi_{k}^c(x) dx.
\end{align*}
Using \eqref{eq:first_decomp}, (7.69)-(7.70), and Theorem 7.11 in \cite{Osipov}, yield, for all $k \not\equiv N [2]$,
$\int_{-1}^{1} \frac{\partial^2 \psi_N^c}{\partial c^2}(x)  \psi_k^c(x) dx=0,$
while, for all $k \equiv N[2]$ and $k\ne N$, using (7.69)-(7.70), Theorem 7.11, (7.99) and the eigenvalues $(\chi_N^c)_{N\in\N_0}$  of the differential operator in (1.1) in \cite{Osipov},
$$\int_{-1}^{1} \frac{\partial^2 \psi_N^c}{\partial c^2}(x)  \psi_k^c(x) dx= \frac{2}{c}\frac{\mu_k^c}{\mu_N^c - \mu_k^c} \int_{-1}^{1} x  \frac{\partial \psi_N^c}{\partial c}(x)\frac{\partial\psi_{k}^c}{\partial x}(x)dx+\Xi_{N,k},$$
\begin{equation*}
\Xi_{N,k} :=   \frac{\psi_N^{c}(1)\psi_k^{c}(1)}{c^2} \left(\frac{\mu_N^{c}\mu_k^{c}(\chi^c_k - \chi^c_N)}{\left(\mu_N^{c} - \mu_k^{c}\right)^2}
-2\frac{\mu_N^{c}\mu_k^{c}}{\left(\mu_N^{c}\right)^2 - \left(\mu_k^{c}\right)^2}\left(
2+\frac{\mu_N^c\left(2\psi_N^{c}(1)^2 - 1\right)}{\mu_N^c - \mu_k^c} \right)\right).
\end{equation*}
Differentiating  (7.114) in \cite{Osipov} in $c$ yields 	
$\int_{-1}^{1} \frac{\partial^2 \psi_N^c}{\partial c^2}(x)  \psi_N^c(x) dx=-\int_{-1}^{1} \left(\frac{\partial  \psi_N^{c}}{\partial c }(x)\right)^2 dx$. 
Also, by \eqref{b1}, for all $k\equiv N[2]$,
\begin{equation}\label{e:majr} 
\frac{\abs{\mu_N^c}}{\abs{\mu_N^c - \mu_k^c}}\le 1\ \textrm{if }k<N\ \textrm{and else } \frac{\abs{\mu_N^c}}{\abs{\mu_N^c - \mu_k^c}}\le 2.
\end{equation}
We obtain, using Lemma \ref{rem:psi} and $N\ge c$ for the first term, 
\begin{align}\notag
\left\|\frac{\partial^2  \psi_N^{c}}{\partial c^2 }\right\|_{L^{\infty}([-1,1])}  \leq & \sqrt{N+\frac{1}{2}}
\int_{-1}^{1} \left(\frac{\partial  \psi_N^{c}}{\partial c }(x)\right)^2 dx + \sum_{k \equiv N[2],\ k\neq N} \left|\Xi_{N,k}\right| \|\psi_k^c\|_{L^{\infty}([-1,1])} \\ &  + \sum_{k \equiv N[2],\ k\neq N}   \frac{2\abs{\mu_k^c}}{c\abs{\mu_N^c - \mu_k^c}} \left|\int_{-1}^{1} x  \frac{\partial \psi_N^c}{\partial c}(x) \frac{\partial\psi_{k}^c}{\partial x}(x) dx \right|   \|\psi_k^c\|_{L^{\infty}([-1,1])}. \label{eq:ref}
\end{align}
For the first term on the right-hand side of \eqref{eq:ref}, using Lemma \ref{lem:first},  we obtain
$$
\sqrt{N+\frac{1}{2}}\int_{-1}^{1} \left(\frac{\partial  \psi_N^{c}}{\partial c }(x)\right)^2 dx  \leq 8\left(c_f\left(C_1( c) +C_2( c) \right)C_3( c)\right)^2C_3( c)   \left(\frac{\Pi(c)}{c} \right)^2N^{3/2}.
$$
For the second term in \eqref{eq:ref}, using that for all $ k\equiv N [2] $, $\mu_N^c/\mu_k^c\in\R$ and \eqref{b1} we obtain
$$\abs{\Xi_{N,k}} \leq  \frac{\abs{\psi_N^{c}(1)}\abs{\psi_k^{c}(1)}}{c^2}\left( g\left( \rho_k\right)(\chi^c_k - \chi^c_N) + 2\left(
2+\frac{\left|2\psi_N^{c}(1)^2 - 1\right|\left|\mu_N^c\right|}{\left|\mu_N^c - \mu_k^c\right|} \right)   f\left( \rho_k \right)\right),$$
where $\rho_k=\mu_N^c/\mu_k^c$ when $k<N$ and $\rho_k=\mu_k^c/\mu_N^c$ when $k>N$. 
Using $N\geq c$, \eqref{e:majr}, $\abs{\chi_N^c - \chi_k^c}  \leq \abs{N-k}(k+N+1) + c^2$ (see (13) in \cite{PSWF}), \eqref{ubfg}, and $ \left|2\psi_N^{c}(1)^2 - 1\right| \leq 2N $ (by Lemma \ref{rem:psi}) for the first inequality,  $(N-k)(k+N+1)  \leq N(N+1)$ for all $0< k<N$,  \eqref{b1}, and \eqref{sumod} for the second, $(k-N)(k+N+1)=k(k+1)-N^2-N$ for the third, 
the computations in \eqref{eq:lem4}, 
$\sum_{k=1}^{\infty} k^2 2^{-k} = 6$ and  $\sum_{k=1}^{\infty} k^3 2^{-k}  = 26$, and Euclidean division for the fourth, yield
\begin{align*}
&\sum_{k \equiv N[2],\ k\neq N}\left|\Xi_{N,k}\right|  \|\psi_k^c\|_{L^{\infty}([-1,1])} \leq \frac{c_g\sqrt{4N+2}\indic\{N\equiv 0 [2]\}}{\abs{c}\pi}  \abs{\frac{\mu_N^c}{\mu_0^c}} \left(  N(N+1)  + c^2+  \frac{4 c_f}{c_g}\left(
N+1 \right)   \right) \\
&\quad +  \frac{c_g\sqrt{4N+2}}{2c^2} \sum_{k \equiv N[2],\ 0< k<N} \left(k+\frac{1}{2}\right)l(c)\abs{\frac{\mu_N^c}{\mu_k^c}} \left(  (N-k)(k+N+1)   + c^2+   \frac{4 c_f}{c_g}\left(
N+1 \right)   \right) \\ 
& \quad + \frac{c_g\sqrt{4N+2}}{2c^2} \sum_{k \equiv N[2],\ k>N} \left(k+\frac{1}{2}\right)\abs{\frac{\mu_k^c}{\mu_N^c}} \left(   \abs{N-k}(k+N+1)  + c^2+    \frac{4 c_f}{c_g}(2N +1)  \right) \\
& \leq \frac{c_g\sqrt{4N+2}}{2c^2} \left(  N(N+1)  + c^2+  \frac{4 c_f}{c_g}\left(N+1 \right)   \right) \left(  \frac{2\abs{c}}{\pi}+ \frac{l(c)N(N-1)}{4} \right)\frac{\Pi(c)}{N^2} \\
& \quad + \frac{c_g\sqrt{4N+2}}{2c^2} \sum_{k \equiv N[2],\ k>N} \frac{k+1/2}{2^{(k-N)/2}} \left(\frac{\sqrt{2\Pi(c)}}{N+2}\right)^{k-N} \ \left(    (k-N)(k+N+1)  + c^2+    \frac{4 c_f}{c_g}(2N +1)  \right) \\
&\leq \frac{c_g\sqrt{4N+2}\Pi(c)}{2c^2} \left( N(N+1) + c^2 +  \frac{4 c_f}{c_g}(N+1) \right) \left(\frac{2\abs{c}}{\pi H(c)(H(c)-1)} + \frac{l(c)}{4}     
\right) \\
& \quad + \frac{c_g\sqrt{4N+2}}{2c^2}  \frac{2\Pi(c)}{(N+2)^2}\sum_{l \equiv 0[2],\ l\ge2} \frac{l + N+1/2}{2^{l/2}} \left(c^2 +   \frac{4 c_f}{c_g}(2N +1) -N - N^2\right)\\ 
& \quad  + \frac{c_g\sqrt{4N+2}}{2c^2} \frac{2\Pi(c)}{(N+2)^2} \sum_{l \equiv 0[2],\ l\ge2}\left(l+N+\frac{1}{2}\right)(l+N)(l+N+1) \frac{1}{2^{l/2}} \\
&\leq \frac{c_g\sqrt{4N+2}\Pi(c)}{2c^2}\Bigg[ C_2( c)\left( N(N+1) + c^2 +  \frac{4 c_f}{c_g}(N+1)\right)
+C_1( c)\left(c^2 +   \frac{4 c_f}{c_g}(2N +1) -N - N^2\right)\\
&  \quad\quad\quad\quad\quad\quad\quad\quad\quad+   2N+19+\frac{85}{N+2}+\frac{246}{(N+2)^2}
\Bigg]\\
&\leq \frac{c_g \Pi(c)}{c^2} C_3( c) \Bigg[ N^{5/2}\left(C_2( c) - C_1( c)\right) 
+ N^{3/2}\left( \left(1+\frac{4 c_f}{c_g} \right)C_2( c)+\left(\frac{8 c_f}{c_g}-1 \right)C_1( c)+2\right) \\ 
&  \quad\quad\quad\quad\quad\quad\quad+ \sqrt{N}\left( (C_1(c)+C_2( c))\left( c^2+ \frac{4 c_f}{c_g}\right)+19\right) + \frac{85}{(H(c)+2)^{1/2}} + \frac{246}{(H(c)+2)^{3/2}} \Bigg].
\end{align*}
For the third term in \eqref{eq:ref},  using \eqref{eq:first_decomp}, the triangle inequality, and (7.74) in \cite{Osipov} for the first inequality and using $\abs{\mu_m^c}/\abs{\mu_m^c + \mu_k^c}\leq 1$ for the second, we obtain
\begin{align}
\left| \int_{-1}^{1} x  \frac{\partial \psi_{N}^c}{\partial c}(x) \frac{\partial  \psi_k^{c}}{\partial x}(x) dx \right| &\leq \frac{4 \left| \psi_{N}^c(1)\right| \left| \psi_{k}^c(1)\right|  }{\abs{c}} \sum_{m\neq N,\ m \equiv N[2]}   \frac{\left|\mu_N^{c}\right| \left|\mu_m^{c}\right|  \left|\psi_{m}^c(1)\right|^2 }{\left|\left(\mu_m^{c}\right)^2 - \left(\mu_N^{c}\right)^2\right|}\frac{\abs{\mu_m^c}}{\abs{\mu_m^c + \mu_k^c}}  \notag\\
&\leq \frac{4 \left| \psi_{N}^c(1)\right| \left| \psi_{k}^c(1)\right|  }{\abs{c}} \mathcal{C}(f,N,c), \notag
\end{align}
hence, using \eqref{eq:lem4} for the first inequality and 
\eqref{ubfg} and \eqref{eq:lem4}  replacing  $c_f$ by $c_h$ for the third,  
\begin{align*}
&\sum_{k \equiv N[2],\ k\neq N}\frac{2\abs{\mu_k^c}}{c\abs{\mu_N^c - \mu_k^c}} \left|\int_{-1}^{1} x  \frac{\partial \psi_{N}^c}{\partial c}(x) \frac{\partial  \psi_k^{c}}{\partial x}(x)dx \right| \|\psi_k^c\|_{L^{\infty}([-1,1])} \\
& \leq 
4c_f\sqrt{4N+2} ( C_1( c)+C_2( c)) \frac{\Pi(c)}{c^2} \sum_{k \equiv N[2],\ k\neq N}\frac{\abs{\mu_k^c}}{\abs{\mu_N^c - \mu_k^c}}  \left| \psi_{k}^c(1)\right| \|\psi_k^c\|_{L^{\infty}([-1,1])}  \\
& \leq 
4c_f\sqrt{4N+2} ( C_1( c)+C_2( c)) \frac{\Pi(c)}{c^2}   \mathcal{C}(h,N,c)\\
& \leq  4c_h c_f\sqrt{4N+2} ( C_1( c)+C_2( c))^2\frac{\Pi(c)^2}{c^2}  \leq  8c_h c_f C_3( c)    ( C_1( c)+C_2( c))^2 \frac{\Pi(c)^2}{c^2}\sqrt{N}.\quad\square
\end{align*}

\begin{lemma}\label{lem:secondH} For all $u,x_0,R>0$, $t\in\R$, $q\in\{1,\infty\}$, $\lambda$ from \eqref{eq:phi} and $N(Rx_0U)$ and $\widetilde{\mt{N}}(q)$ from \eqref{eq:Nq}, for all 
	$ N \geq N(Rx_0U)$, we have
	\begin{align*} 
	&\underset{\mt{b} \in[-R,R]^p }{\sup}\abs{ \frac{\partial^2}{\partial t^2}  \left(\left(\frac{Rx_0t}{2\pi} \right)^{p/2}  \lambda(t) \psi_{\widetilde{\mt{N}}(q)}^{Rx_0t}\left(  \frac{\mt{b}}{R}\right) \right) }  \leq \indic\{U/2\leq \abs{t} \leq U\} C_{8}(Rx_0U, p, U) N^{2+k_q/2},\\
	&C_{8}(Rx_0U, p, U):= \left(\frac{Rx_0U}{\pi}\right)^{p/2}\frac{C_3(Rx_0U)^pN(Rx_0U)^{(p-1)/(2q)}}{N(Rx_0U)^2}\Bigg(\frac{p|p-2|}{U^2}+ C_9(U)\frac{2p}{ U} + C_{10}(U)\\
	&\hspace{3.2cm}
	+   \left( \frac{ 2p}{U}+2C_9(U)\right) \frac{pC_{16}(Rx_0U)}{N(Rx_0U)^2}+ \frac{p(p-1)C_{16}(Rx_0U)}{N(Rx_0U)^2}  + pC_{11}(Rx_0U) \Bigg),\\
	&C_9(U): =\sup_{t\in [U/2,  U]}  \abs{\lambda'(t)},\ C_{10}(U) := \sup_{t\in [U/2,  U]} \abs{\lambda''(t)},\\
	&C_{11}(Rx_0U) := \frac{\left(Rx_0\right)^2\Pi(Rx_0U)}{(Rx_0U)^2}\left(C_4( Rx_0U)  + \frac{C_5(Rx_0U)}{N(Rx_0U)}+ \frac{ C_6(Rx_0U) }{N(Rx_0U)^2}+\frac{C_7(Rx_0U)}{ N(Rx_0U)^{5/2}} \right),\\
	&C_{16}(Rx_0U) :=  2c_f Rx_0 \left( C_1(Rx_0U)+ C_2(Rx_0U) \right)C_3(Rx_0U)\frac{\Pi(Rx_0U)}{Rx_0U}.
	\end{align*} 
\end{lemma}
\noindent {\bf Proof.}
Let $q=1$. By 
symmetry, we take $t\in[U/2,U]$, $\mt{b} \in [-R,R]^p$, and $c>0$. We have 
\begin{align*}
&R(t,\mt{b}) := \abs{ \frac{\partial^2}{\partial t^2}  \left(\left(\frac{R x_0t}{2\pi} \right)^{p/2}  \lambda(t) \psi_{\widetilde{\mt{N}}(q)}^{Rx_0t}\left(  \frac{\mt{b}}{R}\right) \right)}\\
&\le  \left(\frac{Rx_0}{2\pi}\right)^{p/2} t^{p/2} \Bigg [ \left(\frac{p\left|p-2\right|}{4t^2} \lambda(t)   + \frac{p}{t}\left|\lambda'(t)\right|   + \left|\lambda''(t)\right|\right)  \abs{ \psi_{\widetilde{\mt{N}}(q)}^{Rx_0t}\left(  \frac{\mt{b}}{R}\right) }  \\  
&  \quad \quad+ Rx_0\left(  \frac{p}{t} \lambda(t) + 2|\lambda'(t)| \right) \abs{ \left.\frac{\partial  \psi_{\widetilde{\mt{N}} (q)}^{c}}{\partial c}\right|_{c=Rx_0t} \left(  \frac{\mt{b}}{R}\right)}    +     \left(Rx_0\right)^2 \lambda(t) \abs{ \left. \frac{\partial^2 \psi_{\widetilde{\mt{N}}(q)}^{c}}{\partial c^2} \right|_{c=Rx_0t}\left(  \frac{\mt{b}}{R}\right) } \Bigg],\\
&\frac{\partial  \psi_{\widetilde{\mt{N}}(q)}^{c}}{\partial c} \left(  \frac{\mt{b}}{R}\right)  =\sum_{j=2}^p \psi_{N}^{c} \left(  \frac{\mt{b}_1}{R}\right)   \frac{\partial  \psi_{N(Rx_0U) }^{c}}{\partial c} \left(  \frac{\mt{b}_j}{R}\right) \prod_{\underset{l\neq j}{l=2}}^p \psi_{ N(Rx_0U) }^{c} \left(  \frac{\mt{b}_l}{R}\right)\\
& \hspace{2.8cm} +   \frac{\partial  \psi_{N}^{c}}{\partial c} \left(  \frac{\mt{b}_1}{R}\right) \prod_{l=2}^p \psi_{ N(Rx_0U) }^{c} \left(  \frac{\mt{b}_l}{R}\right) ,  \\
& \frac{\partial^2 \psi_{\widetilde{\mt{N}}(q)}^{c}}{\partial c^2} \left(  \frac{\mt{b}}{R}\right) = 2\sum_{j=2}^p \frac{\partial  \psi_{ N(Rx_0U)}^{c}}{\partial c} \left(  \frac{\mt{b}_j}{R}\right)\frac{\partial  \psi_{N}^{c}}{\partial c} \left(  \frac{\mt{b}_1}{R}\right)  \prod_{\underset{l\neq j}{l=2}}^p \psi_{N(Rx_0U)}^{c} \left(  \frac{\mt{b}_l}{R}\right)  \\ 
& +  \sum_{k=2}^p\sum_{\underset{j \neq k}{j=2}}^p \psi_{N}^{c} \left(  \frac{\mt{b}_1}{R}\right) \frac{\partial  \psi_{ N(Rx_0U)}^{c}}{\partial c} \left(  \frac{\mt{b}_j}{R}\right)\frac{\partial  \psi_{ N(Rx_0U) }^{c}}{\partial c} \left(  \frac{\mt{b}_k}{R}\right)  \prod_{\underset{l\neq j, l\neq k}{l=2}}^p \psi_{N(Rx_0U)}^{x_0t} \left(  \frac{\mt{b}_l}{R}\right) 
\\ 
& +   \frac{\partial^2  \psi_{N}^{c}}{\partial c^2} \left(  \frac{\mt{b}_1}{R}\right) \prod_{l=2}^p \psi_{N(Rx_0U) }^{c} \left(  \frac{\mt{b}_l}{R}\right)
+ \sum_{j=2}^p  \psi_{N}^{c} \left(  \frac{\mt{b}_l}{R}\right)  \frac{\partial^2  \psi_{N(Rx_0U) }^{c}}{\partial c^2} \left(  \frac{\mt{b}_j}{R}\right) \prod_{\underset{l\neq j}{l=2}}^p \psi_{ N(Rx_0U) }^{c} \left(  \frac{\mt{b}_l}{R}\right). 
\end{align*}
We conclude using 
$N\geq  Rx_0U$ (by the discussion before Lemma \ref{lem:ratio}), the third assertion of Lemma \ref{rem:psi}, and Lemma \ref{lem:second}. The case $q=\infty$ is obtained with $N(Rx_0U)=N$.\hfill $\square$  
\begin{lemma}\label{lem:cst_f_0}
	For all $R, x_0>0$, $2\sigma>k_q+4$, $q\in \{1,\infty\}$, $2\tau \geq  \left(3e^{1/2} Rx_0/4\right)\vee 1$, 
	we have
	\begin{align*}&\int_{\R}e	^{- 2\tau\abs{t}} \sum_{\boldsymbol{m}\in\N_0^p}  |\boldsymbol{m}|_{q}^{2\sigma}\left(c^P_{\mt{m}}(t)\right)^2  dt \leq \frac{C_{12}(\sigma,p)}{\tau p^{2\sigma/q}},\\
	&C_{12}(\sigma,p):=  \Gamma(2\sigma+p+1/2)\left( \frac{2^{p-1}p }{2\sigma+p}\left(\frac{8}{3e^{1/2}}\right)^{2\sigma+p}  +\frac{\pi e^{3}p2^{p}\sqrt{3}}{9} \right). 
	\end{align*}
\end{lemma}
\noindent {\bf Proof.} 
When $q=1$, we use $  \abs{\mt{m}}_{1} \leq p \abs{\mt{m}}_{\infty} $. Let $q=\infty$, $R, x_0,\sigma,\tau$ as in the lemma.  
Because $P_0= \indic\{\abs{\cdot}_{\infty}\leq 1\}/2^{p/2} $, for all $m\in\N_0$, $ \left| \braket{P_0,\psi_m^{c}}_{L^2([-1,1])}\right| \leq 1$, and, for all $m > \abs{c}$,  
$\left|\braket{
	P_0,\psi_m^{c}}_{L^2([-1,1])}\right| \leq \abs{\mu_m^c}/\sqrt{2} $ (see Proposition 3 and (13) in \cite{PSWF}) 
we obtain, for all $t\neq0$, 
\begin{align}\label{es}
&\sum_{\boldsymbol{m}\in\N_0^p} |\boldsymbol{m}|_{\infty}^{2\sigma} \left(c^P_{\mt{m}}(t)\right)^2  \leq   \sum_{ \abs{\boldsymbol{m}}_{\infty}  \le Rx_0\abs{t}  }   |\boldsymbol{m}|_{\infty}^{2\sigma} \indic\{Rx_0\abs{t} \geq 1\} +   \sum_{ \abs{\boldsymbol{m}}_{\infty} > Rx_0\abs{t} } \frac{  |\boldsymbol{m}|_{\infty}^{2\sigma} \abs{\mu_{|\boldsymbol{m}|_{\infty}}^{Rx_0t}}^2}{2}  .
\end{align} 
Using \eqref{es}, Lemma \ref{upper_bound}, and $\sum_{ |\boldsymbol{m}|_{\infty} =k}1 \leq p(k+1)^{p-1}$ for the first inequality, $m+1\leq 2m$ when $m\geq 1$ for the second, and $2m+1\leq 3m$,   $(Rx_0t+1)^{2\sigma+p}\leq(2Rx_0t)^{2\sigma+p}  $ when $m, Rx_0t\geq 1$, and 
(1.3) in \cite{inq_gamma2}
for the third, we have 
\begin{align*}
&\int_{\R}e^{- 2\tau\abs{t}}\sum_{\boldsymbol{m}\in\N_0^p} |\boldsymbol{m}|_{\infty}^{2\sigma} \left(c^P_{\mt{m}}(t)\right)^2 dt \\
& \leq \int_{0}^{\infty}2pe^{-2\tau t}\left(\sum_{  m\leq Rx_0t  }(m+1)^{p-1} m^{2\sigma} \indic\{Rx_0t \geq 1\} + \frac{\pi e^3}{9}   \sum_{  m>Rx_0t  } (m+1)^{p-1}m^{2\sigma} \left( \frac{eRx_0t}{4m}\right)^{2m} \right)dt\\   
& \leq \int_{0}^{\infty} 2^ppe^{-2\tau t} \int_{1}^{Rx_0t+1}u^{2\sigma+p-1}du \indic\{Rx_0t \geq 1\}dt + \frac{\pi e^3p2^{p}}{9}  \sum_{  m\ge1  } m^{2\sigma+p-1}\left( \frac{eRx_0}{4m}   \right)^{2m}    \int_{0}^{\infty} \frac{t^{2 m}}{ e^{2\tau t}}   dt  \\ 
& \leq\frac{2^{2(\sigma+p)}p}{2\sigma+p} \int_{1/(Rx_0)}^{\infty} e^{-2\tau t}(Rx_0t)^{2\sigma+p} dt + \frac{\pi e^{3}p2^{p}\sqrt{3}}{9\tau}  \sum_{  m\ge1  } m^{2\sigma+p-1/2} e^{2m\ln( 3Rx_0/(8\tau))}\\
& \leq \frac{2^{p-1}p \Gamma(2\sigma+p+1) }{(2\sigma+p)\tau}  \left(\frac{8}{3e^{1/2}}\right)^{2\sigma+p}+\frac{\pi e^{3}p2^{p}\sqrt{3}}{9\tau} \int_{0}^{\infty} e^{-t}t^{2\sigma+p-1/2} dt  
\leq  \frac{C_{12}(\sigma,p)}{\tau p^{2\sigma/q}}.\quad\square
\end{align*}

\begin{lemma}\label{lem:I12}
	For all $ N \geq H(Rx_0U)$, $R, U>0$, $q\in\{1,\infty\}$, 
	and $F$ from \eqref{eq:H_N}, we have
	\begin{align}
	I_1& := \int_{[-1,1]^p} \int_{\R} \left|\partial_t \mathcal{F}\left[F\right](t, x_0t\mt{x}) \right|^2   d\mt{x} dt \leq R^pC_{17}(Rx_0U,p,U)  N^2  \rho_{\widetilde{\mt{N}}(q)}^{W_{[-1,1]},Rx_0U} \label{eq:I1}\\
	C_{17}(Rx_0U,p,U) &:= C_{15}(Rx_0U,p,U) + \frac{2pUC_{16}(Rx_0U)^2}{N(Rx_0U)},\notag\\
	C_{15}(Rx_0U,p,U) &:=  \frac{25p^2}{8U}  \left(1+ \frac{2(Rx_0U)^2}{3^{3/2}}\right)^{4}  +  \frac{U  C_{9}(U)^2}{N(Rx_0U)^2} +  \frac{5p C_{9}(U)\ln(2)}{2N(Rx_0U)}  \left(1+ \frac{2(Rx_0U)^2}{3^{3/2}}\right)^{2} ,\notag\\
	I_2 &:= \int_{[-1,1]^p}  \int_{\R} \left|\mathcal{F}\left[F\right](t, x_0t\mt{x}) \right|^2  d\mt{x} dt   \leq R^p U \rho_{\widetilde{\mt{N}}(q)}^{W_{[-1,1]},Rx_0U}. \label{eq:I2}
	\end{align}
\end{lemma}

\noindent {\bf Proof.} 
Let $ N \geq H(Rx_0U)\ge2$. For simplicity of notations, we omit $W_{[-1,1]}$ from $\rho$.  We have
\begin{align}\notag
\mathcal{F}\left[F\right](t, x_0t\mt{x}) 
& = \left(\frac{Rx_0\abs{t}}{2\pi}\right)^{p/2} R^{p/2} \lambda(t)\mathcal{F}_{Rx_0t}\left[ \psi_{\widetilde{\mt{N}}(q)}^{Rx_0t}\right](\mt{x})  \notag \\
& = R^{p/2} i^{\left|\widetilde{\mt{N}}(q)\right|_1}\lambda(t) \sqrt{\rho_{\widetilde{\mt{N}}(q)}^{Rx_0t} } \psi_{\widetilde{\mt{N}}(q)}^{Rx_0t} \left(\mt{x}\right) \ \left(\text{because}\  \mu_{m}^{Rx_0t}=i^{m}\left(\frac{2\pi}{Rx_0\abs{t}}\right)^{1/2} \sqrt{\rho_{m}^{Rx_0t}} \right).\label{eq:bound_2}
\end{align} 
This yields
\begin{align*}
I_1 \leq   & R^{p} \int_{\R} \int_{[-1,1]^p}  \left( \left(\frac{\mathrm{d}\sqrt{\rho_{\widetilde{\mt{N}}(q)}^{Rx_0t}}}{\mathrm{d}t}   \lambda(t)    +   \lambda'(t)\sqrt{\rho_{\widetilde{\mt{N}}(q)}^{Rx_0t}}  \right)  \psi_{\widetilde{\mt{N}}(q)}^{Rx_0t}(\mt{x} )  +  \sqrt{\rho_{\widetilde{\mt{N}}(q)}^{Rx_0t}} \lambda(t) \frac{\partial \psi_{\widetilde{\mt{N}}(q)}^{Rx_0t}(\mt{x} )}{\partial t} \right)^2  dtd\mt{x} . 
\end{align*}
Using (7.114) in \cite{Osipov},  cross-products terms are zero and $I_1\le R^p(I_{11}+I_{12})$, where
\begin{align*}
I_{11} &=  \int_{\R} \left( \lambda(t)^2 \left( \frac{\mathrm{d}\sqrt{\rho_{\widetilde{\mt{N}}(q)}^{Rx_0t}}}{\mathrm{d}t}\right)^2 +   \lambda'(t)^2 \rho_{\widetilde{\mt{N}}(q)}^{Rx_0t} + 2 \lambda(t) |\lambda'(t)|\sqrt{\rho_{\widetilde{\mt{N}}(q)}^{Rx_0t}} \frac{\mathrm{d} \sqrt{\rho_{\widetilde{\mt{N}}(q)}^{Rx_0t}}}{\mathrm{d} t} 
\right) dt , \\ 
I_{12}&=     \int_{\R}  \lambda(t)^2\rho_{\widetilde{\mt{N}}(q)}^{Rx_0t} \left( \int_{[-1,1]^p}\left(\frac{ \partial \psi_{\widetilde{\mt{N}}(q)}^{Rx_0t}(\mt{x} )} {\partial t}\right)^2 d\mt{x} \right)    dt . 
\end{align*}
Then, using (7.100) in \cite{Osipov} for the second equality 
yields, for all $t\ne0$,
\begin{align*}
\frac{\mathrm{d}\sqrt{\rho_{N}^{Rx_0t}}}{\mathrm{d} t}
& = \frac{x_0R}{2\sqrt{\rho_{N}^{Rx_0t}}}\left.\frac{\mathrm{d}\rho_{N}^{c}}{\mathrm{d}c}\right|_{c=Rx_0t}= \frac{\sqrt{\rho_{N}^{Rx_0t}}}{|t|}\left(\psi_N^{Rx_0t}(1)\right)^2, 
\end{align*}
in particular $\rho_{N}^{Rx_0t}$ is increasing in $t$ and, by the first assertion of Lemma \ref{rem:psi},  
\begin{align}
\forall U/2 \leq \abs{t} \leq U,\ 
\frac{\mathrm{d} \sqrt{\rho_{N}^{Rx_0t}}}{\mathrm{d} t} 
& \leq \frac{(N+1/2)\sqrt{ \rho_{N}^{Rx_0t}}}{\abs{t}}\left(1+ \frac{2(Rx_0U)^2}{3^{3/2}}\right)^{2}  \label{I_1_A}. 
\end{align}
When $q=1$, using 
$N+1/2\leq 5N/4$ for  all $N\geq 2$ and 
\begin{align*}
\notag \frac{\mathrm{d} \sqrt{\rho_{\widetilde{\mt{N}}(q)}^{Rx_0t}}}{\mathrm{d} t} & =(p-1) \left(\sqrt{\rho_{N(Rx_0U)}^{Rx_0t}}\right)^{p-2} \sqrt{\rho_{N}^{Rx_0t}}\left(\frac{\mathrm{d} \sqrt{\rho_{N(Rx_0U)}^{Rx_0t}}}{\mathrm{d} t}\right) + \left(\sqrt{\rho_{N(Rx_0U)}^{Rx_0t}}\right)^{p-1}   \left(\frac{\mathrm{d} \sqrt{\rho_{N}^{Rx_0t}}}{\mathrm{d} t} \right),
\end{align*}
we have
\begin{equation}
\frac{\mathrm{d} \sqrt{\rho_{\widetilde{\mt{N}}(q)}^{Rx_0t}}}{\mathrm{d} t} \leq \frac{5pN}{4\abs{t}}  \left(1+ \frac{2(Rx_0U)^2}{3^{3/2}}\right)^{2} \sqrt{\rho_{\widetilde{\mt{N}}(q)}^{Rx_0t}}.\label{eq:sameq}
\end{equation}
The same inequality holds for $q=\infty$ (there $N=N(Rx_0U)$).
This yields, for all $q\in\{1,\infty\}$, 
\begin{align*}
I_{11}& \leq \left(\frac{25p^2N^2}{8} \int_{U/2}^{U} \frac{dt}{t^2} \left(1+ \frac{2(Rx_0U)^2}{3^{3/2}}\right)^{4}  +  U  C_{9}(U)^2+  \frac{5pN C_{9}(U)}{2}  \left(1+ \frac{2(Rx_0U)^2}{3^{3/2}}\right)^{2} \int_{U/2}^{U} \frac{dt}{t}\right)\rho_{\widetilde{\mt{N}}(q)}^{Rx_0U} \\
&\leq  C_{15}(Rx_0U,p,U) N^2  \rho_{\widetilde{\mt{N}}(q)}^{Rx_0U}.
\end{align*}
Then, by (7.114) in \cite{Osipov} 
and  Lemma \ref{lem:first}, 
we have, for all $  U/2 \leq \abs{t} \leq U $,
\begin{align*}
\int_{[-1,1]^p} \left(\frac{\partial  \psi_{\widetilde{\mt{N}}(q)}^{Rx_0t}(\mt{x} )}{\partial t}\right)^2d\mt{x}   = &  (Rx_0)^2 \int_{[-1,1]} (p-1)\left(\left. \frac{\partial  \psi_{N(Rx_0U)}^{c}(x)}{\partial c}\right|_{c=Rx_0t}\right)^2dx \\
&+ (Rx_0)^2 \int_{[-1,1]} \left(\left. \frac{\partial  \psi_N^{c}(x )}{\partial c}\right|_{c=Rx_0t}\right)^2 dx 
\\  
 \leq &  2p\left(C_{16}(Rx_0U)\right)^2N \quad (\text{using }  N \geq N(Rx_0U)).
\end{align*} 
The same holds for $q=\infty$  (there $N=N(Rx_0U)$). This and 
$ N \geq N(Rx_0U)$ yield \eqref{eq:I1}.\\
\eqref{eq:I2} follows from 
\eqref{eq:bound_2} and the fact that  $ c\in (0,\infty) \mapsto \rho_{m}^{c} $ is nondecreasing.\hfill $\square$

\renewcommand{\theequation}{C.\arabic{equation}}
\renewcommand{\thelemma}{C.\arabic{lemma}}
\renewcommand{\thecorollary}{C.\arabic{corollary}}
\renewcommand{\thedefinition}{C.\arabic{definition}}
\renewcommand{\theproposition}{C.\arabic{proposition}}
\renewcommand{\theremark}{C.\arabic{remark}}
\renewcommand{\thetheorem}{C.\arabic{theorem}}
\renewcommand{\theassumption}{C.\arabic{assumptio}}
\renewcommand{\thesubsection}{C.\arabic{subsection}}
\setcounter{equation}{0}  
\setcounter{lemma}{0}
\setcounter{corollary}{0}
\setcounter{proposition}{0}
\setcounter{remark}{0}
\setcounter{definition}{0}
\setcounter{lemma}{0}
\setcounter{theorem}{0}
\setcounter{assumption}{0}
\setcounter{subsection}{0}
\setcounter{footnote}{0}
\setcounter{figure}{0}

\section{Estimation of the marginal $f_{\mt{\beta}}$}\label{Amarginals}
For all $\left(\omega_m\right)_{m\in\N_0}$ increasing, $\omega_0 = 1$, $l,M>0$, $q\in\{1,\infty\}$, consider 
$$ \mathcal{H}^{q,\omega}_{w,W}(l,M) := \left\{ f:  \norm{f}_{L^2\left(w\otimes W^{\otimes p}\right)} \leq M, \  
\sum_{k\in\N_0} \omega_{k}^2 \left\|\theta_{q,k} \right\|^2_{L^2(\R)} \leq 2\pi l^2  \right\}.$$
For brevity, 
we present the slow rates and 
the estimator $\widehat{f}_{\beta}^{q,N, \epsilon}:=\sum_{|\mt{m}|_{q}\leq N(\epsilon)} \widehat{c}_{\mt{m}}(\epsilon)\varphi_{\mt{m}}^{W,\epsilon x_0}/\sigma_{\mt{m}}^{W,\epsilon x_0}$.
It is based on $f_{\mt{\beta}} = \mathcal{F}_{1\mathrm{st}}\left[f_{\alpha,\mt{\beta}}\right](0,\cdot_2)$. 
\begin{proposition}\label{upper:marginal}  Let $W=W_{[-R,R]}$. For all $q\in \{1,\infty\}$, $l,M,R>0$, $ \sigma > 2 $, $\mathbb{S}_{\beta}\subseteq[-R,R]^p$, 
	$\underline{N}$ solution of $ 2(1+\sigma)k_q\underline{N}\ln(\underline{N}) + p(1-\sigma)\ln(\underline{N}) + \ln(\omega_{\underline{N}}^2) = \ln(n_e)$,  $\epsilon=\theta/\omega_{\underline{N}}$, $(\omega_k)_{k\in \N_0} =\left(k^{\sigma}\right)_{k\in \N_0} $, and $w$ such that $\int_{\R}a^2/w(a)da<\infty$, 
	we have
	\begin{equation*}
\left(   \frac{\ln(n_e)}{\ln_2(n_e)}\right)^{2\sigma}	\underset{
		f_{\mt{\beta}} \in \mathcal{H}^{q,\omega}_{w,W}(l,M)\cap{D},\ {f_{\mt{X}|\mathcal{X}} \in \mathcal{E}}}{\sup} \E \left[ \left\|  \widehat{f}^{q,N,\epsilon}_{\mt{\beta}}-f_{\mt{\beta}} \right\|_{L^2(\R^{p})}^2
	\right]  = O_p(1). 
	\end{equation*}
\end{proposition}
\noindent {\bf Proof.}
We assume $ f_{\mt{X}|\mathcal{X}} $ is known. The general case can be handled like in the proof of 
(T\ref{theo:compact}.\ref{t_comp_1}). 
Use $ f_{\mt{\beta}}^{\epsilon} :=  \mathcal{F}_{1\mathrm{st}}\left[f_{\alpha,\mt{\beta}}\right](\epsilon,\cdot) $  and define  $f_{\mt{\beta}}^{q,\epsilon,N} $ like $ \widehat{f}_{\mt{\beta}}^{q,\epsilon,N}$ with  $ \widetilde{c}_{\mt{m}}(t) $ (see Lemma \ref{lem:def_est}) instead of $ \widehat{c}_{\mt{m}}(t)$. Use 
$ \left\|  \widehat{f}^{q,N,\epsilon}_{\mt{\beta}}-f_{\mt{\beta}} \right\|_{L^2(\R^{p})}^2 \leq 3 \sum_{j=1}^{3} \left\|R_j\right\|_{L^2(\R^{p})}^2$, where $ R_1 :=  \widehat{f}_{\mt{\beta}}   - f_{\mt{\beta}}^{q,N,\epsilon} $,  $ R_2 :=  f_{\mt{\beta}}^{q,N,\epsilon}- f_{\mt{\beta}}^{\epsilon} $, and $ R_3 := f_{\mt{\beta}}^{\epsilon}- f_{\mt{\beta}}$. Let $n\geq e^e$ large enough so that $N\geq 1\bigvee ((\sigma-1)p-\sigma)/(2k_q)$. By similar arguments from \eqref{eq:R_0_w}, \eqref{eq:Nk}, 
$N\leq \underline{N}$, and $(\underline{N}+1)^{2k_q\underline{N}}\le e^{2k_q}\underline{N}^{2k_q\underline{N}}$, we have
\begin{align}
\E\left[\left\| R_1 \right\|_{L^2(\R^{p})}^2 \right]\leq &   \frac{Q_q c_{\mt{X}}e^{2k_q}}{\pi^pn} \epsilon^p \underline{N}^{p}\left(1\bigvee \frac{\theta\underline{N}}{\epsilon}\right)^{2k_q\underline{N}} =  \frac{Q_q c_{\mt{X}}e^{2k_q}\theta^p}{\pi^pn} \underline{N}^{p(1-\sigma)+2(1+\sigma)k_q\underline{N}}.\label{eq:marg_R1}
\end{align}
We also obtain $\left\|R_2\right\|^2_{L^2(\R^{p})} \leq 2\pi l^2 /\omega_{\underline{N}}^2$ and 
\begin{align}
\left\|R_3\right\|^2_{L^2(\R^{p})}  
\leq  & \int_{[-R,R]^p} \left(\int_{\R}  \abs{e^{i\epsilon a} -1 } f_{\alpha,\mt{\beta}}(a,\mt{b}) da  \right)^2 d\mt{b} \notag \\
\leq &  \epsilon^2 \int_{[-R,R]^p} \left(\int_{\R} \abs{a} f_{\alpha,\mt{\beta}}(a,\mt{b}) da \right)^2 d\mt{b} 
\leq    \frac{\theta^2M^2}{\omega_{\underline{N}}^2} \int_{\R}\frac{a^2}{w(a)}da<\infty
. \label{eq:marg_R2}
\end{align}
Then, using 
$ \ln(n) = 2(1+\sigma)k_q\underline{N}\ln(\underline{N}) + p(1-\sigma)\ln(\underline{N}) + \ln(\omega_{\underline{N}}^2) \geq 2\sigma k_q\underline{N}\ln(\underline{N})$ and  $\mathcal{W}(x) \leq  \ln(x +1 )$, we have $\underline{N} \leq \ln(n)/(2\sigma k_q\ln_2(n_e)(1+\ln(1+2\sigma k_q/e))$.
The result follows from the definition of $\underline{N}$, \eqref{eq:marg_R1}, and  \eqref{eq:marg_R2}. \hfill $\square$\vspace{0,3cm}

Similar ideas apply for the estimation of $f_{\mt{\beta}_j}$ for $j=1,\dots,p$.

\renewcommand{\theequation}{D.\arabic{equation}}
\renewcommand{\thelemma}{D.\arabic{lemma}}
\renewcommand{\thecorollary}{D.\arabic{corollary}}
\renewcommand{\thedefinition}{D.\arabic{definition}}
\renewcommand{\theproposition}{D.\arabic{proposition}}
\renewcommand{\theremark}{D.\arabic{remark}}
\renewcommand{\thetheorem}{D.\arabic{theorem}}
\renewcommand{\theassumption}{D.\arabic{assumptio}}
\renewcommand{\thesubsection}{D.\arabic{subsection}}
\setcounter{equation}{0}  
\setcounter{lemma}{0}
\setcounter{corollary}{0}
\setcounter{proposition}{0}
\setcounter{remark}{0}
\setcounter{definition}{0}
\setcounter{lemma}{0}
\setcounter{theorem}{0}
\setcounter{assumption}{0}
\setcounter{subsection}{0}
\setcounter{footnote}{0}
\setcounter{figure}{0}

\section{Talagrand inequality for complex functions}\label{ATalagrand}

\begin{lemma}\label{Talagrand}
	Let $X_1,\dots,X_n$ $n$ independent random vectors, $\Lambda:= (\sqrt{1 +\cdot}-1)\wedge 1$, $\mathcal{U}$ a countable set of complex measurable functions, and, for all $u\in\mathcal{U}$, $\nu_n(u):= \sum_{i=1}^n \left(u(X_i) - \E\left[u(X_i)\right]\right)/n$. If there exist $M, H, v$ $>0$ such that 
	\begin{align*}
	& \sup_{u \in \mathcal{U}} \left\| u \right\|_{L^{\infty}(\R^p)} \leq M,\ \E\left[ \underset{u \in \mathcal{U}}{\sup} \left| \nu_n(u) \right| \right] \leq H,\ 
	\sup_{u \in \mathcal{U}} \frac{1}{n}\sum_{i=1}^{n}  \text{Var}\left(\mathfrak{R}(u(X_i))\right) \bigvee
	\text{Var}\left(\mathfrak{I}(u(X_i))\right)  \leq v, 
	\end{align*}
	then, for all $\eta >0$,
	\begin{align*}
	\E\left[ \left( \sup_{u \in \mathcal{U}}  \left| \nu_n(u) \right|^2 - 4(1+2\eta)H^2 \right)_+  \right]
	\le48\left( \frac{v}{n}e^{-\eta \frac{nH^2}{6v}}+ \frac{294 M^2}{ \Lambda(\eta)^2 n^2} e^{-  \frac{\sqrt{2} \Lambda(\eta) \sqrt{\eta} }{42} \frac{nH}{M}} \right).
	\end{align*}
\end{lemma}
\noindent {\bf Proof.}  The result follows from Theorem 7.3 in \cite{comte2018regression} and 
\begin{align*}
&\E\left[ \left( \sup_{u \in \mathcal{U}}  \left| \nu_n(u) \right|^2 - 4(1+2\eta)H^2 \right)_+  \right] \leq  \E\left[ \left( \sup_{u \in \mathcal{U}}  \mathfrak{R}(\nu_n(u))^2 + \sup_{u \in \mathcal{U}}   \mathfrak{I}(\nu_n(u))^2 - 4(1+2\eta)H^2 \right)_+  \right] \\
&\leq  \E\left[ \left( \sup_{u \in \mathcal{U}}  \mathfrak{R}(\nu_n(u))^2 - 2(1+2\eta)H^2 \right)_+ \right] + \E\left[ \left(\sup_{u \in \mathcal{U}}   \mathfrak{I}(\nu_n(u))^2- 2(1+2\eta)H^2 \right)_+   \right].\quad\quad \square
\end{align*}

\renewcommand{\theequation}{E.\arabic{equation}}
\renewcommand{\thelemma}{E.\arabic{lemma}}
\renewcommand{\thecorollary}{E.\arabic{corollary}}
\renewcommand{\thedefinition}{E.\arabic{definition}}
\renewcommand{\theproposition}{E.\arabic{proposition}}
\renewcommand{\theremark}{E.\arabic{remark}}
\renewcommand{\thetheorem}{E.\arabic{theorem}}
\renewcommand{\theassumption}{E.\arabic{assumptio}}
\renewcommand{\thesubsection}{E.\arabic{subsection}}
\setcounter{equation}{0}  
\setcounter{lemma}{0}
\setcounter{corollary}{0}
\setcounter{proposition}{0}
\setcounter{remark}{0}
\setcounter{definition}{0}
\setcounter{lemma}{0}
\setcounter{theorem}{0}
\setcounter{assumption}{0}
\setcounter{subsection}{0}
\setcounter{footnote}{0}
\setcounter{figure}{0}

\section{Approximation by PSWF in Sobolev ellipsoids.}\label{ASob}
For all $\sigma, s,l>0$ and $q\in\{1,\infty\}$,  denote by $\left(\phi_{\mt{m}}\left(\cdot/R\right)\right)_{\mt{m}\in\Z^p} := \left(  e^{i\pi \mt{m}^{\top} \cdot/R}/(2R)^{p/2} \right)_{\mt{m}\in\Z^p}$, 
$\mathcal{F}[f](\star,\mt{k}) := \int_{\R} e^{i \star a} \int_{[-R,R]^p} e^{i\pi \mt{k}^{\top} \mt{b}/R} \ f(a,\mt{b}) da d\mt{b}/(2R)^{p/2} $, and  
\begin{align*} 
H^{q,s,\sigma}(l):=\left\{ f:   \int_{\R} \sum_{\mt{k} \in \Z^p} \left|\mathcal{F}[f](t,\mt{k})\right|^2  (1\vee t^{2s})  dt \bigvee \int_{\R} \sum_{\mt{k} \in \Z^p} \left|\mathcal{F}[f](t,\mt{k})\right|^2  \left(1\vee|\mt{k}|_{q}^{2\sigma}\right)dt  \leq 2\pi l^2  \right\}.
\end{align*}
Denote, for all $N\in\N$ and $c\ne0$,  
by $P^N_{c}$ (resp. $\mathcal{E}^N$) the projector in $L^2\left(W_{[-R,R]}^{\otimes p}\right) $ onto the vector space spanned by $\left(\psi_{\mt{m}}^{c}\left(\cdot /R\right)/R^{p/2} \right)_{|\mt{m}|_{\infty}<N}$ (resp. $\left(\phi_{\mt{m}}\left(\cdot/R\right)\right)_{{|\mt{m}|_{\infty}<N}}$).  For all $ t\neq	 0$ and $\left(n,m,N, \widetilde{N}\right)\in \N_0^4$, denote 
by $ \varphi^{t}: = \mathcal{F}_{1\mathrm{st}}\left[f\right](t, \cdot_2)$, $\beta_n^{m}(t):= \braket{\psi_{m}^{t} ,P_n }_{L^2([-1,1])} $, $J_j$ the Bessel function of the first kind and order $j>-1$, $ K^{N,\widetilde{N}}_{t}: =\left\| \mathcal{E}^{ \widetilde{N}} \varphi^t - P^N_{x_0t} \mathcal{E}^{ \widetilde{N}} \varphi^t  \right\|^2 $, and 
$I_{N, \widetilde{N}}: =\sum_{\mt{k}\in\Z^p:\ \abs{\mt{k}}_{\infty}< \widetilde{N}}\sum_{\abs{\mt{m}}_{\infty}\ge N}
\abs{\braket{\phi_{\mt{k}}\left(\cdot/R\right) , \psi_{\mt{m}}^{x_0\star}\left(\cdot/R\right)} }^2$. 
\begin{proposition}\label{lem:comp_Sobol} 
	For all $\sigma ,l,M,R>0$, 
	$q\in \{1,\infty\}$,  and 
	$s \geq \sigma+p/2$, 
	we have, for all $N\geq10$, 
	\begin{align}\label{eq:fin1}
	&\int_{\R} \left\|  \mathcal{F}_{1\mathrm{st}}\left[f\right](t,\cdot_2) - P^N_{x_0t}  \mathcal{F}_{1\mathrm{st}}\left[f\right](t,\cdot_2) \right\|^2   dt \leq \frac{2\pi A l^2 }{N^{2\sigma} 	},  
	\end{align}
	\begin{align*}
	A:=& 2\left(  \left(\dfrac{1}{1/(\pi e)-1/10}\right)^{2\sigma} + c \left(\frac{p +2\sigma}{be}\right)^{p+2\sigma}  + \left(\frac{2eRx_0}{\pi}\right)^p(e^2x_0)^{2\sigma}  \right),\\
	a:=	& \frac{\sqrt{5e^3}(e^2+ 1/e^2)^{5/8}}{3(\ln(2)+2)2^{11/4}}, b:=  p\left(\frac{5}{8}\ln\left(\frac{21}{10}\right)-\frac{1}{e}\right),\\
	c :=&\frac{p(4R^2/(\pi e))^p }{2R} \left( a^{2p}  \frac{8^p(p-1)^{p-1} }{(3p)^{p}e^{p-1}}+   \left(\frac{(2p-1)8}{5p e}\right)^{2p-1} \frac{5^{p-1}8}{p16^p\ln\left(21/10\right)} \right).
	\end{align*}
\end{proposition}
Proposition \ref{lem:comp_Sobol} is an analogue of \eqref{eq:R_2_w} with $N$ constant and  $\left(\omega_k\right)_{k\in\N_0} =  \left(k^{\sigma}\right)_{k\in\N_0}$. It shows that the approximation error  when we use a truncated series expansion in the PSWF basis is of order $N^{-2\sigma}$ whether we work on the class $H^{q,s,\sigma}(l)$ or $ \mathcal{H}^{q,\phi,\omega}_{w,W_{[-R,R]}}(l, M)$ with $\phi=1\vee \abs{\cdot}^s$. 
\eqref{eq:R3_w} can be obtained using the first inequality in the definition of $H^{q,s,\sigma}(l)$ and, for all $t\neq0$, $\sum_{\mt{k} \in \Z^p} \left|\mathcal{F}[f](t,\mt{k})\right|^2 = \left\|\mathcal{F}_{1\text{st}}[f](t,\cdot_2)\right\|^2_{L^2\left(W_{[-R,R]}^{\otimes p}\right)}= \sum_{\mt{m} \in \N_0^p} \left|b_{\mt{m}}(t)\right|^2$. 
Thus 
(T\ref{theo:compact}.\ref{t_comp_1}) also holds for functions in the intersection of $H^{q,s,\sigma}(l)\bigcap\left\{f : \  \|f\|_{L^2\left(w\otimes W_{[-R,R]}^{\otimes p}\right)}\le M \right\}$. The proof below uses techniques from the proof of Lemma 11 in \cite{PSWF}.\\ 
\noindent {\bf Proof.} In this proof,  $\braket{\cdot,\cdot} $ and $\left\| \cdot \right\|$ denote the scalar product and norm in $L^2([-R,R]^p)$. 
Take $f \in H^{q',s,\sigma'}(l)$. 
Let $N\geq 10$ and $ \widetilde{N}:= \lfloor \tau N \rfloor $, where $ \tau:= 1/(\pi e) $. 
We have 
\begin{align}
\left\| \varphi^t - P^N_{x_0t} \varphi^t \right\|^2&\leq  2 \left(\left\| \varphi^t - \mathcal{E}^{ \widetilde{N}} \varphi^t -  P^N_{x_0t} \left( \varphi^t - \mathcal{E}^{ \widetilde{N}} \varphi^t \right)  \right\|^2  +K^{N,\widetilde{N}}_{t}\right) \notag\\
&\leq  2 \left(\left\| \varphi^t - \mathcal{E}^{ \widetilde{N}} \varphi^t  \right\|^2+K^{N,\widetilde{N}}_{t} \right) . \label{eq:sob1}
\end{align}
Using that  $\left(\psi_{\mt{m}}^{x_0t}\left( \cdot/R\right)/R^{p/2} \right)_{\mt{m}\in\N_0^p}$ are orthonormal in $L^2([-R,R]^p)$  and the Cauchy-Schwarz inequality in the second display yield
\begin{align*}
K^{N,\widetilde{N}}_{t}
=& \left\|  \sum_{\mt{k} \in \Z^p:\  |\mt{k}|_{\infty} < \widetilde{N}}    \braket{   \varphi^t, \phi_{\mt{k}}\left(\frac{\cdot}{R}\right) }   \left( \sum_{\abs{\mt{m}}_{\infty}\ge N}
\braket{\phi_{\mt{k}}\left(\frac{\cdot}{R}\right) , \psi_{\mt{m}}^{x_0t}\left(\frac{\cdot}{R}\right)}  \psi_{\mt{m}}^{x_0t}\left(\frac{\star}{R}\right)\frac{1}{R^p}\right)     \right\|^2   \\
\leq& \sum_{\mt{k} \in \Z^p:\ |\mt{k}|_{\infty} < \widetilde{N}}  \abs{ \braket{  \varphi^t, \phi_{\mt{k}}\left(\frac{\cdot}{R}\right) }     }^2    I_{N, \widetilde{N}}(t)   \leq \left(\sum_{\mt{k} \in \Z^p}   \abs{ \mathcal{F}[f](t,\mt{k})  }^2\right)    I_{N, \widetilde{N}}(t).
\end{align*}
We have, using (18.17.19) in \cite{olver2010nist} for the first equality and for all $k\in \Z$ and $m\in \N_0$,
\begin{align*}
& \left|\braket{\phi_k\left(\frac{\cdot}{R}\right) , \psi_{m}^{x_0t}\left(\frac{\cdot}{R}\right)}\right|^2
= \frac{R}{2}\left| I_{m,k} + O_{m,k}\right|^2\leq R\left( |I_{m,k}|^2 + |O_{m,k}|^2\right) ,\\
&I_{m,k}:=  \sum_{n=0}^{\lfloor 5m/8\rfloor -1}  \beta_n^{m}(x_0t) \braket{e^{i\pi k\cdot}, P_{n }}_{L^2([-1,1])} ,\  
O_{m,k} := \sum_{n\geq \lfloor 5m/8\rfloor } \beta_n^{m}(x_0t) i^n \sqrt{\frac{2}{|k|}}\sqrt{n + \frac{1}{2}} J_{n+1/2}(|k| \pi).
\end{align*}
Using, for all $k\in \Z$, $ \abs{\braket{e^{i\pi k\cdot}, P_{n }}} \leq \sqrt{2} $, Proposition 3 in \cite{PSWF}, and integral test for convergence (indeed, by  (3.4) page 34 in \cite{Osipov}, for all $m \geq 2\vee (e^2x_0\abs{t})$, $2 \sqrt{\chi_m^{x_0t}}/\left(x_0\abs{t}\right)\ge 2e^2>1$),
we obtain, for all $m \geq 2\vee (e^2x_0\abs{t}) $, 
\begin{align*}
\abs{I_{m,k}} &
&\leq \sqrt{\frac{5}{2\pi}}  \int_{0}^{\left\lfloor 5m/8\right\rfloor} \left(\frac{ 2 \sqrt{\chi_m^{x_0t}}}{x_0\abs{t}}\right)^x dx \abs{ \mu_{m}^{x_0t}} 
\leq \frac{\sqrt{5/(2\pi)}  }{\ln\left(2 \sqrt{\chi_m^{x_0t}}/\left(x_0\abs{t}\right)\right)} \ \left(\frac{ 2 \sqrt{\chi_m^{x_0t}}}{x_0\abs{t}}\right)^{ \lfloor 5m/8\rfloor} \abs{ \mu_{m}^{x_0t}}.
\end{align*} 
Let $m \geq 2\vee (e^2x_0\abs{t})$. 
Using Lemma \ref{upper_bound} for the first inequality, we obtain 
\begin{align*}
\abs{I_{m,k}} & \leq \frac{\sqrt{5e^3}}{3} \frac{1}{\ln(2)+2}  \left(\frac{ 2\sqrt{ m(m+1)+x_0^2t^2}}{x_0\abs{t}}\right)^{5m/8} \left( \frac{e x_0\abs{t}}{4(m+3/2)}\right)^m  \\
&  \leq  \frac{\sqrt{5e^3}}{3}  \frac{1}{\ln(2)+2}\left(\frac{\sqrt{e^2+ 1/e^2}(m+1)}{2^{11/5}(m+3/2)}\right)^{5m/8}   \exp\left( - \frac{3m}{8}  \ln\left( \frac{m}{ex_0\abs{t}}\right)\right)
\\& 
\leq a\exp\left( - \frac{3m}{8} \ln\left( \frac{m}{ex_0\abs{t}}\right)\right). 
\end{align*}  
Using, for all $j> -1/2$, $x\in\R$, and $n\in\N_0$, $\abs{J_{j}(x)} \leq \abs{x}^{j}/\left(2^{j}\Gamma(j+1)\right)$ (see 9.1.20 in \cite{abramowitz1965handbook}), $\abs{\beta_n^m(x_0t)}\leq1$, and $\sqrt{n+1/2}< \Gamma(n+3/2)/n!$ (see (5.6.4) in \cite{olver2010nist}) for the first inequality and $m> 8/5$ and $n! \geq \left(n/e\right)^n \sqrt{2\pi n}$ 
for the third, we obtain, for all $k\in \Z$,
\begin{align*}
\abs{O_{m,k}}    \leq& \sum_{n\geq \lfloor  5m/8\rfloor } \frac{\sqrt{\pi}}{n! } \left(\frac{|k| \pi}{2}\right)^n  \\
\leq &\frac{\sqrt{\pi}}{\lfloor 5m/8\rfloor!}\left(\frac{|k| \pi}{2}\right)^{\left\lfloor 5m/8\right\rfloor}\exp\left(\frac{|k| \pi}{2}\right)  
\le\sqrt{\frac{5m}{16 } } \left(\frac{|k|\pi e}{2(5m/8-1)} \right)^{5m/8} \exp\left(\frac{|k|\pi}{2}\right).
\end{align*}
Using   $\left|\braket{\phi_k\left(\cdot/R\right) , \psi_{m}^{x_0t}\left(\cdot/R\right)}\right|^2\leq R$ for all $(k,m)\in \N_0^2$ for the first inequality, $\sum_{|\mt{m}|_{\infty}=j} 1\leq p(j+1)^{p-1} $ for the second, \eqref{eq:Nk} and the convexity of $x\mapsto x^p$ for the fourth inequality, we have, for all $t$ such that $N \geq e^2x_0\abs{t}$,
\begin{align*}
I_{N, \widetilde{N}}(t)  &\leq  R^{p-1} \sum_{\mt{k} \in \Z^p:\ |\mt{k}|_{\infty} < \widetilde{N}}   \sum_{j=N}^{\infty} \sum_{\abs{\mt{m}}_{\infty}=j}  \prod_{l=1}^p \left|\braket{\phi_{\mt{k}_l}\left(\frac{\cdot}{R}\right) , \psi_{j}^{x_0t}\left(\frac{\cdot}{R}\right)}\right|^2 \\
&\leq R^{2p-1} \sum_{\mt{k} \in \Z^p:\ |\mt{k}|_{\infty} < \widetilde{N}}   \sum_{j=N}^{\infty} p(j+1)^{p-1}  \prod_{l=1}^p  \left(|I_{j,\mt{k}_l}|^2 +|O_{j,\mt{k}_l}|^2\right) \\
&\leq pR^{2p-1}  \sum_{\mt{k} \in \Z^p:\ |\mt{k}|_{\infty} < \widetilde{N}}   \sum_{j=N}^{\infty} (j+1)^{p-1}   \left( a^2\left( \frac{ex_0\abs{t}}{j}\right)^{3j/4} +\frac{5je^{|\mt{k}_l|\pi}}{16}\left(\frac{|\mt{k}_l|\pi e }{2(5j/8-1)}\right)^{5j/4}\right)  \\
&\leq \frac{p(4R^2\tau N)^p }{2R} \sum_{j=N}^{\infty} j^{p-1}a^{2p}\left( \frac{ex_0\abs{t}}{j}\right)^{3pj/4} + \left(\frac{5e^{\widetilde{N}\pi}j^2}{16}\right)^p\left(\frac{\widetilde{N}\pi e }{2(5j/8-1)}\right)^{5pj/4}\frac{1}{j}.
\end{align*} 
Using $ \kappa(t) := -3 \ln\left(ex_0\abs{t}/N\right)/8$, $\kappa(t) \geq 3/8$ for  $N \ge 2\vee (e^2x_0\abs{t})$, and  $\sup_{j\geq 1} j^{p-1} e^{-p\kappa(t)j} = (1-1/p)^{p-1}/(\kappa(t)e)^{p-1}$ for the second inequality, we obtain, for all $N \geq e^2x_0\abs{t}$,
\begin{align*}
\sum_{j=N}^{\infty} j^{p-1}\left( \frac{ex_0\abs{t}}{j}\right)^{3pj/4} &\leq \frac{(1-1/p)^{p-1}}{(\kappa(t)e)^{p-1}}\int_{N}^{\infty}e^{-p \kappa(t)j}dj \leq   \frac{(p-1)^{p-1} }{(p\kappa(t))^{p}e^{p-1}}e^{-p \kappa(t)N }.
\end{align*} 
Using $1-8/(5N) \geq 1/5$, 
that   for $j\geq N$, $\widetilde{N}\pi e/(2(5j/8-1)) \leq 10\tau\pi e/21=10/21$, and $\sup_{j\geq 1} j^{2p-1}e^{- 5pj/8} = ((2p-1)8/(5pe))^{2p-1}$ for the first inequality, we obtain
\begin{align*}
\sum_{j=N}^{\infty}j^{2p-1}\left(\frac{\widetilde{N}\pi e }{2(5j/8-1)}\right)^{5pj/4} &\leq  \left(\frac{(2p-1)8}{5p e}\right)^{2p-1} \int_{N}^{\infty}e^{-5pj\ln(21/10)/8}dj \\
& \leq \left(\frac{(2p-1)8}{5p e}\right)^{2p-1} \frac{8}{5p\ln(21/10)}e^{-5pN\ln(21/10)/8}.
&\end{align*}
Using \eqref{eq:sob1} for the first display, using $ \sup_{t: \ \abs{t} \leq N/(e^2x_0)} I_{N, \widetilde{N}} (t)\leq  c N^pe^{- bN}$ (because $5\ln(21/10)/8- \tau \pi  <3/8 \leq \kappa(t)$), using $s\geq \sigma + p/2$, $f \in H^{q,s,\sigma}(l)$, and, for all $\abs{t} > N/(e^2x_0)$, 
$I_{N, \widetilde{N}}(t) \leq  R^{p} \sum_{\mt{k} \in \Z^p: \ |\mt{k}|_{\infty} < \widetilde{N}} \left\| \phi_{\mt{k}} \right\|^2_{L^2([-1,1]^p)}\leq \left(2\tau NR\right)^p$ for the second, we have
\begin{align*}
&\int_{\R} \left\|  \mathcal{F}_{1\mathrm{st}}\left[f\right](t,\cdot_2) - P^N_{x_0t}  \mathcal{F}_{1\mathrm{st}}\left[f\right](t,\cdot_2) \right\|^2   dt  \notag\\
& \leq 2 \left\|  \mathcal{F}_{1\mathrm{st}}\left[f\right] - \mathcal{E}^{ \widetilde{N}} \mathcal{F}_{1\mathrm{st}}\left[f\right] \right\|^2_{L^2(\R\times[-R,R]^p)}    + 2  \int_{-N/(e^2x_0)}^{N/(e^2x_0)}  \sum_{\mt{k} \in \Z^p}   \abs{\mathcal{F}_{1\mathrm{st}}\left[f\right](t, \mt{k}) }^2   dt \sup_{t:\  \abs{t} \leq N/(e^2x_0)} I_{N, \widetilde{N}} (t) \notag\\
&\quad + \frac{2\left(2\tau NR\right)^p}{1\vee (N/(e^2x_0))^{2s}}  \int_{\abs{t} > N/(e^2x_0)} \sum_{\mt{k} \in \Z^p}  \abs{\mathcal{F}_{1\mathrm{st}}\left[f\right](t, \mt{k}) }^2 (1\vee  t^{2s})  dt      \notag \\
&\leq\frac{4\pi l^2}{  (\tau N-1)^{2\sigma}}   + \frac{4\pi l^2cN^p}{e^{bN}} + \frac{4\pi l^2\left(2\tau Re^2 x_0\right)^p(e^2 x_0)^{2\sigma}}{ N^{2\sigma} }. 
\end{align*}
Using $\tau -1/10 >0$ and \eqref{elog} yield the result.\hfill $\square$

\bibliographystyle{abbrv}
\bibliography{Bibliography-MM-MC3}

\end{document}